\documentclass[10pt,a4paper]{article}
\usepackage[utf8]{inputenc}
\usepackage{psfrag}
\usepackage{amsmath}
\usepackage{subfig}
\usepackage{fleqn}
\usepackage{amssymb}
\usepackage{amsmath}
\usepackage{amsthm}

\usepackage[nonumberlist]{glossaries}
\usepackage{graphicx}
\usepackage{hyperref}

\usepackage{subfig}
\usepackage{pdfpages}
\usepackage[sc]{mathpazo}
\usepackage{pdfpages}
\usepackage{tikz}
\usetikzlibrary{graphs}
\usepackage{circuitikz}
\usepackage{multirow}
\usepackage{placeins}
\usepackage{subfig}
\usepackage{bm}
\usepackage{booktabs}
\usepackage{wrapfig}
\usepackage{multicol}
\usepackage{lineno}
\usepackage{url}
\usepackage{longtable}
\usepackage{xcolor}
\usepackage{colortbl}

\usepackage{array}
\usepackage{booktabs}
\usepackage{tabulary}
\usepackage{rotating}

\usepackage{natbib}
\usepackage{authblk} 

\usepackage{algorithmicx} 
\usepackage{algorithm}
\usepackage[noend]{algpseudocode}

\usepackage[normalem]{ulem}  

\title{A generalized Riemann problem-based compact reconstruction method for finite volume schemes}
\author[1]{Gino I. Montecinos}
\author[2]{Eleuterio F. Toro}
\author[3]{Lucas O. Müller}

\affil[1]{Department of Mathematical Engineering, Universidad de la Frontera, Temuco, Chile}
\affil[2]{Laboratory of Applied Mathematics, DICAM, University of Trento, Trento, Italy}
\affil[3]{Department of Mathematics, University of Trento, Trento, Italy}

\begin{document}

\maketitle

\begin{abstract}

We present a Generalized Riemann Problem-based reconstruction method (GRPrec) for high-order finite volume schemes applied to hyperbolic partial differential equations. The method constructs spatial polynomials using cell averages at the current time level and GRP solution data from the previous time level.   The resulting GRPrec  stencil is as compact as that of discontinuous Galerkin (DG) schemes but unlike DG,  our finite volume schemes obey  a generous CFL stability condition that is independent of the order of accuracy.  We assess the method's performance through test problems for smooth and discontinuous solutions of the linear advection equation and the Euler equations of gas dynamics in one space dimension.  Results are compared against exact solutions and against numerical results from well-known spatial reconstruction finite volume and  DG schemes,  with all methods implemented in the fully discrete ADER framework.  The performance of GRPrec is very promising, especially in terms of efficiency,  that is error against CPU cost.

\end{abstract}


\section{Introduction}

Modern finite volume schemes for hyperbolic partial differential equations achieve high-order accuracy by performing a reconstruction step at each time level (or stage in the case of semi-discrete schemes) \cite{Toro:2009a}, \cite{Mecalizzi:2025a}. Reconstruction methods compute a spatial polynomial from available cell averages at the current time level.  Notably, the stencil of the reconstruction method, and thus of the numerical scheme in general, grows with the desired order of accuracy. 

While spatial reconstruction methods are a necessary ingredient of high-order finite volume schemes, they pose several challenges. First, their implementation can become cumbersome, especially in three-dimensional unstructured grids \cite{Kaser:2006a,Kaser:2007a}. Furthermore, reconstruction operators will in general not preserve desirable properties of the underlying numerical scheme, such as well-balancing properties, i.e. the capacity of the numerical scheme to preserve certain steady state solutions of PDEs system being discretized. Restoring well-balancing for standard reconstruction operators can be rather technical \cite{Muller:2013b,PIMENTELGARCIA2023111869,CASTRODIAZ2013242}. Another aspect to consider is that of applications in which one-dimensional (1D) hyperbolic PDE systems are applied on networks. Such applications can regard  traffic flow \cite{Canic:2015a}, pipe networks \cite{Borsche2014b} or the human circulatory system \cite{Alastruey:2011a,Boileau:2015a,Mueller:2014a}. In many of these cases, the length of 1D domains can be highly heterogeneous, with small domains discretized by only a few computational cells. Achieving high-order accuracy with finite volume schemes in these contexts can be challenging and requires substantial adjustments \cite{Muller:2015a,Contarino:2016a}.
Moreover, a distinguishing feature of spatial reconstruction schemes is non-linearity; such property is an attempt  to circumvent Godunov's theorem \cite{Godunov:1959a},  with the aim of computing solutions without, or much reduced, spurious oscillations in the vicinity of large gradients.  In this paper however, we focus fundamentally on reconstruction schemes on compact stencils to obtain schemes of high-order of accuracy.  Still,  some preliminary ideas are put forward to produce a non-linear version of our reconstruction procedure GRPrec. 

A well-established approach  that offers high order of accuracy and compact stencils is that of discontinuous Galerkin (DG) schemes \cite{Reed:1973a,Cockburn:2001a}. In this case, the numerical methods evolve a polynomial approximation of the sought solution in time, without requiring  a spatial reconstruction step at each time level. The stencil for DG schemes of any order of accuracy and in any number of space dimensions reduces to the computational cell of interest and its neighbours, needed to compute numerical fluxes at the cell interfaces. However, in these schemes the allowed time step size depends on the order of accuracy, and thus not only on the mesh spacing and the wave speeds of the PDE system being discretized, as is the case in finite volume schemes.  In fact, the higher the order of accuracy of the scheme, the lower the maximum Courant number allowed for stability. Furthermore, DG schemes need to be equipped with conventional slope limiting strategies \cite{Cockburn:2001a} or a posteriori limiting procedures such as the Multidimensional Optimal Order Detection (MOOD) method \cite{Diot:2013a}, in order to circumvent Godunov's theorem \cite{Godunov:1959a}. 

These considerations reveal that compactness of the stencil and generous linear stability conditions have not been reconciled,  in the context of numerical methods for hyperbolic PDEs.  Our goal is to design a reconstruction method strategy with a compact stencil, similar to that of DG schemes. Furthermore, when applied to standard finite volume schemes, stability requirements should resemble those of finite volume schemes. While the methodology proposed in this work could be applied to any 
fully-discrete or semi-discrete finite volume scheme, it fits naturally in the context of finite volume schemes that reach high-order accuracy in space and time using the Generalized Riemann Problem GRP$_m$ as the building block \cite{Toro:2002a, Toro:2024a}.  For a numerical scheme of order $m+1$,   GRP$_m$  is an initial value problem that admits source terms in the equations and where the initial condition is given by polynomials of degree $m \ge 1$.  Often in the literature, the notation GRP is meant to be GRP$_1$ ($m=1$) and no source terms; in such case the corresponding schemes are at most second-order accurate.  For notational convenience however,  in this paper we shall often use  GRP to mean  GRP$_m$.

In our reconstruction approach, we enrich the information available at time level $t^n$ and construct spatial polynomials by combining information from GRP solutions from time level $t^{n-1}$ and cell averages from current time $t^n$.  We call our reconstruction operator GRPrec, for GRP-based reconstruction.  
Noteworthy, our reconstruction procedure uses information of past time levels, and is therefore related to the Time Reconstruction (TR) scheme for solving the GRP \cite{Dematte:2020a}.   However,  in that work, the resulting ADER-TR scheme had a Courant stability limit smaller than  unity that depended  on the order of accuracy.   For example,  for the linear advection equation solved to fifth-order of accuracy,  the reported Courant stability limit was around 0.45.   In the present case we do not observe a reduction in the CFL stability limit below unity,  in line with all standard ADER finite volume schemes, apart from ADER-TR \cite{Dematte:2020a}.

We will restrict our presentation to the class of methods proposed by Toro and collaborators, called ADER (Arbitrary DERivative Riemann problems) \cite{Toro:2001c, Millington:2001a, Titarev:2002a, Toro:2002a,Schwartzkopff:2002a, Dumbser:2008b,
Toro:2024a}, a fully discrete non-linear method characterized by an accurate computation of integrals for source and flux terms, obtained by solving GRPs.  In particular, we use the ADER-DET solver for the GRP \cite{Dumbser:2008a}. This method reduces the solution of the GRP to (i) numerical evolution of polynomial data via a space-time discontinuous Galerkin approximation and (ii) interaction of evolved data at time integration points via classical Riemann solvers for numerical flux evaluation. Often, the time-evolution step is called the predictor step.  See \cite{Toro:2024a} for a recent review of all methods available to solve the GRP.

Here we compare our results with those obtained using the fully discrete ADER-DG scheme proposed in \cite{Dumbser:2008b}. This scheme shares the space-time evolution step of the DET solver for the GRP  in the finite volume ADER-DET method \cite{Dumbser:2008a},  making a head-to-head comparison with ADER-DET reasonable.   Furthermore,  ADER-DG has similar stability constrains to those of standard semi-discrete DG schemes.
In order to assess the performance of our method in the case of discontinuous solutions for nonlinear problems, we present a non-linear version of GRPrec, called GRPrecNL, that is based on the class of WENO reconstruction operators proposed in \cite{Dumbser:2007a}. While GRPrecNL is only a starting point in the design of non-oscillatory schemes based on GRPrec, it allows us to apply the proposed methodology to standard test problems.  Future work will regard this aspect.  In particular we shall explore the MOOD approach \cite{Diot:2013a} and other methods that have been applied to DG schemes \cite{Wei:2024a}.



The rest of the paper is structured as follows. We introduce general GRP-based high-order finite volume schemes in Section \ref{sec:the-ADER-scheme}. In Section \ref{sec:GRP-based-rec} we present our GRP-based reconstruction procedure in its linear and non-linear versions. Section \ref{sec:numerical-results} follows with numerical results.  Conclusions are drawn in Section \ref{sec:conclusions}.

\section{High-order Finite Volume Methods}\label{sec:the-ADER-scheme}

In this work we deal with the discretization of $N \times N$  hyperbolic systems of balance laws in one space dimension (1D) of the following form
\begin{equation}
\partial_t \mathbf{Q} + \partial_x \mathbf{F}(\mathbf{Q}) = \mathbf{S}(\mathbf{Q})\;, \label{eq_hyperbolic_system}
\end{equation}
where $\mathbf{Q}(x,t)$ is the vector of unknowns, taking values from the space of admissible states $\mathcal{A} \subset \mathbb{R}^N$, $\mathbf{F}(\mathbf{Q}) : \mathcal{A} \rightarrow \mathbb{R}^N$ is the physical flux and $\mathbf{S}(\mathbf{Q})  : \mathcal{A} \rightarrow \mathbb{R}^N$ is the source term. 

The discretization of systems like \eqref{eq_hyperbolic_system} by means of high-order numerical methods has flourished over the last fifty years, with finite volume schemes and, more recently, discontinuous Galerkin schemes, playing a prominent role. Within this context, we will focus on fully discrete finite volume schemes based on the GRP solution, yet to be formally defined.

Without loss of generality, we define the computational cell $\Omega_i = [x_{i-\frac{1}{2}}, x_{i+\frac{1}{2}}]$, with $\Delta x = x_{i+\frac{1}{2}} - x_{i-\frac{1}{2}}$, and the time interval $I^n = [t^n, t^{n+1}]$, with time step $\Delta t^n = t^{n+1} - t^n$. The computational cell, together with the specified time interval, define the space-time control volume $V^n_i = \Omega_i \times I^n$. 

After these preliminary definitions, we proceed with the description of high-order finite volume schemes. These schemes  emanate from exact integration of \eqref{eq_hyperbolic_system} in the control volume
$V^n_i $, namely
\begin{equation}
\mathbf{Q}^{n+1}_i = \mathbf{Q}^n_i - \frac{\Delta t^n}{\Delta x} \left( \mathbf{F}^n_{i + \frac{1}{2}} - \mathbf{F}^n_{i - \frac{1}{2}}\right) + \Delta t^n \mathbf{S}^n_i\;. \label{eq_finite_volume}
\end{equation}
Practical numerical schemes seek approximations to the emerging integrals in
(\ref{eq_finite_volume}), that is
\begin{equation}
 \mathbf{Q}^n_i \approx \frac{1}{\Delta x} \int_{\Omega_i} \mathbf{Q}(x,t^n) dx\;,
\end{equation}
\begin{equation}
 \mathbf{F}^n_{i \pm \frac{1}{2}} \approx \frac{1}{\Delta t^n} \int_{I^n}  \mathbf{F}(\mathbf{Q}(x_{i \pm \frac{1}{2}},t)) dt\; 
 \label{eq:numfluxes}
\end{equation}
and
\begin{equation}   \label{eq:source}
 \mathbf{S}^n_i \approx \frac{1}{\Delta x  \Delta t^n} \int_{I^n} \int_{\Omega_i}  \mathbf{S}(\mathbf{Q}(x,t)) dx dt\;.
\end{equation}
Consequently
\begin{equation}
 \mathbf{Q}^{n+1}_i \approx \frac{1}{\Delta x} \int_{\Omega_i} \mathbf{Q}(x,t^{n+1}) dx\;.
\end{equation}

In order to achieve high-order \sout{of} accuracy in space and time,  finite volume schemes of the form \eqref{eq_finite_volume} need to approximate space, time, and time-space integrals with suitable accuracy. Generally, all finite volume approaches of interest here share the two following building blocks: (i) spatial reconstruction, (ii) high-order accurate approximation of integrals.

In order to compute numerical fluxes \eqref{eq:numfluxes}, an ADER scheme of order $m+1$ in space and time relies on the solution of the GRP, formally defined as
\begin{equation}                                                \label{eq:GRP}
\hbox{GRP$_m$}: \; \left\{
 \left.\begin{array}{ll}
 \hbox{PDEs:}  &   \partial_{t} {\bf Q} +  \partial_{x} {\bf F}({\bf Q}) = {\bf S}({\bf Q})  
 \;,  \hspace{2mm}  x \in \mathbb{R}\;, t>t^n \;,  \\[5pt]
 \hbox{IC:}  & {\bf Q}(x,t^n) = \left\{
       \begin{array}{lll}
            {\bf P}^n_{i}(x)   & \mbox{ if } & x < x_{i+\frac{1}{2}} \;, \\
            {\bf P}^n_{i+1}(x) & \mbox{ if } & x > x_{i+\frac{1}{2}}  \;.\\
       \end{array}\right.
\end{array}\right.\right.
\end{equation}
Here,  $\mathbf{P}^n_i(x)$ denotes a spatial polynomial of degree $m$  at time $t^n$,  constructed from cell averages $\mathbf{Q}^n_j$, with $j \in S_i = \{i-l, \dots, i, \ldots, i+r\}$, with $l,r \in \mathbb{N}$, and where $S_i$ is the stencil of a reconstruction method.   

Note that the solution of \eqref{eq:GRP}  at the fixed interface position $x=x_{i+\frac{1}{2}}$ depends on time and will be denoted as $\mathbf{Q}^n_{i+\frac{1}{2}}(\tau)$; in practice,  local coordinates are used, with new variable $\tau = t - t^n$, so that $\tau \in [0, \Delta t^n]$.  Clearly,  the suffix $i+\frac{1}{2}$ indicates that $\mathbf{Q}^n_{i+\frac{1}{2}}(\tau)$ is evaluated at cell interface $x_{i+\frac{1}{2}}$.  Evaluation of  the time integral (\ref{eq:numfluxes}) to the appropriate order gives the required numerical flux in \ref{eq_finite_volume}.  A similar but simpler approach is used to determine the numerical source from (\ref{eq:source}), thus completely determining the one-step scheme  \ref{eq_finite_volume}.\\

\noindent{\bf Remark:} The generalized Riemann problem  GRP$_m$ in (\ref{eq:GRP}) is a two-fold generalization of the classical, homogeneous  piece-wise constant data Riemann problem \cite{Toro:2002a}. That is (i) the initial conditions are polynomials of arbitrary degree $m$ and (ii) the equations admit source terms.   Regarding notation,  as already pointed out,  we shall  often use  GRP to mean GRP$_m$,  as defined in (\ref{eq:GRP}).

Different ADER schemes adopt various strategies to compute the high-order GRP solution $\mathbf{Q}^n_{i+\frac{1}{2}}(\tau)$.  The Toro-Titarev GRP solver \cite{Toro:2002a} is based on analytical evolution of the  solution by means of Taylor-type series expansions, use of the Cauchy-Kovalevskaya  procedure and solution of conventional linear Riemann problems for spatial derivatives.  In \cite{Castro:2008a}, the Harten,  Engquist, Osher and Chakravarthy  (HEOC) GRP solver was put forward by re-interpreting the method of Harten et al.   \cite{Harten:1987b} as an ADER method,  with a corresponding solver for the GRP.   In this case,  the limiting values of the spatial polynomials from the left and right neighbouring cells with interface $x_{i+\frac{1}{2}}$ are evolved in time locally using Taylor series expansions and the Cauchy-Kovalevskaya  procedure. Then, $\mathbf{Q}^n_{i+\frac{1}{2}}(\tau)$ is computed at any time integration point $\tau$ by solving classical non-linear Riemann problems.  DET is another successful solver for the GRP, due to  Dumbser, Enaux and Toro \cite{Dumbser:2008a}. This solver is a generalisation of HEOC,  in which  the spatial polynomials are evolved in time numerically  through a locally implicit discontinuous Galerkin method,  thus delivering a space-time polynomial defined in $V^n_i$. This polynomial, available in each computational cell,  is then used to compute $\mathbf{Q}^n_{i+\frac{1}{2}}(\tau)$  by solving classical non-linear  Riemann problems at time integration points; this is like the interaction step in HEOC.  Moreover, the space-time polynomial distribution in each cell is also used to compute source terms to the required accuracy.  Importantly,  the time-evolution step via a locally implicit method  results in a high-order ADER method that admits stiff source terms  \cite{Dumbser:2008a}.  

We will compare results from our methods with those obtained from the fully discrete discontinuous Galerkin scheme proposed in \cite{Dumbser:2008b},  since this scheme shares the time-evolution step (space-time predictor step) of our finite volume scheme.  For full details on available solvers for the high-order GRP see for example \cite{Montecinos:2012a,Toro:2024a},  and references therein.
Noteworthy,  apart from the TR approach \cite{Dematte:2020a}, all existing schemes to obtain $\mathbf{Q}^n_{i+\frac{1}{2}}(\tau)$ depend on  polynomials $\mathbf{P}^n_i(x)$ and $\mathbf{P}^n_{i+1}(x)$, which in turn are obtained by applying a reconstruction method to cell averages of the current time level. 

In the next section we propose a novel compact reconstruction method that, besides using current time level cell averages, makes use of previous time level GRP solutions at cell interfaces.

\section{A novel GRP-based reconstruction}\label{sec:GRP-based-rec}

In this section, two new reconstruction procedures are presented. The first one is linear, in the sense of Godunov \cite{Godunov:1959a}, and is called Generalized Riemann Problem-based reconstruction (GRPrec). The second one is called GRPrecNL, since nonlinearity is introduced in the reconstruction operator in attempting to circumvent Godunov's theorem.

\subsection{GRPrec reconstruction procedure}

The presentation of our novel reconstruction method is provided for polynomials of degrees one to four, which will then be used to construct one-step second to fifth-order accurate, in space and time, ADER finite volume schemes. Without loss of generality, we consider a scalar problem, so that, for time $t^n$, the reconstructed polynomial at the $i$-th cell is $p^n_i(x)$ and the GRP solution at interface $i+\frac{1}{2}$ is $q^n_{i+\frac{1}{2}}(\tau)$. Furthermore, we assume that at time $t^n$, the set of cell averages $\{q^n_j\}_{j=i-1}^{i+1}$, and GRP solutions $\{q^{n-1}_{j+\frac{1}{2}}(\Delta t^{n-1})\}_{j=i-1}^{i}$, are available. Notably, all polynomials of degree one to four will be computed using this compact stencil. 

Before describing our reconstruction procedure, we note that throughout this work, when dealing with systems of equations, reconstruction is performed in characteristic variables, as is customary in the literature.

\paragraph{First degree GRPrec polynomial.} A first degree polynomial, $p^n_i(x)$ has two degrees of freedom. A natural way to compute them would be to enforce
\begin{equation}
p^n_i(x_{i-\frac{1}{2}}) = q_{i-\frac{1}{2}}^{n-1}(\Delta t^{n-1})\;, \label{eq:interpleft}
\end{equation}
and 
\begin{equation}
p^n_i(x_{i+\frac{1}{2}}) = q_{i+\frac{1}{2}}^{n-1}(\Delta t^{n-1})\;. \label{eq:interpright}
\end{equation}
Naturally, these two conditions result in a polynomial that violates the conservation property
	\begin{eqnarray}
	\begin{array}{c}
		\displaystyle
		\frac{1}{\Delta x} \int_{x_{i-\frac{1}{2}}}^{x_{i+\frac{1}{2}}} p^n_i(x) dx = q_i^n\;,
	\end{array} \label{eq:conservation}
\end{eqnarray}
To overcome this, we employ a least squares approach and compute the desired degrees of freedom as those that minimize the functional 
\begin{eqnarray}
\begin{array}{c}
f(p_i,\mu) = 
\biggl(
p^n_i(x_{i-\frac{1}{2}}) - q_{i-\frac{1}{2}}^{n-1}(\Delta t^{n-1}) 
\biggr)^2
  +
\biggl(
p^n_i(x_{i+\frac{1}{2}}) - q_{i+\frac{1}{2}}^{n-1}(\Delta t^{n-1}) 
\biggr)^2
\\
+
 \mu \biggl( q_i^n - \frac{1}{\Delta x}\int_{x_{i-\frac{1}{2}}}^{x_{i+\frac{1}{2}}} p^n_i(x) dx 
 \biggr) \;,
\end{array}
\end{eqnarray}
in which $\mu$ is the Lagrange multiplier coefficient. Noteworthy, one could compute the degrees of freedom by simply imposing conservation property \eqref{eq:conservation} and choosing one of the two interpolatory conditions \eqref{eq:interpleft} and \eqref{eq:interpright}.  Fig. \ref{figure:sketch-rec-ord2} provides a schematic representation of the information used by GRPrec to construct a first degree polynomial.

\paragraph{Second degree GRPrec polynomial.} 

In order to define a second degree reconstruction polynomial, we impose the conservation condition \eqref{eq:conservation} and both interpolatory conditions  \eqref{eq:interpleft} and \eqref{eq:interpright}. Fig. \ref{figure:sketch-rec-ord3} provides a schematic representation of the information used by GRPrec to construct a second degree polynomial.

\paragraph{Third degree GRPrec polynomial.} 

In this case we proceed as for the first-degree case, but including two additional conservation properties, namely, we compute the polynomial coefficients as those that minimize the functional
\begin{eqnarray}
\begin{array}{c}
f(p_i,\mu) = 
\displaystyle
\biggl(
		\frac{1}{\Delta x} \int_{x_{i+\frac{1}{2}}}^{x_{i+\frac{3}{2}}} p^n_i(x) dx - q_{i+1}^n
\biggr)^2

+

\biggl(
\frac{1}{\Delta x} \int_{x_{i-\frac{3}{2}}}^{x_{i-\frac{1}{2}}} p^n_i(x) dx - q_{i-1}^n
\biggr)^2
\\
\displaystyle
+
\biggl(
p^n_i(x_{i-\frac{1}{2}}) - q_{i-\frac{1}{2}}^{n-1}(\Delta t^{n-1}) 
\biggr)^2
  +
\biggl(
p^n_i(x_{i+\frac{1}{2}}) - q_{i+\frac{1}{2}}^{n-1}(\Delta t^{n-1}) 
\biggr)^2

\\
\displaystyle
+
 \mu \biggl( q_i^n - \frac{1}{\Delta x}\int_{x_{i-\frac{1}{2}}}^{x_{i+\frac{1}{2}}} p^n_i(x) dx 
 \biggr)
 \;.
\end{array}
\end{eqnarray}
Also in this case one could have chosen to directly impose conservation in the $i$-th cell and in one of its neighbouring cells, as well as interpolatory conditions  \eqref{eq:interpleft} and \eqref{eq:interpright}. Fig. \ref{figure:sketch-rec-ord4} provides a schematic representation of the information used by GRPrec to construct a third degree polynomial.

\paragraph{Fourth degree GRPrec polynomial.} 

Here we adopt the approach used for the second degree case. Namely, we enforce the two interpolatory conditions  \eqref{eq:interpleft} and \eqref{eq:interpright}, as well as the conservation property for cells $i-1$, $i$ and $i+1$. Fig. \ref{figure:sketch-rec-ord5} provides a schematic representation of the information used by GRPrec to construct a fourth degree polynomial.

\subsection{Uniqueness of the GRPrec polynomials}

In this section, we are interested on the issue of uniqueness of the polynomial obtained from the combination of cell averages, $q_j^n$ with  $j\in S_i:=\{ i-1,i,i+1\}$ and GRP values at cell interfaces, $q_{k+\frac{1}{2}}^{n-1}(\Delta t^{n-1})$ with $k\in C_i:=\{i-2,i-1,i, i+1 \}$.  
For instance, let $ p_1(x) = a_0 + a_1x +a_2 x^2 $ and $ p_2(x) = b_0 + b_1x +b_2 x^2 \;,$ be polynomials satisfying the conservation  property for one cell average, let us say $q_{i_1}^n$ with $i_1\in S_i$,  and  the interpolation condition on two GRP values at cell interfaces,  let us say $ q_{i_2+\frac{1}{2}}^n(\Delta t)$ and $ q_{i_3+\frac{1}{2}}^n(\Delta t)$,  with $i_2, i_3\in C_i$. Then,  one can demonstrate that the two polynomials are necessarily identical.  Indeed, let us define 
$$
q(x) := p_1(x) - p_2(x) \;.
$$
This polynomial satisfies

$$
q( x_{{i_2}+\frac{1}{2}} ) = 0
$$
and 
$$
q( x_{{i_3}+\frac{1}{2} }) = 0\;.
$$
Thus, $x_{{i_3}\pm\frac{1}{2} }$ are two roots of $q(x)$ so this can be expressed as
\begin{eqnarray}
\label{eq:uniq-1}
\begin{array}{c}
q(x) = (x-x_{{i_2}+\frac{1}{2}} ) (x-x_{{i_3}+\frac{1}{2}} ) r(x) \;,
\end{array}
\end{eqnarray}
for some function $r(x)$.  

On the other hand, since the $q(x)$ is a continuous function, the primitive 
$$
\displaystyle 
Q(x) := \int_{x_{i_1 - \frac{1}{2}}}^x q(\xi) d\xi \;
$$
is differentiable and satisfies $Q(x_{i_1 - \frac{1}{2}}) = 0$ and  $Q(x_{i_1 + \frac{1}{2}})=0$.  Furthermore from the Mean Value Theorem, there exist some $\eta \in (x_{i_1 - \frac{1}{2}}, x_{i_1 + \frac{1}{2}} )$  such that $Q'(\eta) =0$; since $Q'(x)=q(x)$ we have $r(\eta) = 0$ and so $\eta $ is a root of $r(x)$, then we can write $ r(x) = (x - \eta) d$, with $d$ a constant value. Therefore, $q(x)$ has the form
$$
q(x) = (x-x_{{i_2}+\frac{1}{2}} ) (x-x_{{i_3}+\frac{1}{2}} ) (x-\eta) d \;.
$$
Since $q(x)$ is built from the sum of two polynomials of degree two, $q(x) = c_0 + c_1 x+ c_2 x^2$ with $c_j = a_j - b_j$ also has degree two. Therefore, the only possibility to have the three roots; $x_{{i_2}+\frac{1}{2}}$, $x_{{i_2}+\frac{1}{2}}$ and $\eta$,  is that $d = 0$, or $q(x) \equiv 0$.  Therefore, $c_j = 0 $ and so $a_j = b_j$ for $j=0,1,2$, which means $p_1(x) = p_2(x)$.

The same procedure can be adopted to demonstrate, in general, the uniqueness of polynomials built from combinations of elements of $S_i$ and $C_i$.

In the present approach, the information available around the cell $\Omega_i = [x_{i-\frac{1}{2}}, x_{i+\frac{1}{2}}]$, consists of three cell averages and four GRPs at cell interfaces. Therefore,  this would allow us to build a polynomial of up to sixth degree,  which could generate a seventh order scheme.  However, the construction of a polynomial of a given degree $m$ does not necessarily guarantee that the derived numerical scheme would have order of accuracy $m+1$ in space and time.

\subsection{A non-linear GRPrec reconstruction procedure}

We present a preliminary attempt to construct a non-linear version of the reconstruction procedure introduced in the previous section, hereafter called GRPrecNL.   We consider three polynomials $p^n_{i,L}(x)$, $p^n_{i,C}(x)$ and $p^n_{i,R}(x)$, where $p^n_{i,C}(x)$ corresponds to the linear GRPrec reconstruction presented in the previous section.  In order to obtain the sought non-linear polynomial we use a convex combination of the three aforementioned polynomials, following the approach devised in \cite{Dumbser:2007a}.  Briefly, we look for
\begin{eqnarray}
	\begin{array}{c}
	p^n_i(x)= p^n_{i,L}(x) \beta_L + p^n_{i,C}(x) \beta_C +  p^n_{i,R}(x) \beta_R \;.
	\end{array}
\end{eqnarray}
Here the weights are computed as
\begin{eqnarray}
	\begin{array}{c}
		\beta_L = \frac{\omega_L}{ \bar{\omega} } \;,
		\beta_C = \frac{\omega_C}{ \bar{\omega}  } \;,
		\beta_R = \frac{\omega_R}{ \bar{\omega}  } \;, 
	\end{array}
\end{eqnarray}
with 
\begin{eqnarray}   \label{eq:omega}
\begin{array}{c}
\omega_k = \frac{\lambda_k}{ (\varepsilon + OI_k)^r   } \;,  \hspace{2mm}  k \in \{L,C,R\} \;,  \hspace{2mm}\bar{\omega} = \omega_L + \omega_C + \omega_R\;.
\end{array}
\end{eqnarray}
We take $r =4$ and $\varepsilon=10^{-14}$ in (\ref{eq:omega});  $OI_k$ denotes the oscillation indicator
\begin{eqnarray}
\begin{array}{c}
\displaystyle
OI_k = \sum_{l=1}^m \int_{x_{i-\frac{1}{2}}}^{x_{i+\frac{1}{2}}} \biggl( \frac{d^l p^n_{i,k}}{dx^l}(x)  \biggr)^2 \Delta x^{2l} dx\;,\;\; k\in \{L,C,R\} \;.
\end{array}
\end{eqnarray}
Additionally,  we set $\lambda_L=\lambda_R=1$ and $\lambda_C = 10^9$, following  \cite{Dumbser:2007a}.

A pending issue to determine $p^n_i(x)$,  to be addressed in the next section,  concerns the  definitions for $p^n_{i,L}(x)$ and $p^n_{i,R}(x)$.  In all cases, it is assumed that at time $t^n$, the set of cell averages $\{q^n_j\}_{j=i-1}^{i+2}$, and GRP solutions $\{q^{n-1}_{j+\frac{1}{2}}(\Delta t^{n-1})\}_{j=i-2}^{i+1}$, are available.  Here $q^n_j$ denotes a component (scalar) of the vector 
${\bf Q}^{n}_j$ of conserved variables.

\subsection{Polynomials for schemes of second to fifth order}

Here we specifically construct non-linear polynomials of degree $m$, with $m=1,2,3,4$, giving rise to finite volume schemes of $2nd$ to $5th$ accuracy in space and time.

\paragraph{First degree GRPrecNL polynomial.} 

\begin{figure}
\centering
\subfloat[][Left stencil    ]{ \includegraphics[scale=0.5]{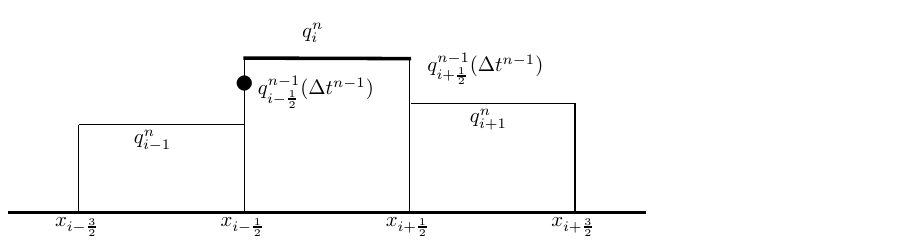}}
\subfloat[][Right stencil   ]{ \includegraphics[scale=0.5]{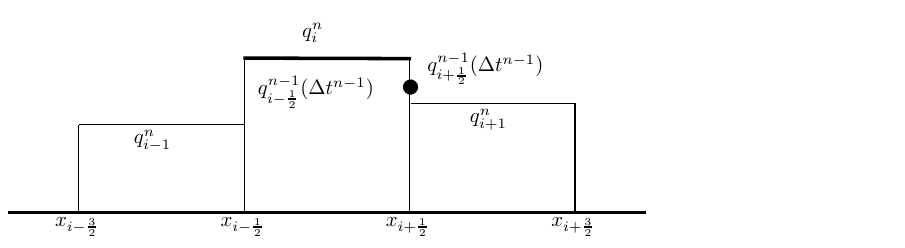}}
\quad
\subfloat[][Centred stencil]{ \includegraphics[scale=0.5]{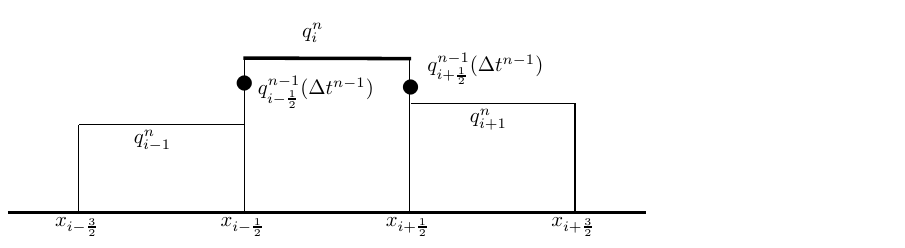}} 
\caption{Sketch for first degree GRPrec and GRPrecNL: Stencils for left, centred and right polynomials. Circles represent cell-interface point values of GRP solutions at the previous time level, full lines indicate the use of cell averages for the conservation property.}
\label{figure:sketch-rec-ord2}
\end{figure}

We compute coefficients for  $p^n_{i,L}(x)$ by enforcing conservation on the $i$-th cell
\begin{equation}
  \displaystyle
\frac{1}{\Delta x}  \int_{ x_{i-\frac{1}{2}}}^{x_{i+\frac{1}{2}}}
p^n_{i,L}( x)dx = q_{i}^n  \;
\end{equation}
and the interpolatory condition \eqref{eq:interpleft}, on interface $i-\frac{1}{2}$, namely
\begin{equation}
p^n_{i,L}(x_{i-\frac{1}{2}}) = q_{i-\frac{1}{2}}^{n-1}(\Delta t^{n-1})\;.
\end{equation}
Similarly, in order to define $p^n_{i,R}(x)$, we again enforce  
conservation on the $i$-th cell
\begin{equation}
  \displaystyle
\frac{1}{\Delta x}  \int_{ x_{i-\frac{1}{2}}}^{x_{i+\frac{1}{2}}}
p^n_{i,R}( x)dx = q_{i}^n  \;
\end{equation}
and the interpolatory condition \eqref{eq:interpright}, on interface $i+\frac{1}{2}$, which is
\begin{equation}
p^n_{i,R}(x_{i+\frac{1}{2}}) = q_{i+\frac{1}{2}}^{n-1}(\Delta t^{n-1})\;.
\end{equation}
Fig. \ref{figure:sketch-rec-ord2} shows a sketch for stencils associated with the left, central and right polynomials.
 
\begin{figure}
\centering
\subfloat[][Left stencil    ]{ \includegraphics[scale=0.5]{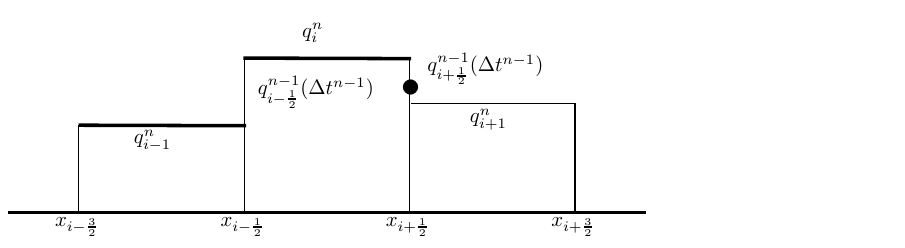}}
\subfloat[][Right stencil   ]{ \includegraphics[scale=0.5]{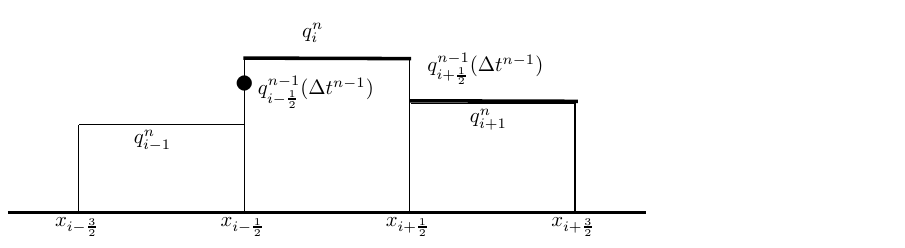}}
\quad
\subfloat[][Centred stencil]{ \includegraphics[scale=0.5]{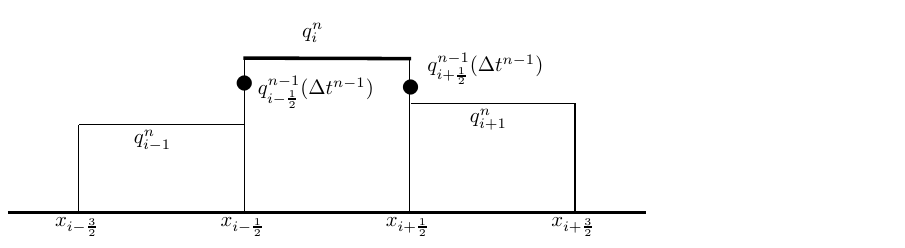}} 
\caption{Sketch for second degree GRPrec and GRPrecNL: Stencils for left, centred and right polynomials. Circles represent cell-interface point values of GRP solutions at the previous time level, full lines indicate the use of cell averages for the conservation property.}
\label{figure:sketch-rec-ord3}
\end{figure}

\paragraph{Second degree GRPrecNL polynomial.} 
In this case, coefficients for $p^n_{i,L}(x)$ are defined by enforcing the conservation property in cells $i$ and $i-1$, as well as the interpolatory condition \eqref{eq:interpright}, on interface $i+\frac{1}{2}$. Similarly, $p^n_{i,R}(x)$ is constructed by requiring the conservation property to be fulfilled on cells $i$ and $i+1$, as well as the interpolatory condition \eqref{eq:interpleft},  at interface $i-\frac{1}{2}$. Fig. \ref{figure:sketch-rec-ord3} shows a sketch for stencils associated with the left, central and right polynomials.

\begin{figure}
\centering
\subfloat[][Left stencil    ]{ \includegraphics[scale=0.5]{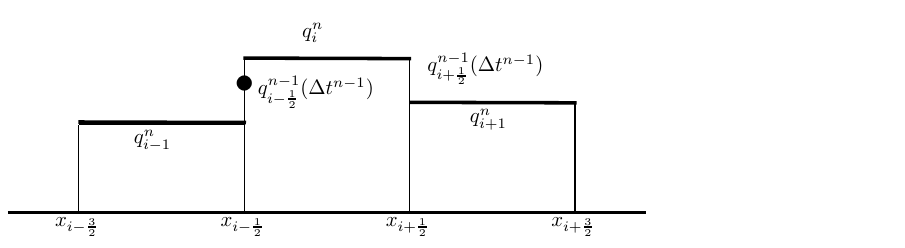}}
\subfloat[][Right stencil   ]{ \includegraphics[scale=0.5]{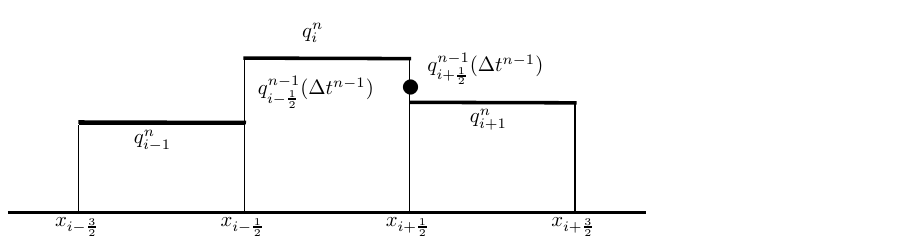}}
\quad
\subfloat[][Centred stencil]{ \includegraphics[scale=0.5]{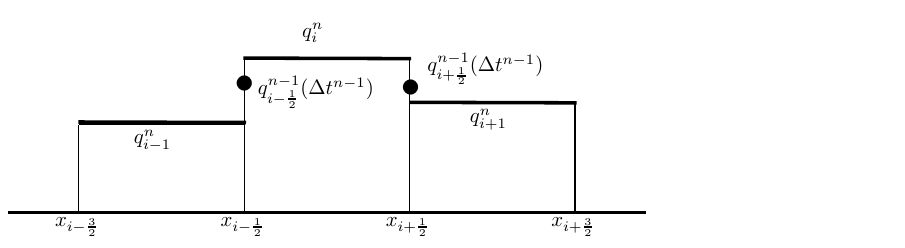}} 
\caption{Sketch for third degree GRPrec and GRPrecNL: Stencils for left, centred and right polynomials. Circles represent cell-interface point values of GRP solutions at the previous time level, full lines indicate the use of cell averages for the conservation property. }
\label{figure:sketch-rec-ord4}
\end{figure}

\paragraph{Third degree GRPrecNL polynomial.} 
In this case we compute coefficients for $p^n_{i,L}(x)$ by imposing conservation in cells $i-1$, $i$ and $i+1$, as well as the interpolatory condition \eqref{eq:interpleft} at interface $i-\frac{1}{2}$. In turn, for coefficients of $p^n_{i,R}(x)$, we use conservation in cells $i-1$, $i$ and $i+1$, as well as the interpolatory condition \eqref{eq:interpright}  at interface $i+\frac{1}{2}$. Fig. \ref{figure:sketch-rec-ord4} shows a sketch for stencils associated with the left, central and right polynomials.

\begin{figure}
\centering
\subfloat[][Left stencil    ]{ \includegraphics[scale=0.5]{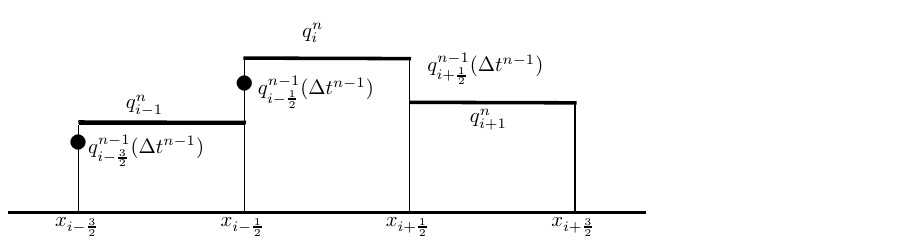}}
\subfloat[][Right stencil   ]{ \includegraphics[scale=0.5]{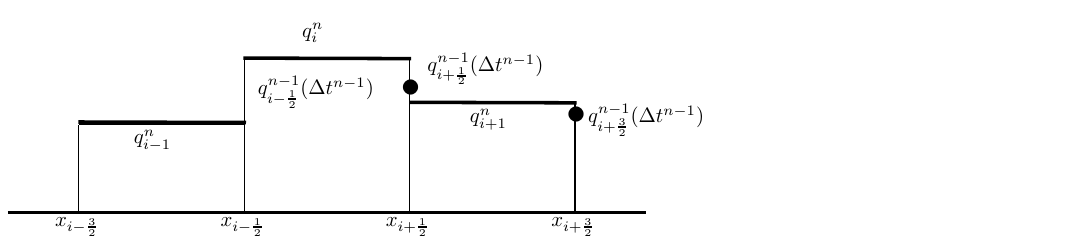}}
\quad
\subfloat[][Centred stencil]{ \includegraphics[scale=0.5]{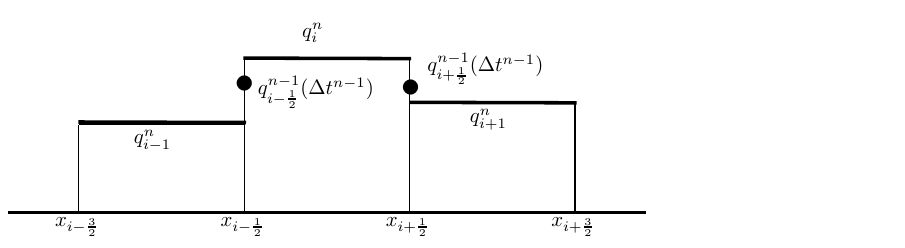}} 
\caption{Sketch for fourth degree GRPrec and GRPrecNL: Stencils for left, centred and right polynomials. Circles represent cell-interface point values of GRP solutions at the previous time level, full lines indicate the use of cell averages for the conservation property.}
\label{figure:sketch-rec-ord5}
\end{figure}

\paragraph{Fourth degree GRPrecNL polynomial.} 
For the fifth degree polynomial $p^n_{i,L}(x)$ we use conservation in cells $i-1$, $i$ and $i+1$, as well as the interpolatory condition \eqref{eq:interpleft} at interfaces $i-\frac{1}{2}$ and $i-\frac{3}{2}$. Similarly, for $p^n_{i,R}(x)$ we enforce conservation in cells $i-1$, $i$ and $i+1$, as well as the interpolatory condition \eqref{eq:interpright} at interfaces $i+\frac{1}{2}$ and $i+\frac{3}{2}$. Fig. \ref{figure:sketch-rec-ord5} shows a sketch for stencils associated with the left, central and right polynomials.

\subsection{Additional considerations on GRPrec}

It is worth noting that left and right stencils used in GRPrecNL could be used as linear reconstruction stencils. While any of these choices would deliver accurate reconstruction methods, the central stencil is always the most compact one. In fact, for first and second degree polynomials, this stencil uses information from the current cell and its interfaces. 

Another interesting aspect of the proposed reconstruction method is that, up to second degree polynomials, it will deliver well-balanced reconstructed polynomials if the GRP solver used to compute $q^n_{i+\frac{1}{2}}(\tau)$ is well-balanced.  Here we refer to a scheme as well-balanced if it can preserve some steady state solutions of the underlying model.  The construction of well-balanced versions of novel schemes of this paper is a pending issue and an in-depth discussion of this property is left for future work; nonetheless we  stress that this well-balanced property is a particularly valuable feature in the context of 1D blood flow in complex blood vessel networks \cite{Muller:2013b,Muller:2015a,PIMENTELGARCIA2023111869},  a topic of special interest to the authors.

We conclude this section by noting that a potential drawback of GRPrec is the need for GRP intermediate states,  which for methods different from ours might not be readily available.   For our methods however,  these are readily available; they provide interface knots and are used to compute polynomials coefficients.   

In the tests presented in the next section for the Euler equations, we show that simple approximate solvers can be used to compute such states.

\section{Numerical Results}\label{sec:numerical-results}

This section is devoted to the assessment of the numerical methods presented in this paper,  as compared to exact solutions and numerical solutions from other similar numerical schemes in the literature.

\subsection{Numerical schemes}

Table \ref{tab:reconstruction_summary} summarises all the numerical methods implemented in this paper, which include three finite volume methods (first 3 rows) and one fully-discrete discontinuous Galerkin method (fourth row).
All finite volume methods have Courant linear stability limit $1$; the fully discrete  DG scheme has the same stability condition as the semi-discrete version,  and depends on the order of accuracy $m+1$, or alternatively the degree $m$ of the underlying polynomial.   For example,  for a second-order of accurate DG scheme  ($m=1$),  the stability limit is $1/3$ and for a fith-order method the stability limit is $1/9$.
All methods use the DET solver for the generalized Riemann problem \cite{Dumbser:2008a},  which automatically renders the schemes suitable for dealing with stiff source terms, though this capability is not demonstrated here.
\begin{table}[ht]
\centering
\begin{tabular}{|l|l|l|l|}
\hline Method notation  & Reconstruction type & GRP$_{m}$ solver & Stability limit \\
\hline
FV+GRPrec & GRPrec & DET  & $1$ \\
\hline
FV+GRPrecNL & GRPrecNL & DET &   $1$ \\
\hline
FV+WENO-DK & WENO-DK & DET &  $1$ \\
\hline
DG & No reconstruction & DET &   $ \frac{1}{2m+ 1}$ \\
\hline
\end{tabular}
\caption{ADER fully discrete methods in the finite volume (first 3 rows) and discontinuous Galerkin finite elements (fourth row) frameworks considered in this paper.  The first two methods are new.  All schemes use the solution of the generalized Riemann problem GRP$_{m}$ (also denoted as GRP in the paper) as the building block; the approximate solution of GRP$_{m}$ is obtained here through the DET solver.  All schemes have Courant-number stability limit 1, except for the DG scheme,  for which the stability limit depends on the degree $m$ of the underlying polynomial for the scheme of accuracy $m+1$.}
\label{tab:reconstruction_summary}
\end{table}
%


We recall that the DET solver for the GRP includes a numerical time-evolution step within each cell (the so-called predictor step) and an evolved-data interaction step at the cell interface at integration points, through classical non-linear Riemann problem solutions; this second step is used for evaluating the DET numerical flux.  In this paper, this interaction step is also used for providing reconstruction data at the interface through 
\begin{itemize}
\item the exact solution of classical Riemann problems,  for calculations concerned with the linear advection equation,
\item the two-rarefaction approximatation of the classical Riemann problem,  for applications concerned with the Euler equations.
\end{itemize}

The time step size for finite volume schemes  (\ref{eq_finite_volume}) (FV+GRPrec/GRPrecNL/WENO-DK) is computed using the  Courant stability condition
\begin{eqnarray}
\label{dt-FV}
    \begin{array}{c}
\Delta t_{FV} = C_{CFL}\frac{\Delta x }{S_\mathrm{max}}\;,
    \end{array}
\end{eqnarray}
where $S_\mathrm{max}$ is an estimate of the maximum wave speed at time $t^n$. The time step for DG  is obtained according to \cite{Dumbser:2008b}, namely
\begin{eqnarray}
\label{dt-DG}
    \begin{array}{c}
\Delta t_{DG} = \frac{C_{CFL}}{2m+1} \frac{\Delta x }{ S_\mathrm{max} }\;,
    \end{array}
\end{eqnarray}
where $m$ is the degree of the  underlying polynomial for the DG scheme of accuracy $m+1$.  \\

\noindent{\bf Remark: Courant number of the computation. }  It is noted that for the finite volume schemes  the
Courant number of the computations is $C_{CFL}$, while for the DG schemes is $\frac{C_{CFL}}{2m+1}$.
Here unless otherwise stated,  the $C_{CFL}$ safety coefficient is set to $0.9$.

As to estimates for $S_\mathrm{max}$, strictly speaking, these must constitute bounds for the true maximal wave speed,  as discussed in \cite{toroBoundsWaveSpeeds2020}.   Here however,  we use  the following simple estimates:
\begin{itemize}
\item for the linear advection equation we set $S_\mathrm{max}=|\lambda|$, with $\lambda$ the characteristic speed of the PDE (exact),
\item for the Euler equations we set $S_\mathrm{max}= \max_{i , j} |\lambda_j(\mathbf{Q}^n_i)|$, with $i=1,\ldots,M$ and $j=1,\ldots,N$, with $N$ the number of eigenvalues of the non-linear system under consideration.
\end{itemize}

\subsection{The linear advection equation}

In this section, we consider the  linear advection equation
\begin{eqnarray}
\begin{array}{c}
\partial_t q + \lambda \partial_x q = 0\;, 
\end{array}
\end{eqnarray}
with $x \in [-1,1]$,  $\lambda=1$ and periodic boundary conditions.  We consider initial conditions for two test problems,  a solution featuring a complex profile,  evolved for long simulation times,  a square wave test, and an empirical convergence rates test for a smooth solution. 

\subsubsection{Multi-wave test}

We solve the classical multi-wave test of Jiang and Shu \cite{Jiang:1995a}, with the initial condition given by
\begin{eqnarray}
\begin{array}{c}
q(x,0) = \left\{
\begin{array}{cc}
\exp\biggl(-\ln(2)\frac{ (x+0.7)^2 }{0.0009}\biggr) &, -0.8 \leq x \leq -0.6 \\
1                                &, -0.4 \leq x \leq -0.2 \\
1-|10x-1|                        &, 0 \leq x \leq 0.2 \\
\biggl( 1- 100(x-\frac{1}{2})^2 \biggr)^\frac{1}{2}                      &, 0.4 \leq x \leq 0.6 \\
0                                & otherwise \;.
\end{array}

\right.
\end{array}
\label{eq:MW}
\end{eqnarray}
This test is employed to evaluate the performance of schemes for long simulation times,  for which diffusion and dispersion errors propagated in time and space become more evident,  thus revealing more clearly the limitations/strengths  of the methods under study.

\begin{figure}
\centering
\subfloat[][2nd order]{ \includegraphics[scale=0.22]{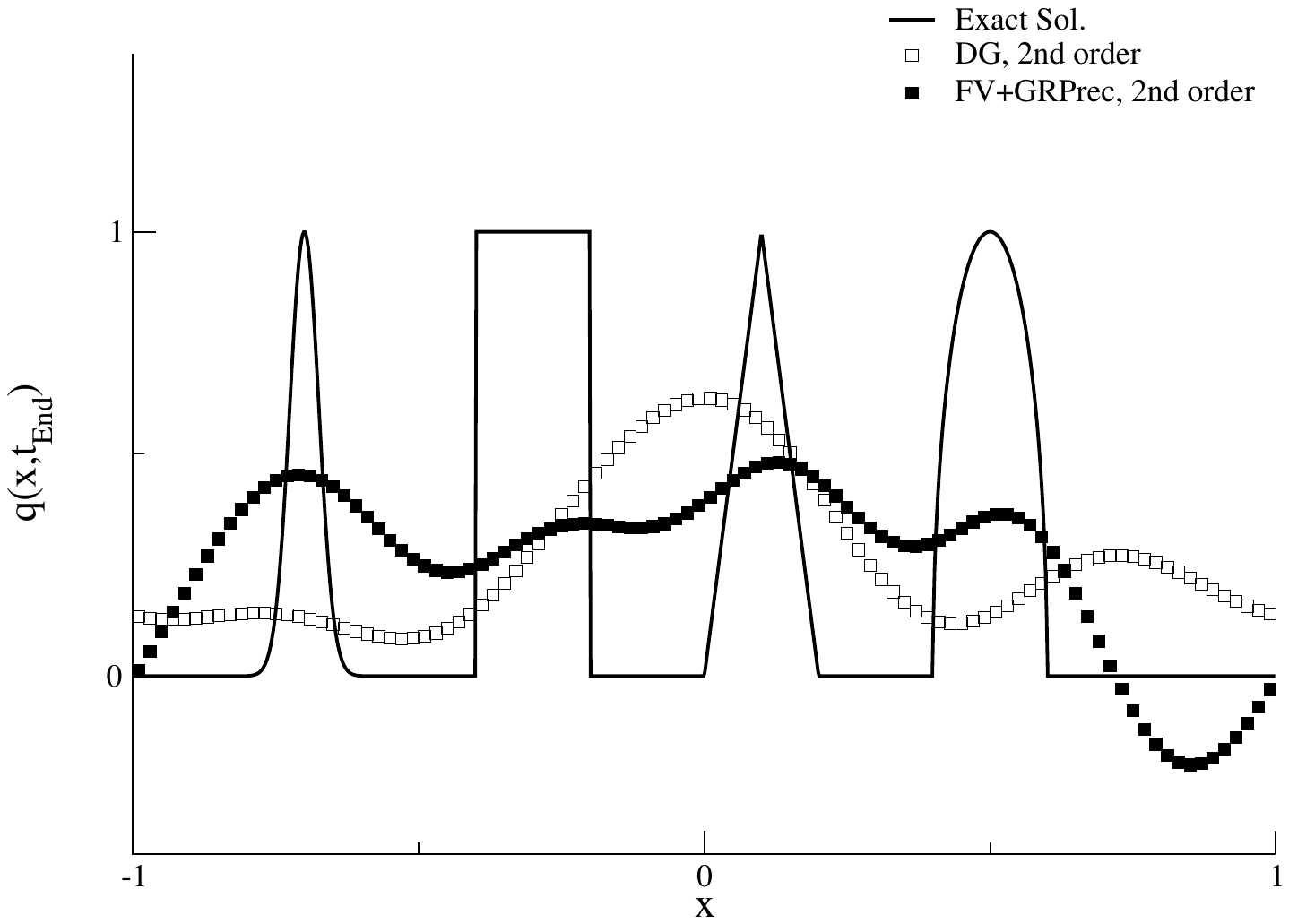}}
\subfloat[][3rd order]{ \includegraphics[scale=0.22]{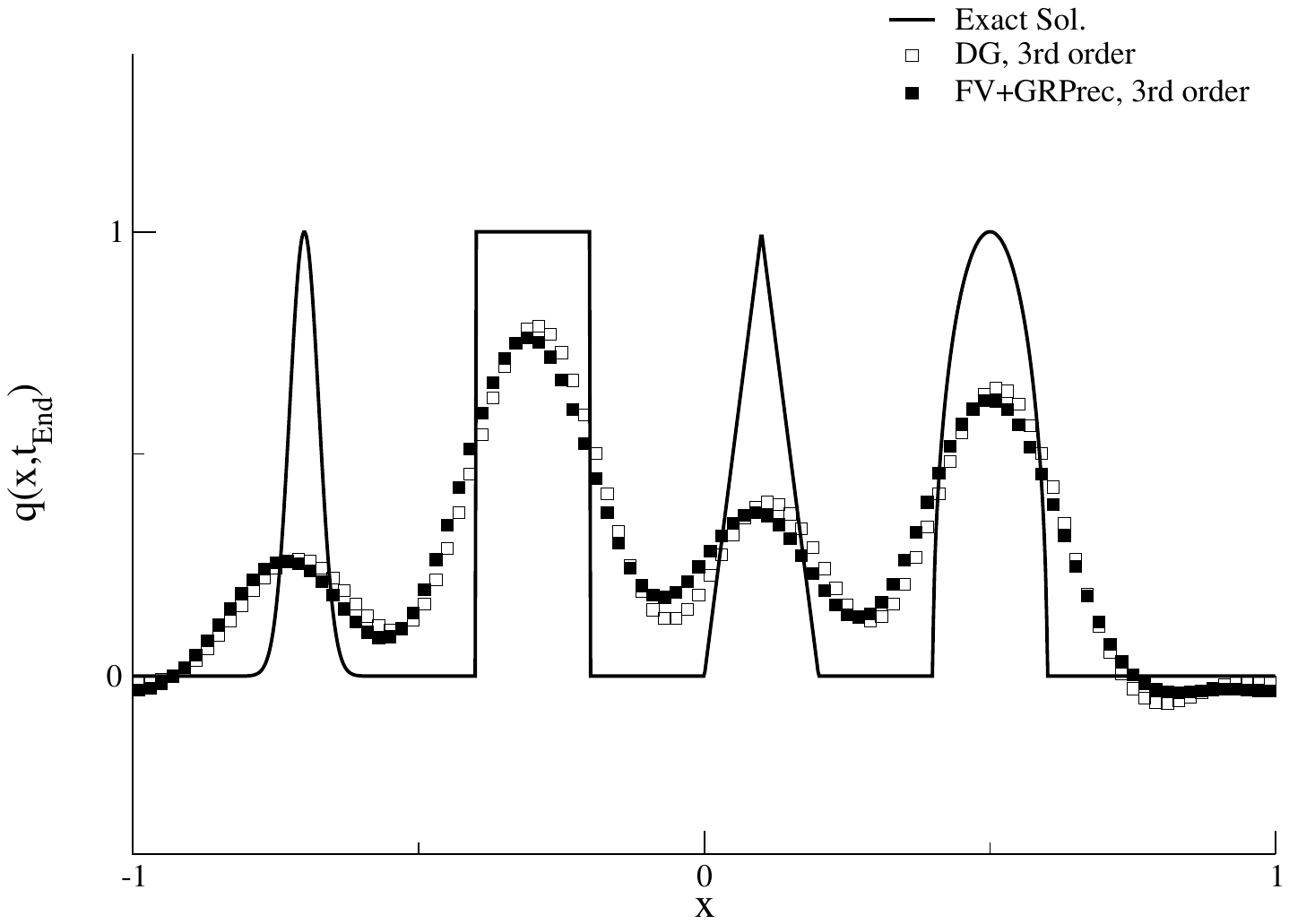}} 
\quad
\subfloat[][4th order]{\includegraphics[scale=0.22]{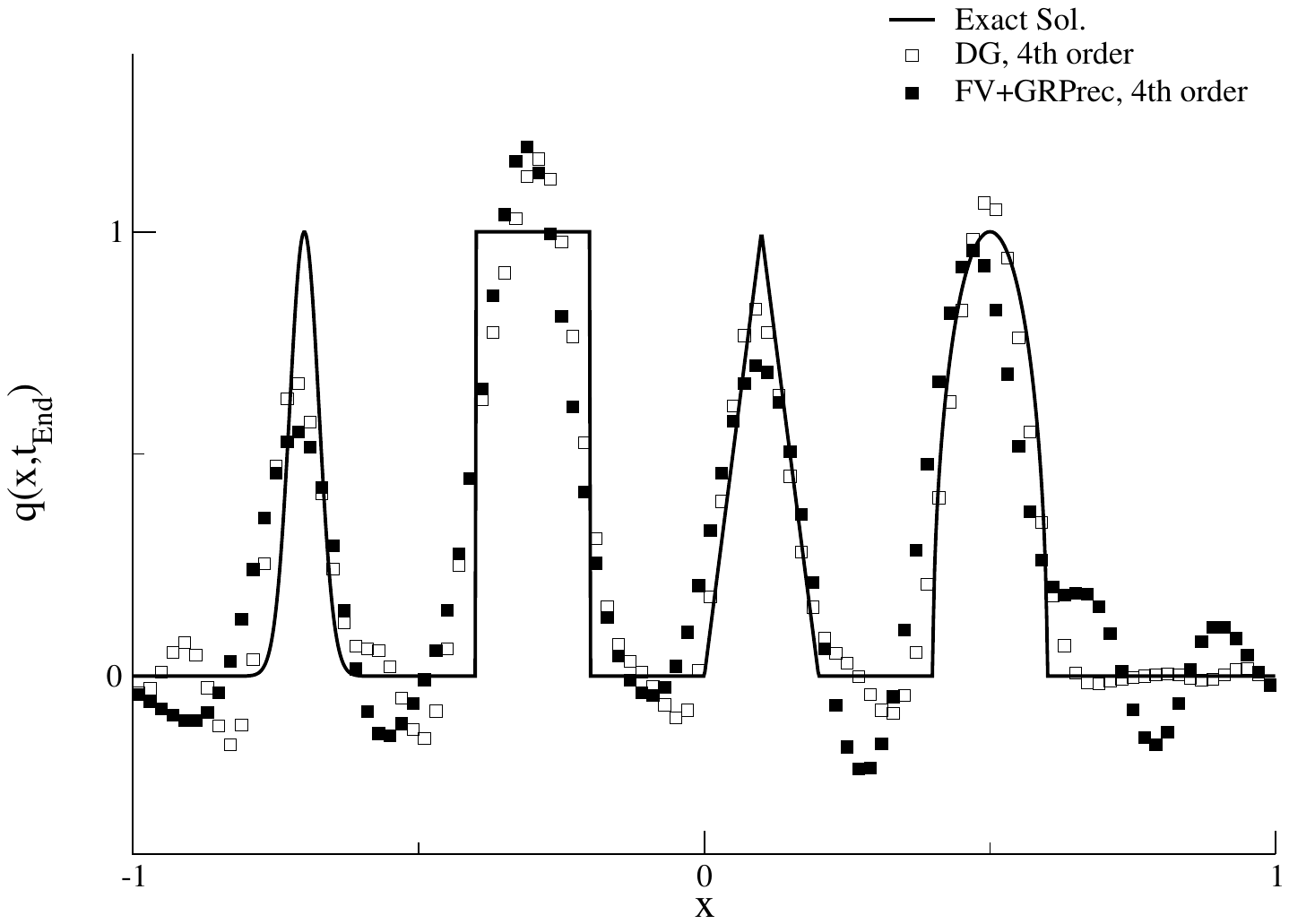}}
\subfloat[][5th order]{ \includegraphics[scale=0.22]{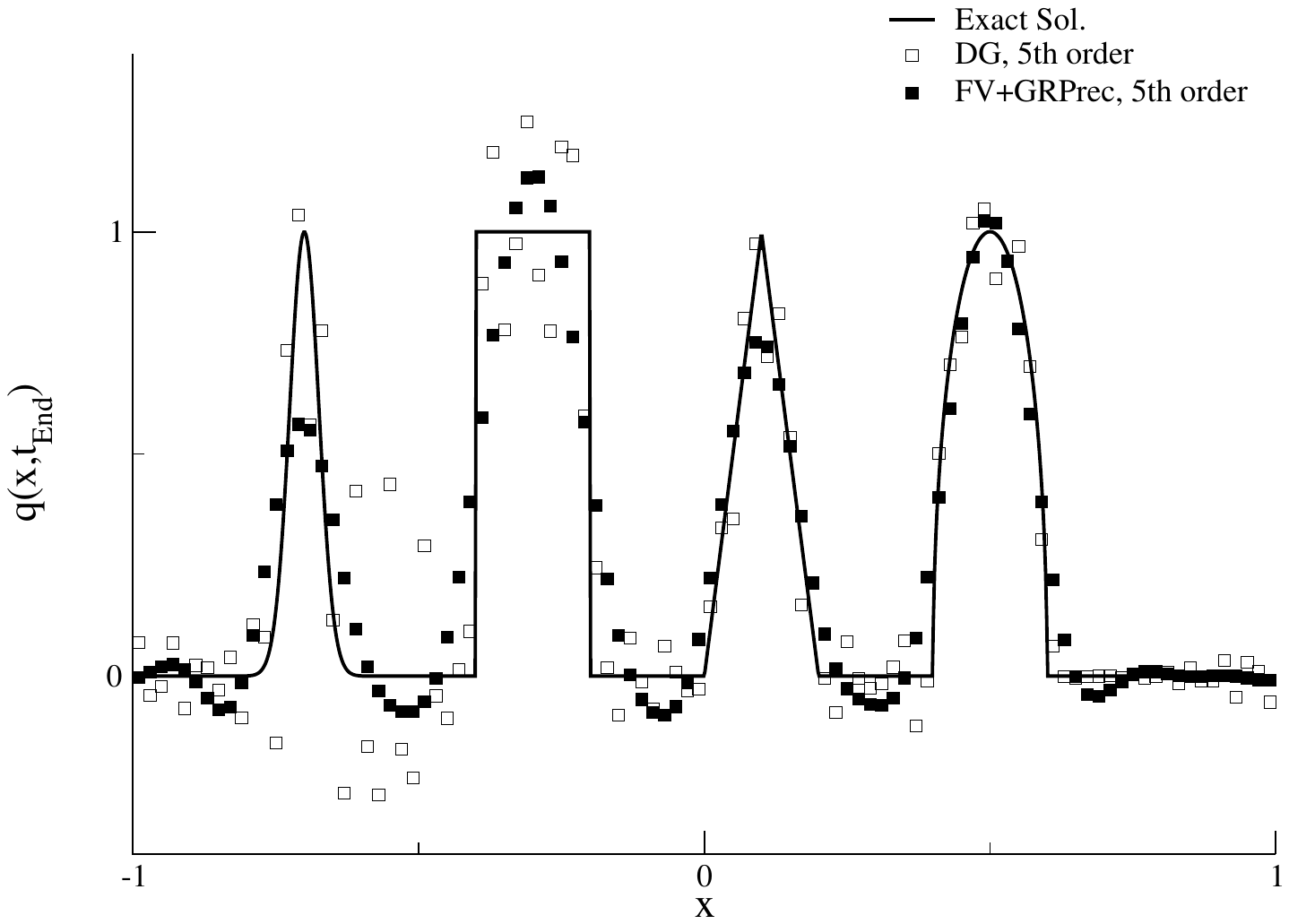}} 

\caption{Linear advection, multiwave-wave test with initial condition (\ref{eq:MW}): profiles for second, third, fourth with and fifth order of accuracy for DG and FV+GRPrec. Parameters: $t_{End} = 2000$,  $100$ cells and $C_{CFL}=0.9$.}
\label{figure:LA-multiwave-Linear}
\end{figure}

Fig. \ref{figure:LA-multiwave-Linear} shows results for DG and FV+GRPrec. We observe that results obtained with both linear schemes are similar, providing empirical evidence on the fact that the proposed reconstruction method results in a stable scheme, with similar performance to DG and allowing for large, finite volume-like CFL conditions.

%
%
 
%

%
%

\subsubsection{The square wave test}

\begin{figure}
	\centering
	\subfloat[][2nd order]{ \includegraphics[scale=0.22]{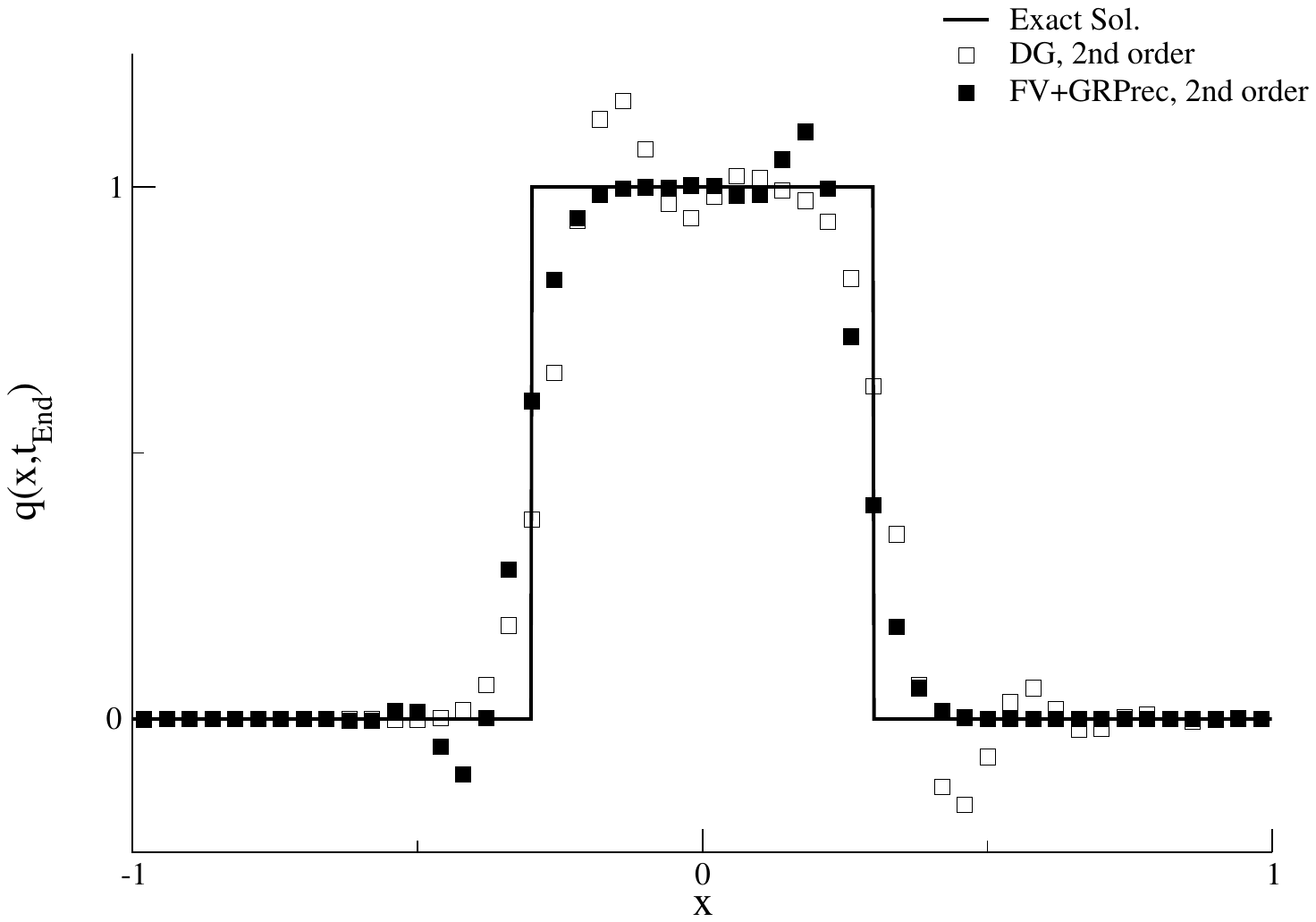}}
	\subfloat[][3rd order]{ \includegraphics[scale=0.22]{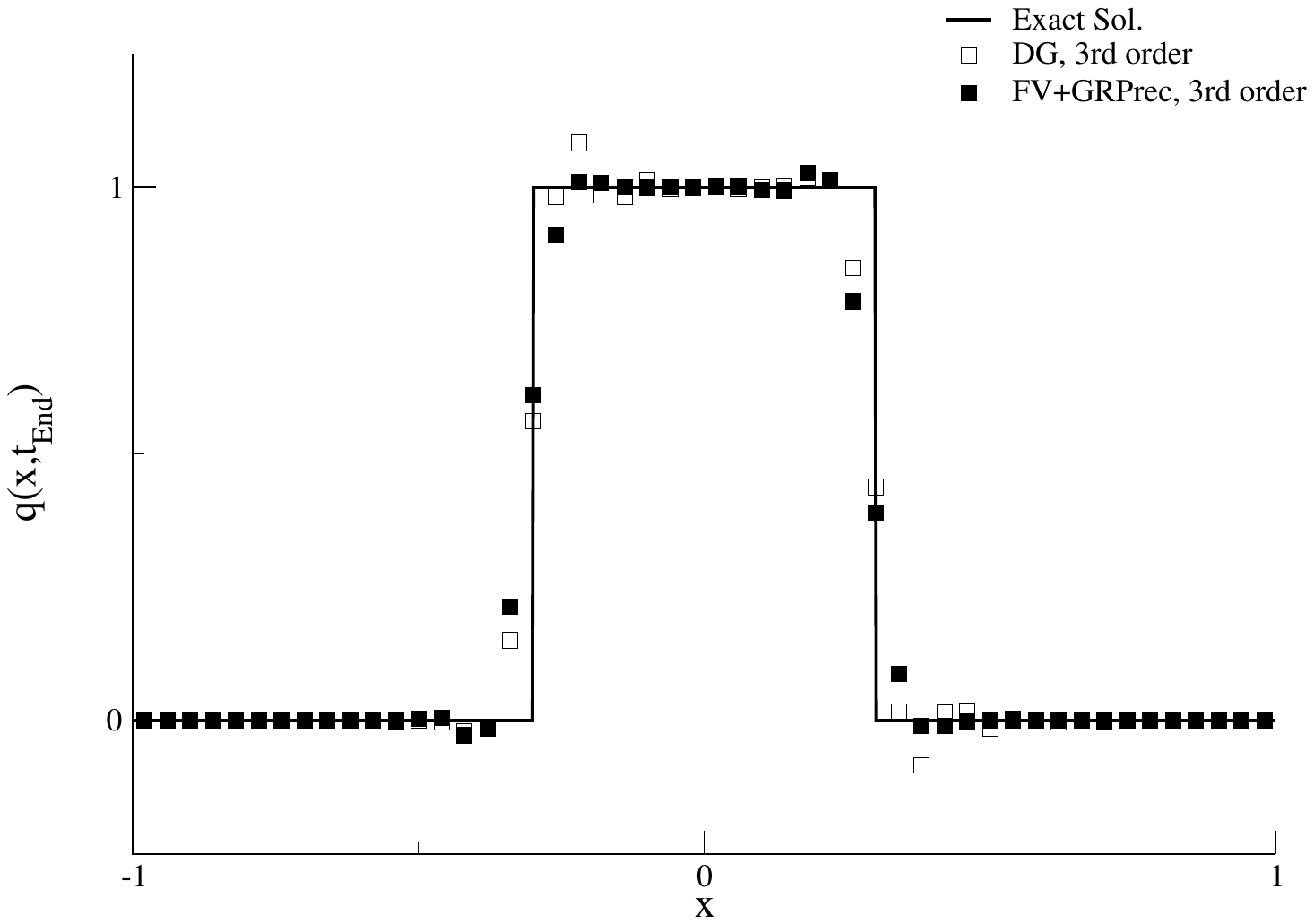}} 
	\quad
	\subfloat[][4th order]{\includegraphics[scale=0.22]{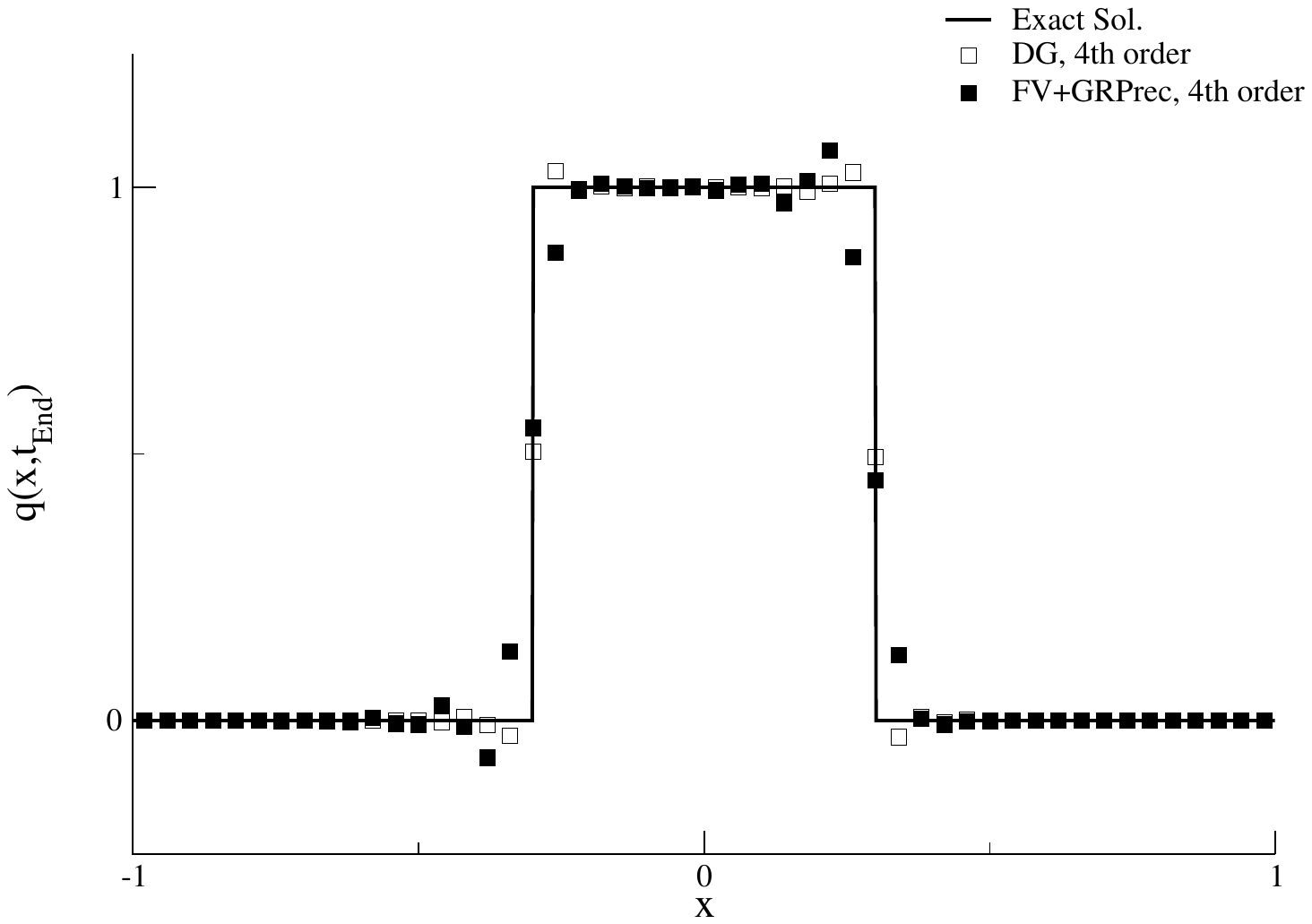}}
	\subfloat[][5th order]{ \includegraphics[scale=0.22]{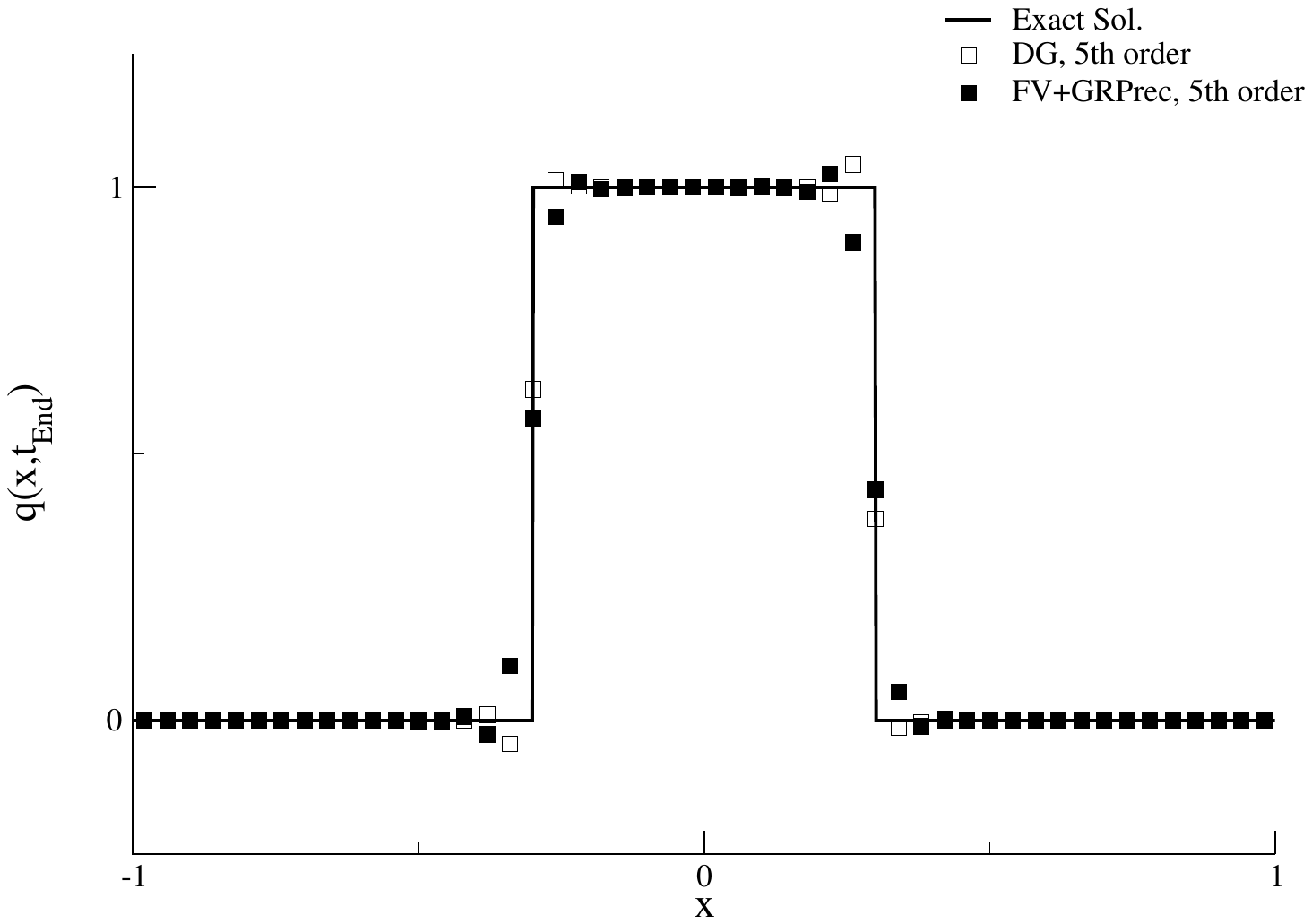}} 

	\caption{Linear advection, square-wave test with initial condition (\ref{eq:SW}): profiles for second, third, fourth and fifth order of accuracy for DG and FV+GRPrec. Parameters: $t_{End} = 4$,   $50$ cells and $C_{CFL}=0.9$.}
	\label{figure:LA-squarewave-Linear}
\end{figure}

\begin{figure}
	\centering
	\subfloat[][2nd order]{ \includegraphics[scale=0.22]{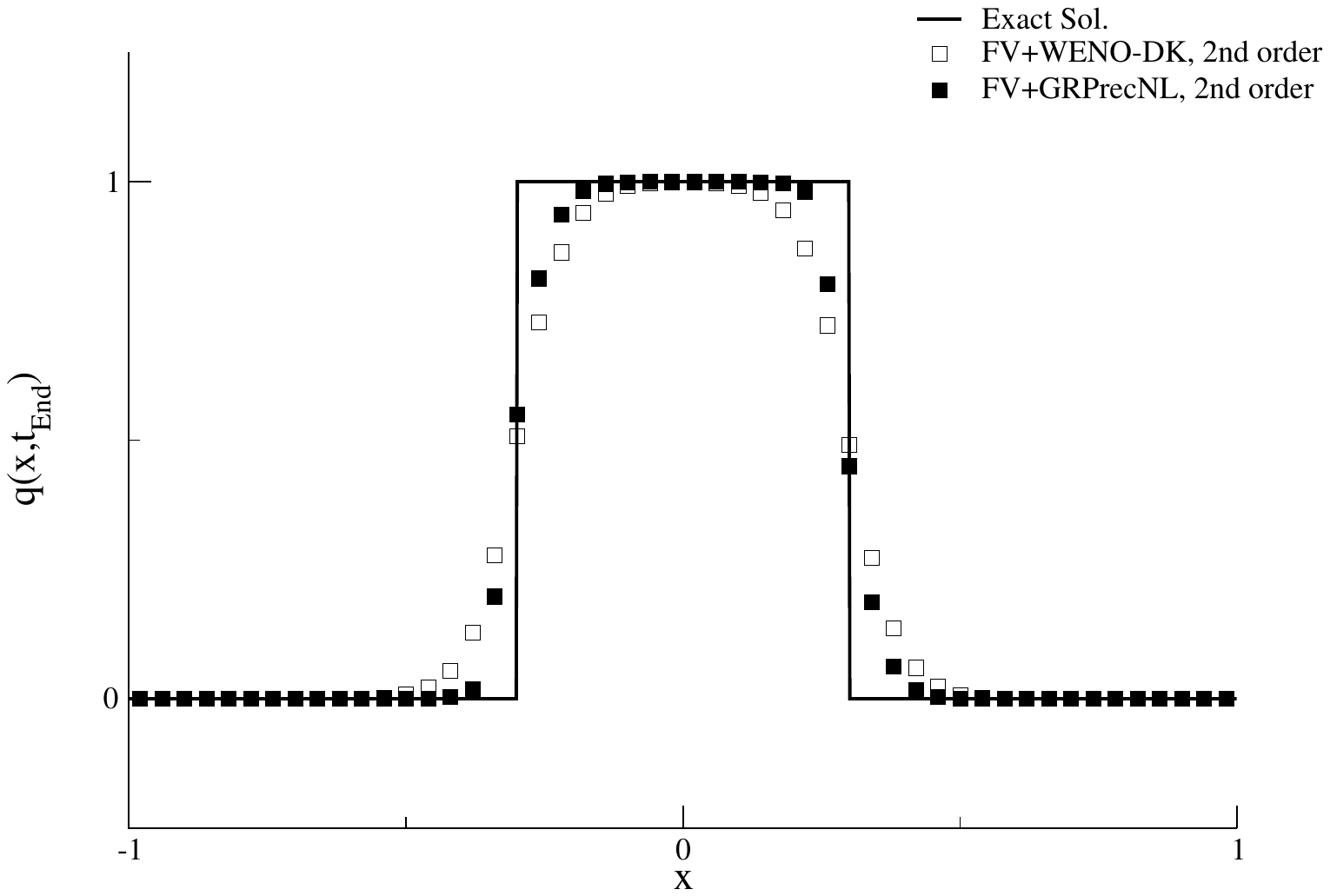}}
	\subfloat[][3rd order]{ \includegraphics[scale=0.22]{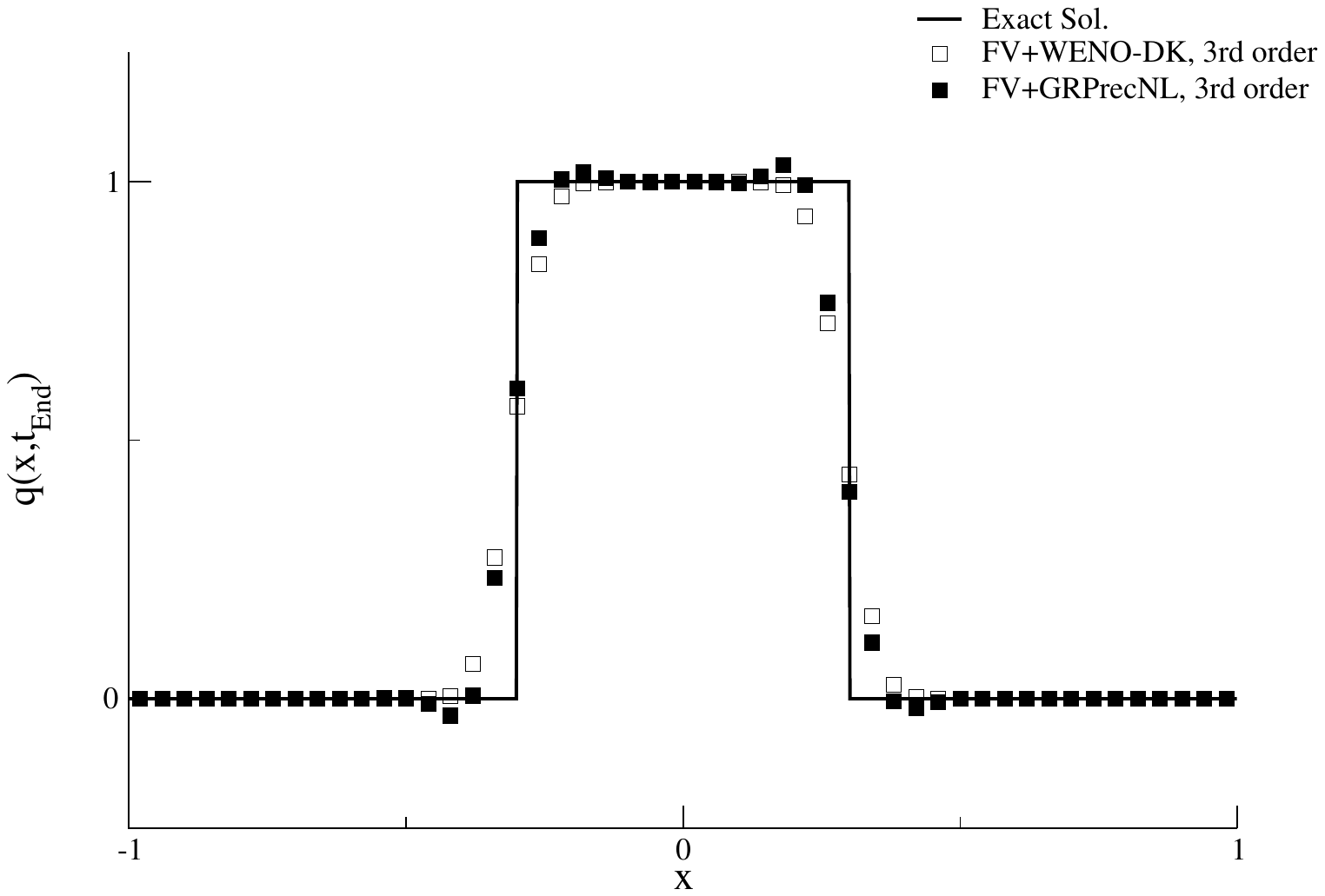}} 
	\quad
	\subfloat[][4th order]{\includegraphics[scale=0.22]{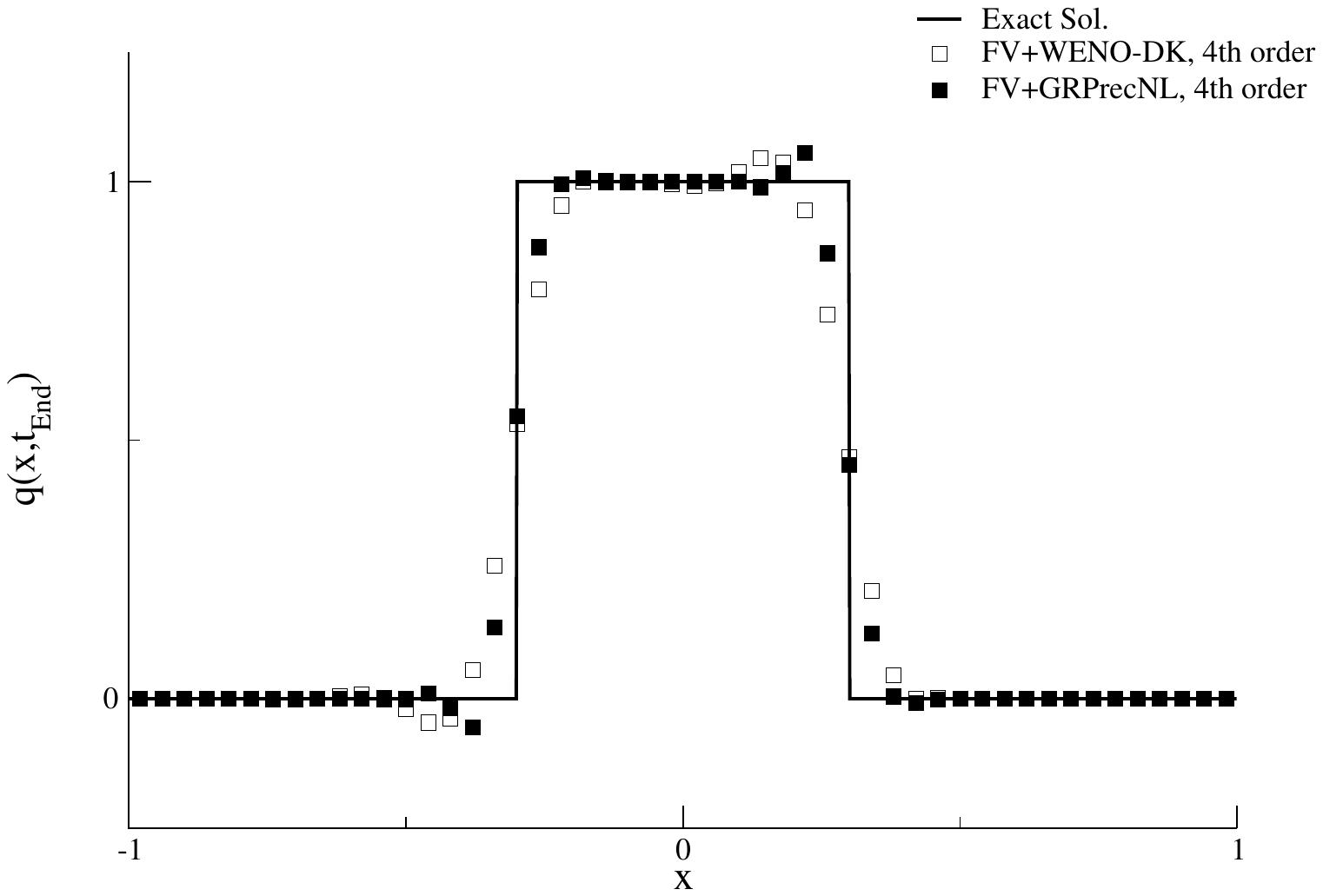}}
	\subfloat[][5th order]{ \includegraphics[scale=0.22]{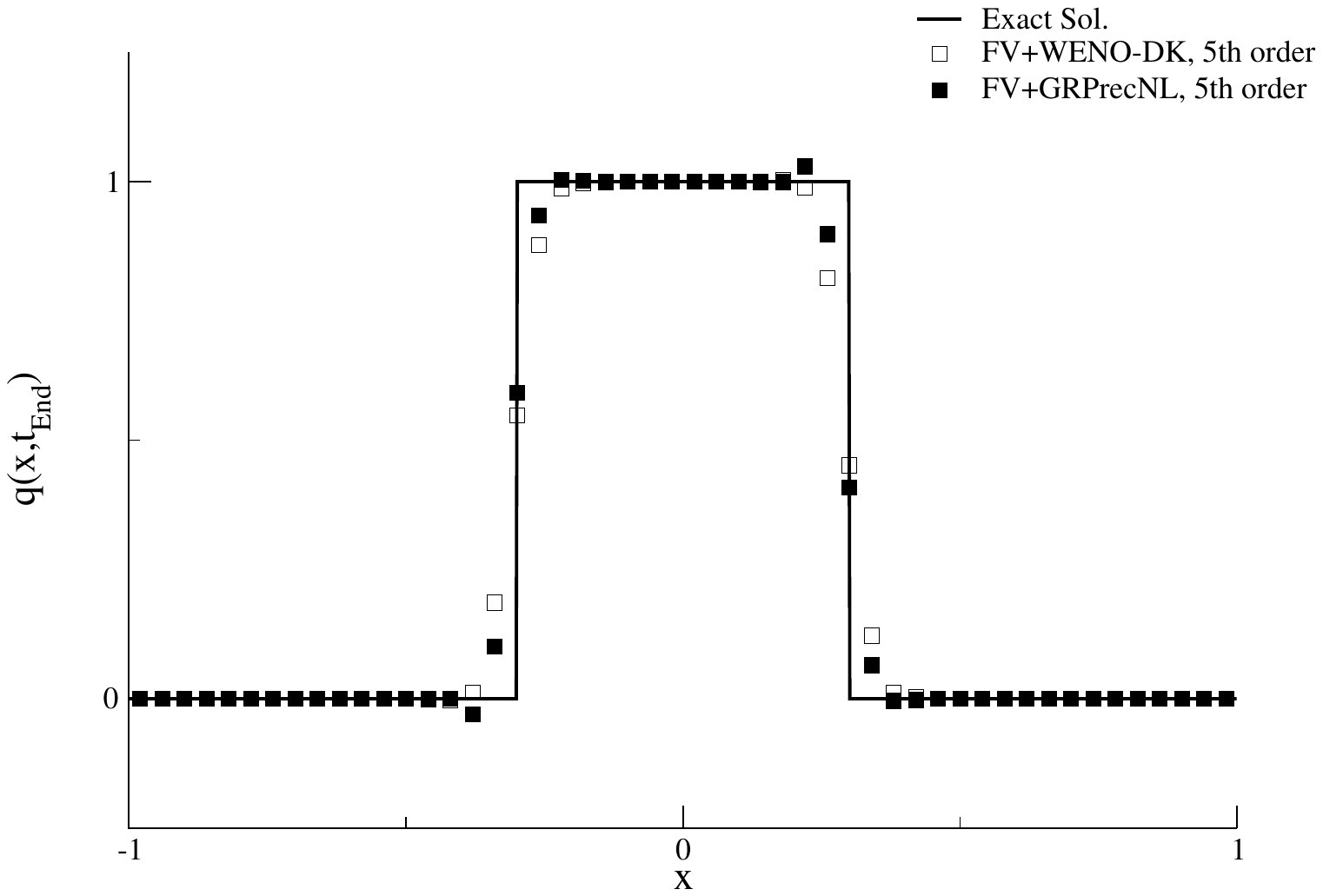}} 

	\caption{Linear advection, square-wave test with initial condition (\ref{eq:SW}): profiles for second, third, fourth and fifth order of accuracy for FV+WENO-DK and FV+GRPrecNL. Parameters: $t_{End} = 4$,   $50$ cells and $C_{CFL}=0.9$.}
	\label{figure:LA-squarewave-NonLinear}
\end{figure}

Here we solve the linear advection equation with initial condition 
\begin{eqnarray}
	\begin{array}{c}
		q(x,0) = \left\{
		\begin{array}{cc}

			1                        &, -0.3 \leq x \leq 0.3\;, \\

			0                                & otherwise\;.
		\end{array}
		
		\right.
	\end{array}
	\label{eq:SW}
\end{eqnarray}
See \cite{Montecinos:2022a}, for example. 
This test is employed to evaluate the ability of schemes to reproduce discontinuous profiles in which the Gibbs phenomenon affects high order schemes.  Fig. \ref{figure:LA-squarewave-Linear} shows results for DG and FV+GRPrec. Furthermore, Fig. \ref{figure:LA-squarewave-NonLinear} presents results for the FV+WENO-DK and FV+GRPrecNL. We see that results obtained with DG and FV+GRPrec are similar, with typical under- and overshoots near discontinuities, exhibiting an antisymmetric pattern in the location of spurious oscillations. We also observe that the use of GRPrecNL effectively attenuates spurious oscillations for the second order scheme, but is less effective in doing so than  WENO-DK, for higher-order schemes, with marginal beneficial effects. These results evidence that there is ample room for improvement in nonlinear version of our reconstruction method.

\subsubsection{Convergence rates and efficiency for linear advection} \label{sec:conv_lae}

Here we study empirical convergence rates for the linear advection equation using the  demanding test 
problem suggested by Shu \cite{Shu:1990a}.  The aim is to demonstrate empirically that the proposed schemes attain the theoretically expected orders of accuracy.  The initial condition is given by
\begin{eqnarray}
	\begin{array}{c}
		q(x,0) = \sin^4(\pi x) \;.
	\end{array}
	\label{eq:QS}
\end{eqnarray}
As for all other tests of this section, the spatial domain is $x \in [-1,1]$ with periodic boundary conditions,  and $\lambda=1$.

Tabs. \ref{Table-LAE-Quartic-MinimalRec-linear}  and \ref{Table-LAE-Quartic-MinimalRec-NonLinear} show computed empirical convergence rates for FV+GRPrec and FV+GRPrecNL. Furthermore, we provide Tabs. \ref{Table-LAE-Quartic-WENO-DK} and \ref{Table-LAE-Quartic-ADER-DG} in Section \ref{sec:app_conv_lae} of the Appendix, which contain empirical convergence rate results for FV+WENO-DK and DG schemes. We  observe that the expected order of accuracy is always obtained for all schemes and orders, with the exception of second-order FV+GRPrecNL. Interestingly, the same suboptimal behaviour is observed for FV+WENO-DK (see Tab. \ref{Table-LAE-Quartic-WENO-DK}). 

\begin{figure}
\centering
\subfloat[][2nd order]{ \includegraphics[scale=0.22]{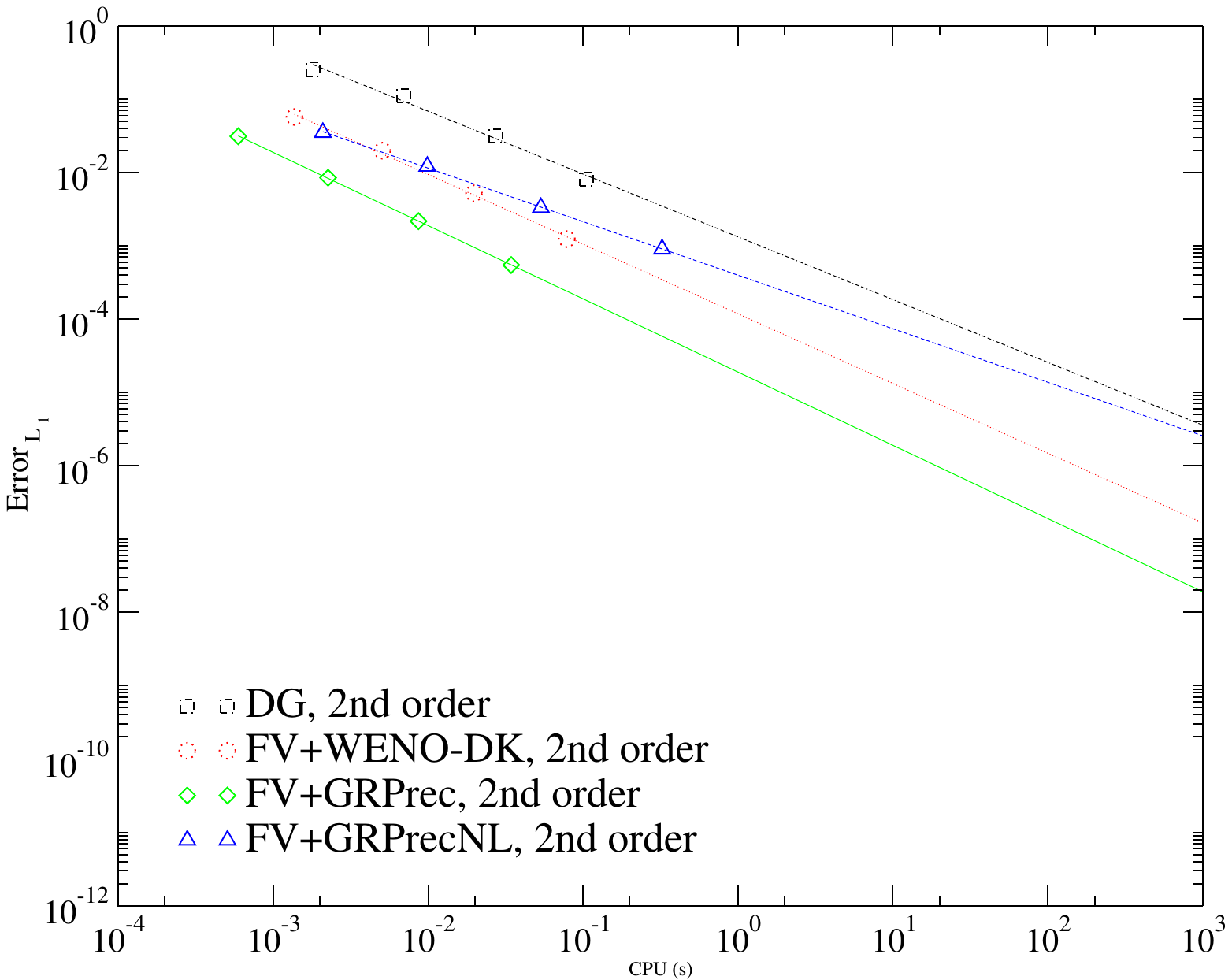}}
\subfloat[][3rd order]{ \includegraphics[scale=0.22]{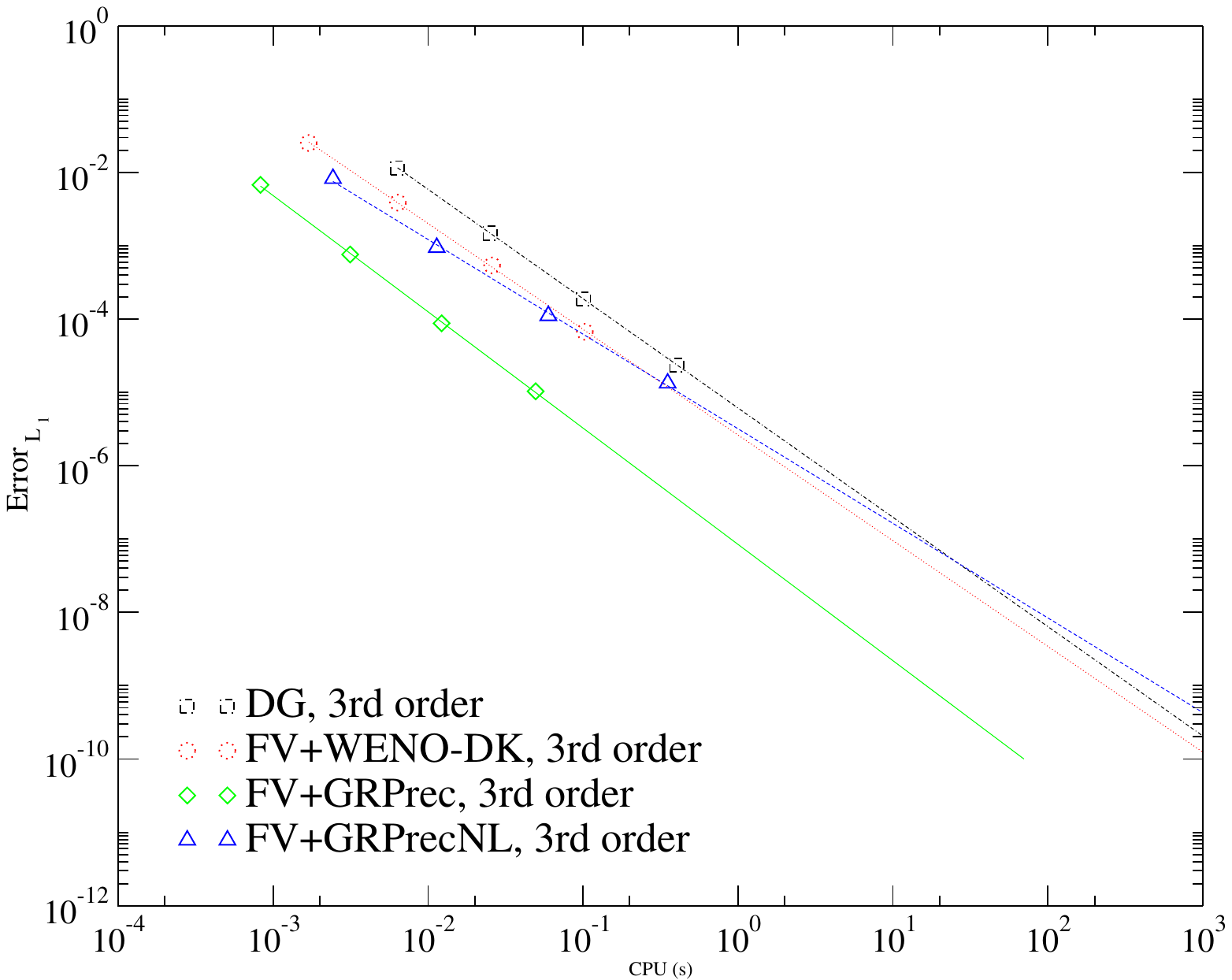}}
\quad
\subfloat[][4th order]{ \includegraphics[scale=0.22]{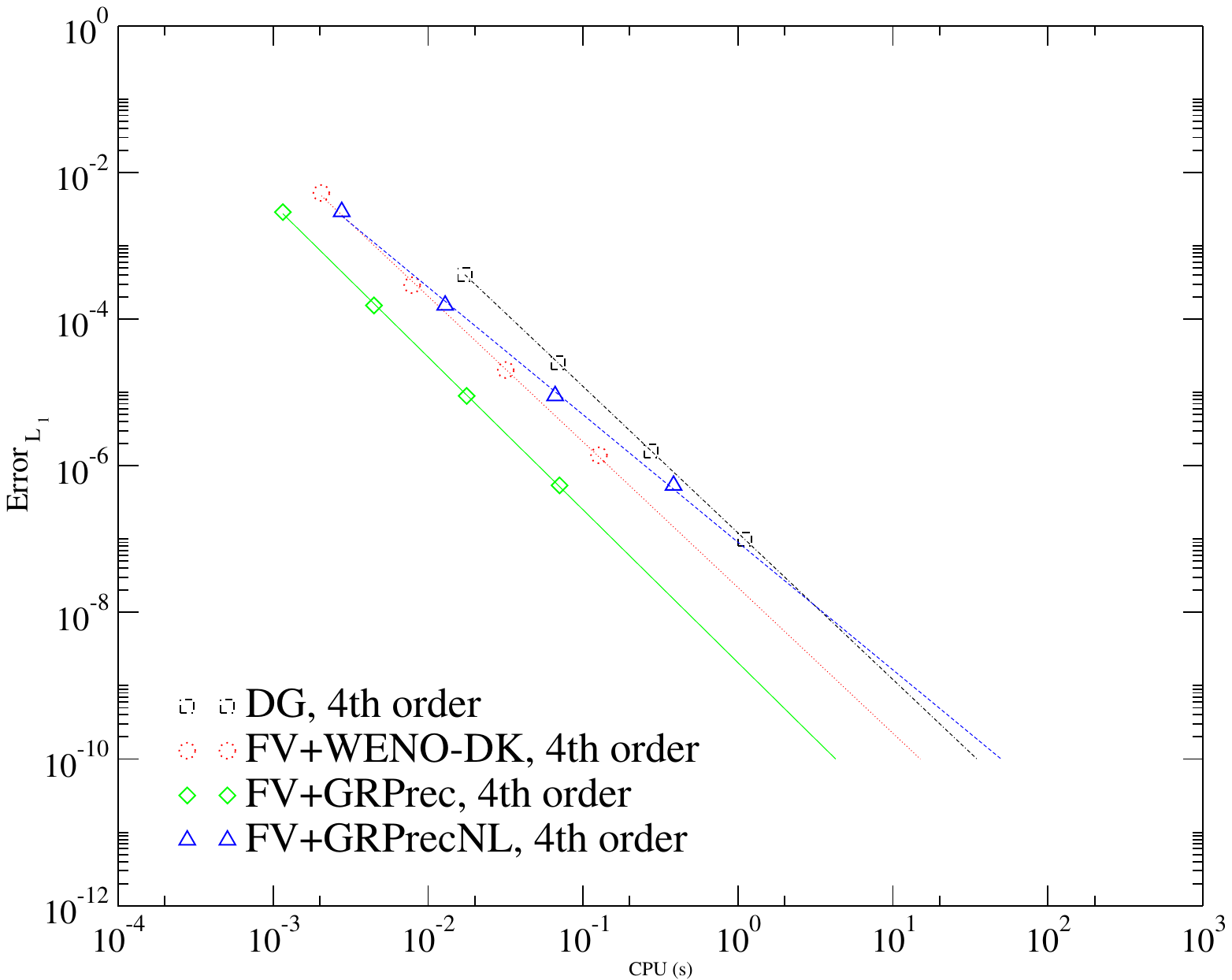}}
\subfloat[][5th order]{ \includegraphics[scale=0.22]{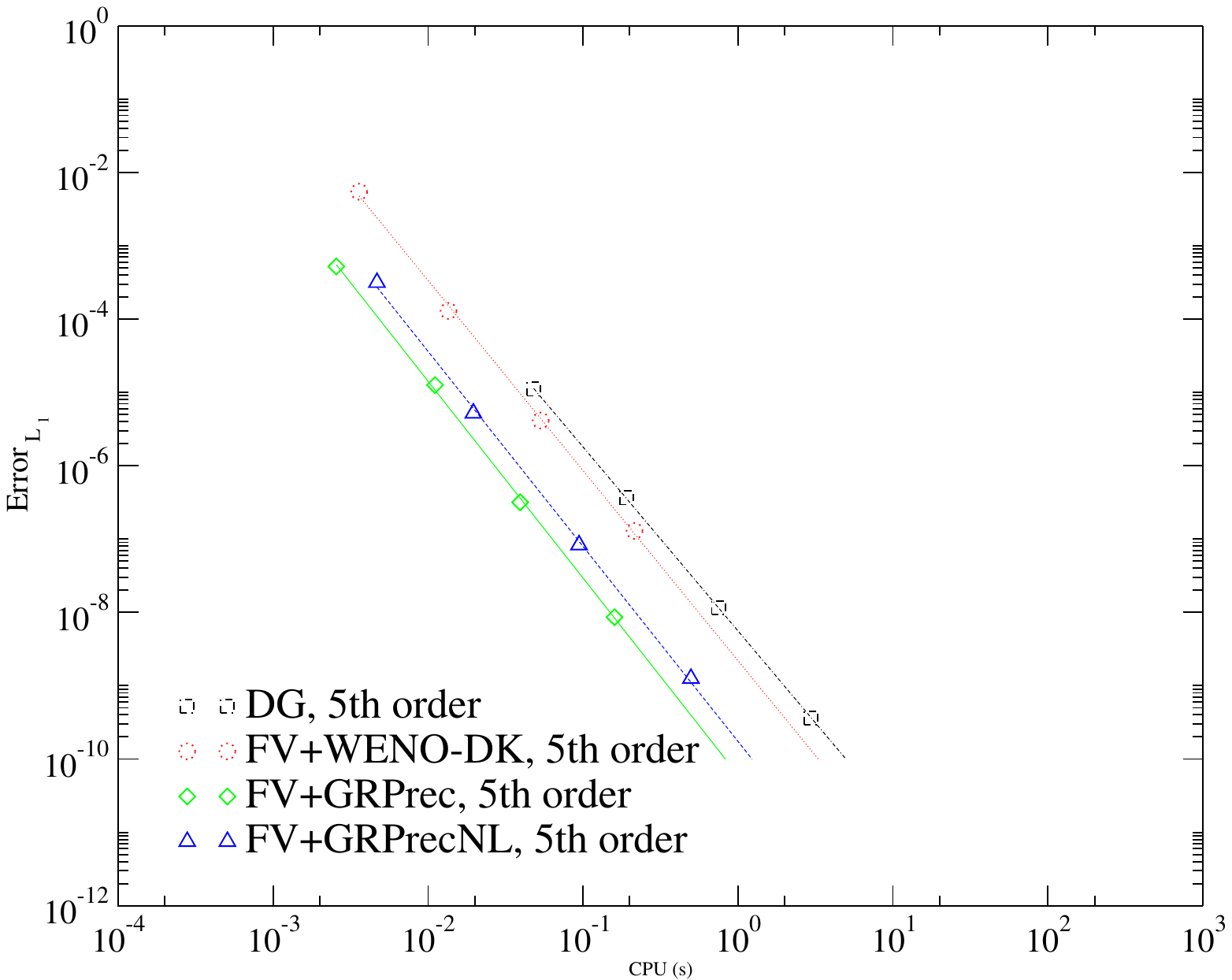}}

\caption{Efficiency plot for linear advection with initial condition 	(\ref{eq:QS}): $L_1-$error vs CPU time for second, third, fourth and fifth order of accuracy. Parameters: $t_{End} = 12$,  $C_{CFL}=0.9$.}
\label{figure:LA-Quartic-Efficiency}
\end{figure}

\begin{figure}
\centering
\includegraphics[scale=0.5]{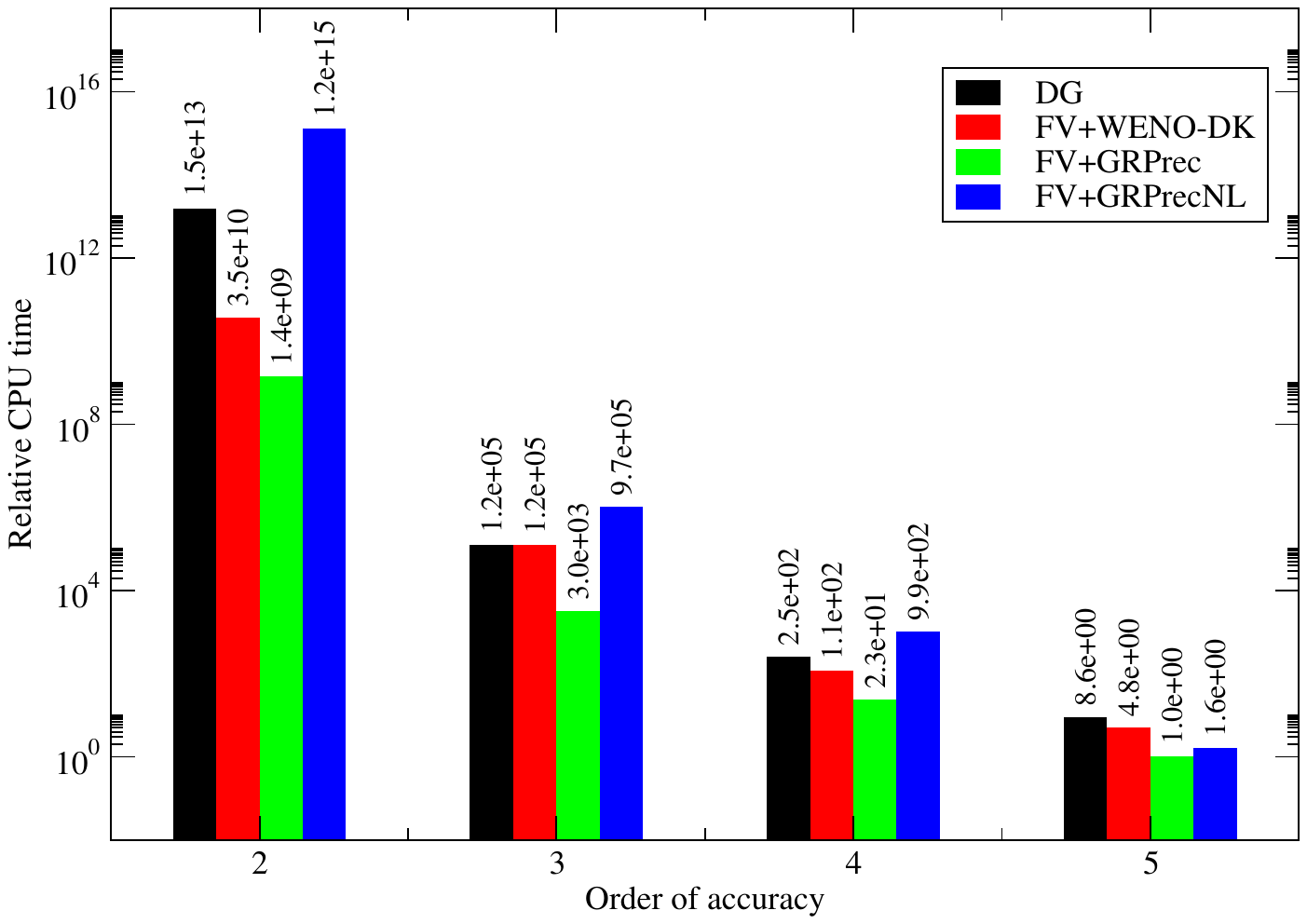}

\caption{Computational cost against order of accuracy for  linear advection with initial condition (\ref{eq:QS}): Relative CPU time with respect to fifth-order FV+GRPrec for an $L_1-$error of $10^{-16}$ for second-, third-, fourth- and fifth-order schemes. Parameter: $C_{CFL}=0.9$.}
\label{figure:LA-Quartic-Efficiency-Barplots}
\end{figure}

Fig. \ref{figure:LA-Quartic-Efficiency} is an efficiency plot; it shows the $L_1$-error versus  CPU time for FV+WENO-DK (empty circles), DG (full squares), FV+GRPrec (full circles), and FV+GRPrecNL (empty squares).  We note that GRPrec yields the best results, followed by the other schemes considered here. Notably, DG is always less efficient that FV+GRPrec. These observations are confirmed by results reported in Fig. \ref{figure:LA-Quartic-Efficiency-Barplots}, which shows relative CPU times with respect to that needed by the fifth-order FV+GRPrec to obtain an $L_1-$error of $10^{-16}$. Importantly, CPU times for such a small error were not directly computed, but extrapolated from data used to produce the efficiency plots shown in Fig. \ref{figure:LA-Quartic-Efficiency}. Besides evidencing the higher efficiency of higher order schemes over lower order ones, Fig. \ref{figure:LA-Quartic-Efficiency-Barplots} shows that, for all orders of accuracy considered, FV+GRPrec turns out to be significantly more efficient than the others. We note that GRPrecNL is the least efficient among all the schemes tested. This highlights a further limitation of our proposed non-linear variant of GRPrec, which merits future investigation. Nonetheless, we emphasize that for the Euler non-linear system, our GRPrecNL scheme yields encouraging results for solutions with discontinuities. More effective strategies for obtaining non-oscillatory schemes will be explored in future works; see also  related Sections \ref{sec:Euler}  and \ref{sec:conclusions}.


%
%

%
%


\begin{table}
\begin{center}
Theoretical order : 2, GRPrec \\
\begin{tabular}{cccccccc} 

\hline
\hline 
Mesh  & $L_\infty$ - ord  & $L_\infty$ - err  & $L_1$ - ord  & $L_1$ - err  & $L_2$ - ord  & $L_2$ - err  & CPU  \\  
\hline

    16  &  - &$  9.12e-0 2$&  -  &$  8.18e-0 2$&  -  &$  7.33e-0 2$  &0.0002\\
    32  &  1.50  &$  3.23e-0 2$&  1.39  &$  3.11e-0 2$&  1.49  &$  2.61e-0 2$  &0.0006\\
    64  &  1.91  &$  8.59e-0 3$&  1.88  &$  8.47e-0 3$&  1.89  &$  7.01e-0 3$  &0.0023\\
   128  &  1.97  &$  2.19e-0 3$&  1.97  &$  2.17e-0 3$&  1.97  &$  1.79e-0 3$  &0.0086\\
   256  &  1.99  &$  5.52e-0 4$&  1.99  &$  5.46e-0 4$&  1.99  &$  4.51e-0 4$  &0.0350\\




 \hline
 \\
\end{tabular}
\\
Theoretical order : 3, GRPrec \\
\begin{tabular}{cccccccc} 

\hline
\hline 
Mesh  & $L_\infty$ - ord  & $L_\infty$ - err  & $L_1$ - ord  & $L_1$ - err  & $L_2$ - ord  & $L_2$ - err  & CPU  \\  
\hline

    16  &  -  &$  4.00e-0 2$&  -  &$  5.26e-0 2$&  -  &$  3.99e-0 2$  &0.0002\\
    32  &  2.67  &$  6.29e-0 3$&  2.96  &$  6.77e-0 3$&  2.89  &$  5.38e-0 3$  &0.0008\\
    64  &  3.10  &$  7.34e-0 4$&  3.15  &$  7.61e-0 4$&  3.14  &$  6.09e-0 4$  &0.0031\\
   128  &  3.12  &$  8.46e-0 5$&  3.12  &$  8.77e-0 5$&  3.12  &$  6.99e-0 5$  &0.0122\\
   256  &  3.08  &$  1.00e-0 5$&  3.08  &$  1.03e-0 5$&  3.08  &$  8.25e-0 6$  &0.0495\\




 \hline
 \\
\end{tabular}
\\
Theoretical order : 4, GRPrec\\
\begin{tabular}{cccccccc} 

\hline
\hline 
Mesh  & $L_\infty$ - ord  & $L_\infty$ - err  & $L_1$ - ord  & $L_1$ - err  & $L_2$ - ord  & $L_2$ - err  & CPU  \\  
\hline

    16  &  -  &$  3.49e-0 2$&  -  &$  4.74e-0 2$&  -  &$  3.48e-0 2$  &0.0003\\
    32  &  3.89  &$  2.36e-0 3$&  4.04  &$  2.87e-0 3$&  3.96  &$  2.23e-0 3$  &0.0011\\
    64  &  4.17  &$  1.31e-0 4$&  4.23  &$  1.53e-0 4$&  4.20  &$  1.21e-0 4$  &0.0045\\
   128  &  4.10  &$  7.64e-0 6$&  4.10  &$  8.96e-0 6$&  4.10  &$  7.06e-0 6$  &0.0180\\
   256  &  4.06  &$  4.59e-0 7$&  4.06  &$  5.39e-0 7$&  4.06  &$  4.25e-0 7$  &0.0708\\




 \hline
 \\
\end{tabular}
\\
Theoretical order : 5, GRPrec \\
\begin{tabular}{cccccccc} 

\hline
\hline 
Mesh  & $L_\infty$ - ord  & $L_\infty$ - err  & $L_1$ - ord  & $L_1$ - err  & $L_2$ - ord  & $L_2$ - err  & CPU  \\  
\hline

    16  &  - &$  1.12e-0 2$&  -  &$  1.78e-0 2$&  -  &$  1.28e-0 2$  &0.0007\\
    32  &  4.82  &$  3.98e-0 4$&  5.09  &$  5.22e-0 4$&  4.99  &$  4.01e-0 4$  &0.0026\\
    64  &  5.26  &$  1.04e-0 5$&  5.37  &$  1.26e-0 5$&  5.34  &$  9.86e-0 6$  &0.0111\\
   128  &  5.29  &$  2.64e-0 7$&  5.31  &$  3.17e-0 7$&  5.30  &$  2.50e-0 7$  &0.0393\\
   256  &  5.20  &$  7.17e-0 9$&  5.21  &$  8.60e-0 9$&  5.21  &$  6.76e-0 9$  &0.1610\\




 \hline

\end{tabular}
\end{center}

\caption{Convergence rates for  linear advection with initial condition (\ref{eq:QS}).  Solution obtained with FV+GRPrec for second, third, fourth and fifth order of accuracy. Output time  $t_{out} = 4$, with  $C_{CFL}= 0.9$.}

  \label{Table-LAE-Quartic-MinimalRec-linear}
\end{table}
%


\begin{table}
\begin{center}
Theoretical order : 2, GRPrecNL\\
\begin{tabular}{cccccccc} 

\hline
\hline 
Mesh  & $L_\infty$ - ord  & $L_\infty$ - err  & $L_1$ - ord  & $L_1$ - err  & $L_2$ - ord  & $L_2$ - err  & CPU  \\  
\hline

    16  &  -  &$  9.19e-0 2$&  -  &$  8.66e-0 2$&  -  &$  7.58e-0 2$  &0.0005\\
    32  &  0.95  &$  4.77e-0 2$&  1.30  &$  3.50e-0 2$&  1.21  &$  3.29e-0 2$  &0.0021\\
    64  &  1.06  &$  2.28e-0 2$&  1.53  &$  1.21e-0 2$&  1.51  &$  1.15e-0 2$  &0.0099\\
   128  &  1.13  &$  1.04e-0 2$&  1.87  &$  3.31e-0 3$&  1.62  &$  3.76e-0 3$  &0.0535\\
   256  &  1.26  &$  4.33e-0 3$&  1.88  &$  9.02e-0 4$&  1.70  &$  1.16e-0 3$  &0.3201\\

  


 \hline
 \\
\end{tabular}
\\
Theoretical order : 3, GRPrecNL \\
\begin{tabular}{cccccccc} 

\hline
\hline 
Mesh  & $L_\infty$ - ord  & $L_\infty$ - err  & $L_1$ - ord  & $L_1$ - err  & $L_2$ - ord  & $L_2$ - err  & CPU  \\  
\hline

    16  & -  &$  5.03e-0 2$&  - &$  6.44e-0 2$&  -  &$  4.91e-0 2$  &0.0005\\
    32  &  2.70  &$  7.73e-0 3$&  2.96  &$  8.25e-0 3$&  2.90  &$  6.57e-0 3$  &0.0024\\
    64  &  3.08  &$  9.15e-0 4$&  3.12  &$  9.46e-0 4$&  3.11  &$  7.60e-0 4$  &0.0114\\
   128  &  3.09  &$  1.08e-0 4$&  3.09  &$  1.11e-0 4$&  3.10  &$  8.89e-0 5$  &0.0597\\
   256  &  3.06  &$  1.29e-0 5$&  3.06  &$  1.34e-0 5$&  3.06  &$  1.07e-0 5$  &0.3547\\

  


 \hline
 \\
\end{tabular}
\\
Theoretical order : 4, GRPrecNL\\
\begin{tabular}{cccccccc} 

\hline
\hline 
Mesh  & $L_\infty$ - ord  & $L_\infty$ - err  & $L_1$ - ord  & $L_1$ - err  & $L_2$ - ord  & $L_2$ - err  & CPU  \\  
\hline

    16  &  -  &$  3.74e-0 2$&  -  &$  5.20e-0 2$&  -  &$  3.82e-0 2$  &0.0006\\
    32  &  3.97  &$  2.38e-0 3$&  4.16  &$  2.92e-0 3$&  4.08  &$  2.25e-0 3$  &0.0028\\
    64  &  4.18  &$  1.31e-0 4$&  4.24  &$  1.54e-0 4$&  4.22  &$  1.21e-0 4$  &0.0130\\
   128  &  4.10  &$  7.65e-0 6$&  4.10  &$  8.99e-0 6$&  4.10  &$  7.08e-0 6$  &0.0671\\
   256  &  4.06  &$  4.60e-0 7$&  4.06  &$  5.40e-0 7$&  4.05  &$  4.26e-0 7$  &0.3836\\




 \hline
 \\
\end{tabular}
\\
Theoretical order : 5, GRPrecNL \\
\begin{tabular}{cccccccc} 

\hline
\hline 
Mesh  & $L_\infty$ - ord  & $L_\infty$ - err  & $L_1$ - ord  & $L_1$ - err  & $L_2$ - ord  & $L_2$ - err  & CPU  \\  
\hline

    16  &  -  &$  1.25e-0 2$& -  &$  1.69e-0 2$&  -  &$  1.25e-0 2$  &0.0010\\
    32  &  5.55  &$  2.68e-0 4$&  5.75  &$  3.15e-0 4$&  5.69  &$  2.43e-0 4$  &0.0046\\
    64  &  5.89  &$  4.52e-0 6$&  5.92  &$  5.22e-0 6$&  5.88  &$  4.12e-0 6$  &0.0196\\
   128  &  6.03  &$  6.92e-0 8$&  5.97  &$  8.31e-0 8$&  5.98  &$  6.53e-0 8$  &0.0931\\
   256  &  6.05  &$  1.04e-0 9$&  6.05  &$  1.25e-0 9$&  6.05  &$  9.85e-010$  &0.5016\\




 \hline

\end{tabular}
\end{center}
\caption{Convergence rates for linear advection with initial condition (\ref{eq:QS}).  Solution obtained with FV+GRPrecNL for second, third, fourth and fifth order of accuracy. Output time  $t_{out} = 4$, with  $C_{CFL}= 0.9$.}
  \label{Table-LAE-Quartic-MinimalRec-NonLinear}
\end{table}

\subsection{The Euler equations}
\label{sec:Euler}

In this section we apply our numerical methods to the one-dimensional Euler equations of gas dynamics 
\begin{eqnarray}
\label{eq:Euler-gen}
\begin{array}{c}
\partial_t \mathbf{Q} + \partial_x \mathbf{F}(\mathbf{Q}) = \mathbf{0} \;, 
\end{array}
\end{eqnarray}
with conserved variables and physical flux vectors, respectively given as
\begin{eqnarray}
\label{eq:Euler-variables}
\begin{array}{c}
\mathbf{Q} = \left(
\begin{array}{c}
\rho \\
\rho u \\
E
\end{array}
\right)\;,

\mathbf{F}(\mathbf{Q}) = \left(
\begin{array}{c}
\rho u \\
\rho u^2 + p \\
u(E+p)
\end{array}
\right)
\;.
\end{array}
\end{eqnarray}
Here $\rho$ is  density, $u$ is velocity, $p$ is pressure and $E$ is total energy given as
\begin{equation}
   E = \rho ( \frac{1}{2} u^{2} + e), 
\end{equation}
where $e$ is the specific internal energy given by an equation of state (a closure condition)
\begin{equation}
   e = e(\rho, p) \;.
\end{equation}
Other forms of the equation of state exist \cite{Toro:2009a}.  Here we adopt the ideal gas equation of state
\begin{equation}
   e = e(\rho, p) = \frac{p}{(\gamma-1) \rho} \;,
\end{equation}
where $\gamma$ is the ratio of specific heats, which under appropriate physical conditions may be taken as $\gamma =1.4$.

\subsubsection{Convergence rate study for the Euler equations} \label{sec:conv_euler}

We solve the Euler equations in the spatial domain $[0,1]$,  with periodic boundary conditions and the following initial conditions
\begin{eqnarray}
\begin{array}{c}
\rho(x,0) = 1 + 0.2 \sin(2 \pi x) \;,  u(x,0) = 1 \;, p(x,0) = 2 \;.
\end{array}
\end{eqnarray}
This problem has the exact solution 
\begin{eqnarray}
\begin{array}{c}
\rho(x,t) = 1 + 0.2 \sin(2 \pi (x-t)) \;, u(x,t) = 1, p(x,t) = 2 \;.
\end{array}
\end{eqnarray}
As for the linear advection equation, we compare results for FV+GRPrec and FV+GRPrecNL, with those for FV+WENO-DK, as well as for DG.  As in all other tests, the Courant number for FV+GRPrec, GRPrecNL, and WENO-DK is close to 1. For the DG scheme, the admissible Courant number depends on the order of accuracy (see Table \ref{tab:reconstruction_summary}).

\begin{figure}
\centering
\subfloat[][2nd order]{ \includegraphics[scale=0.22]{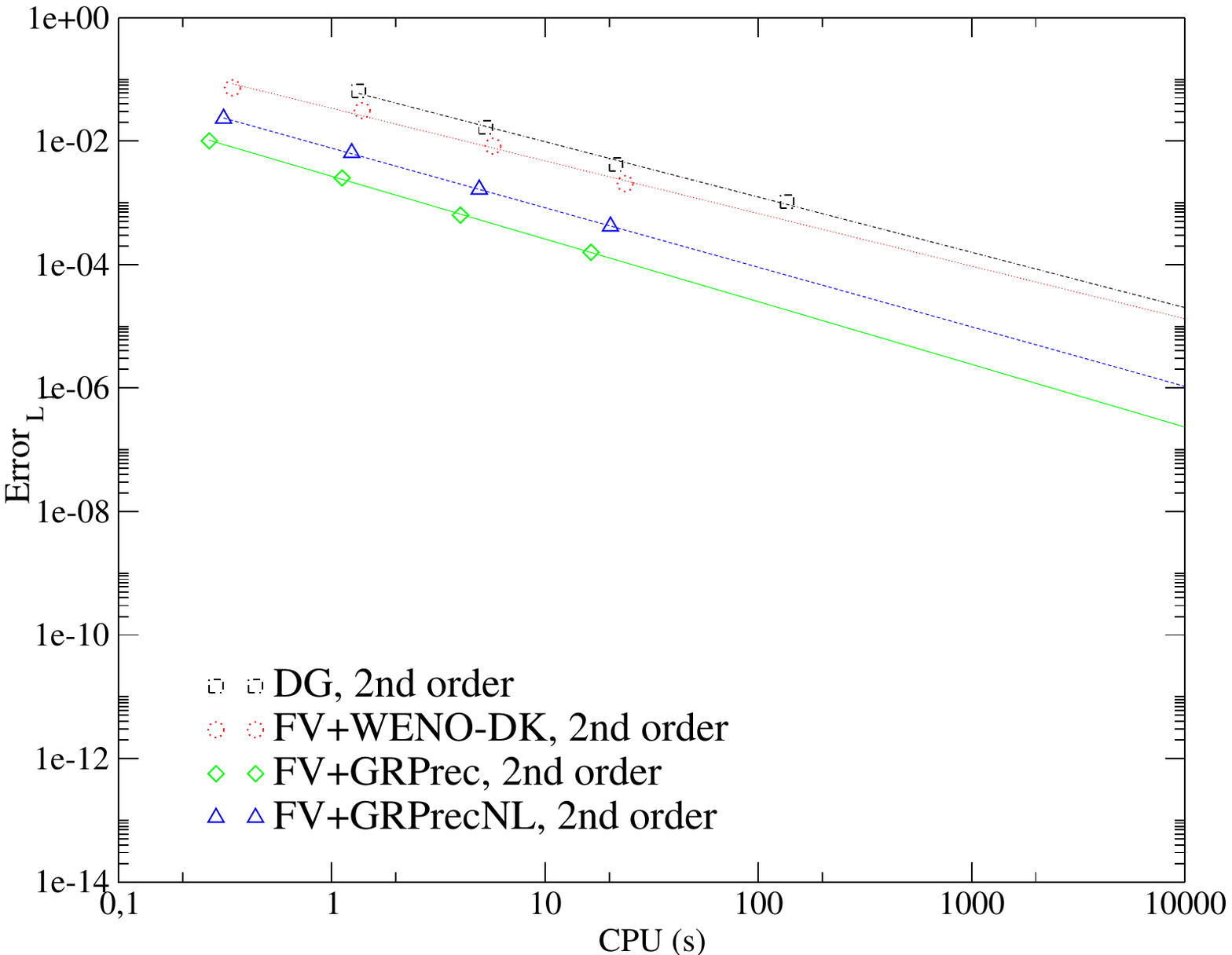}}
\subfloat[][3rd order]{ \includegraphics[scale=0.22]{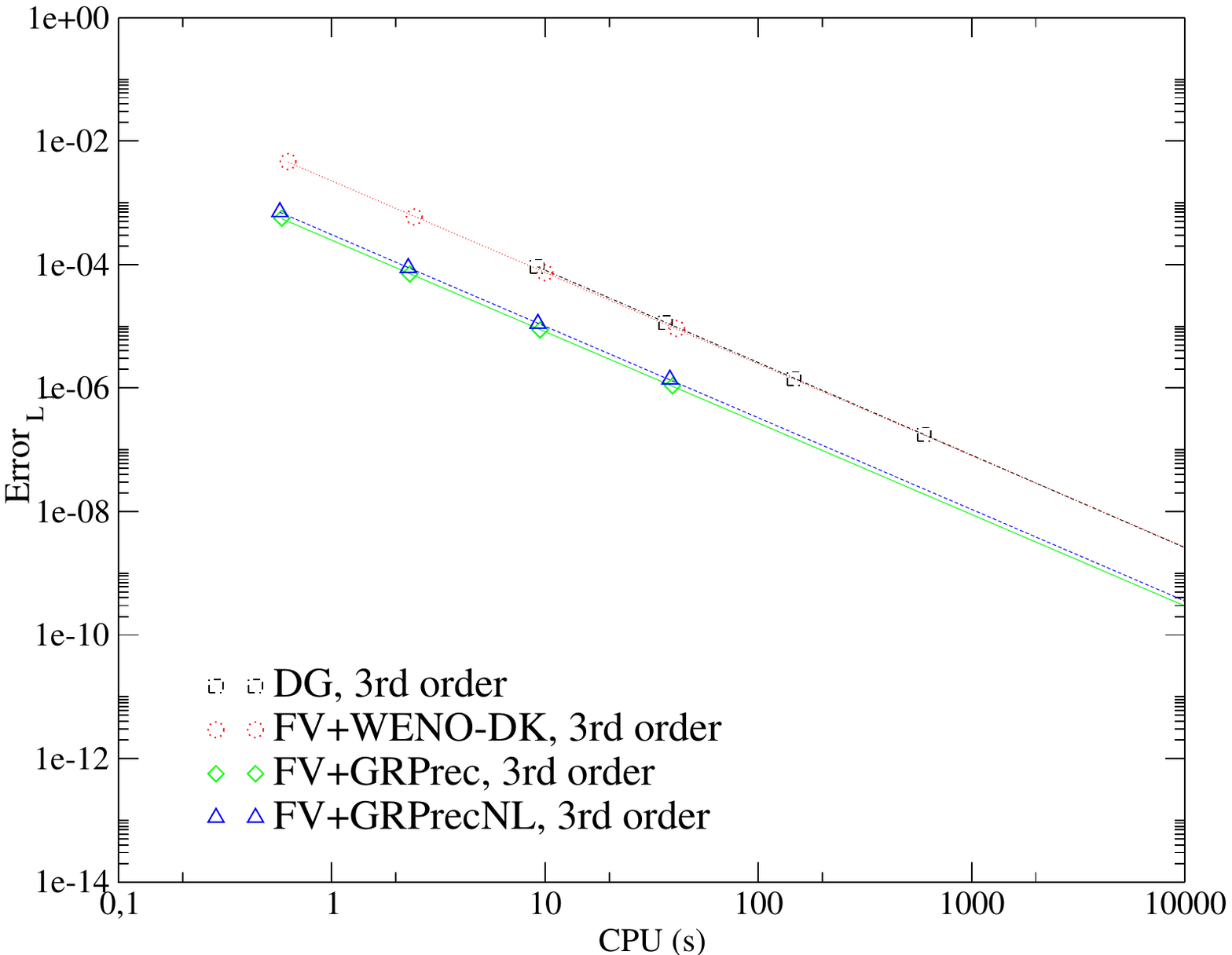}} 
\quad
\subfloat[][4th order]{\includegraphics[scale=0.22]{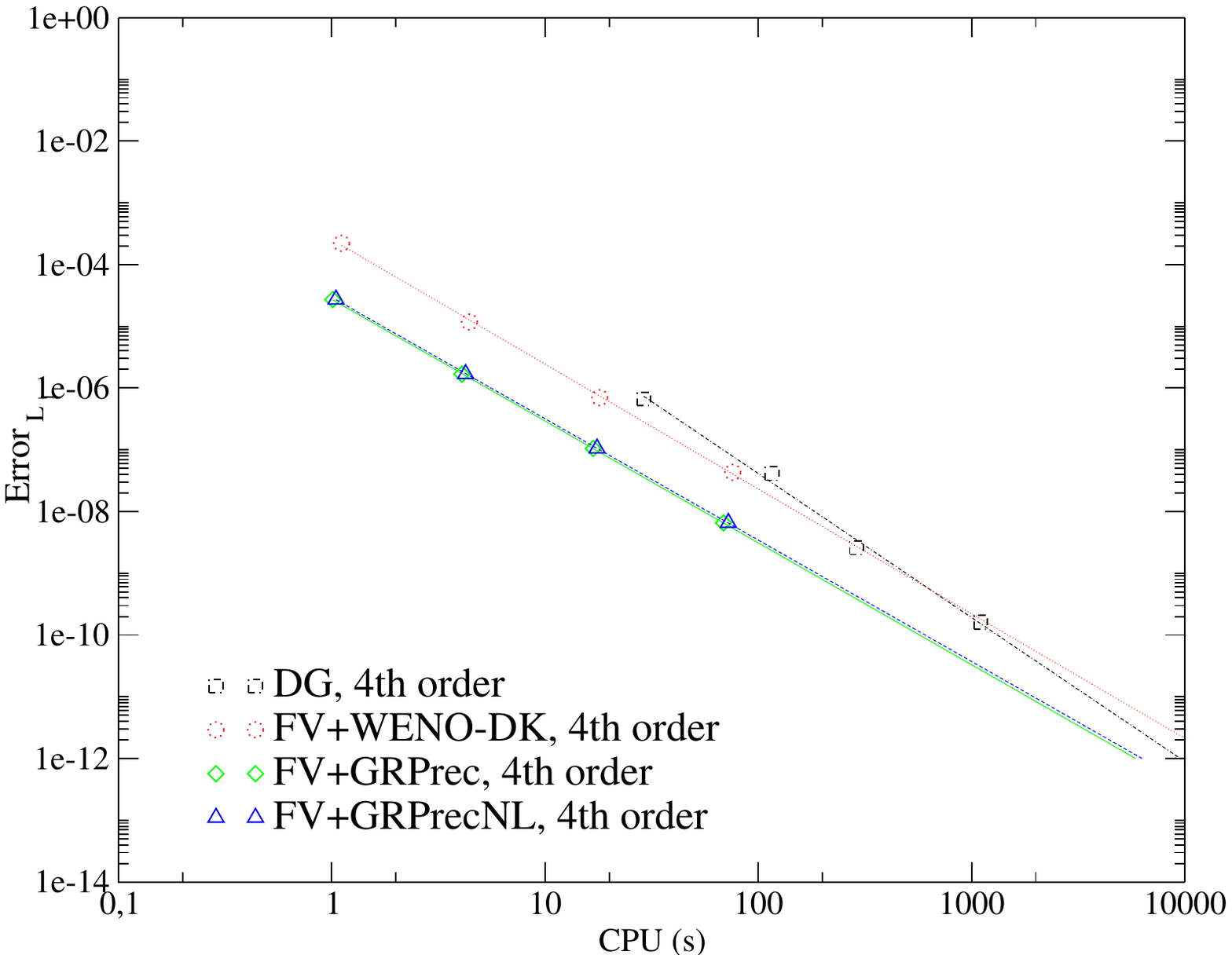}}
\subfloat[][5th order]{ \includegraphics[scale=0.22]{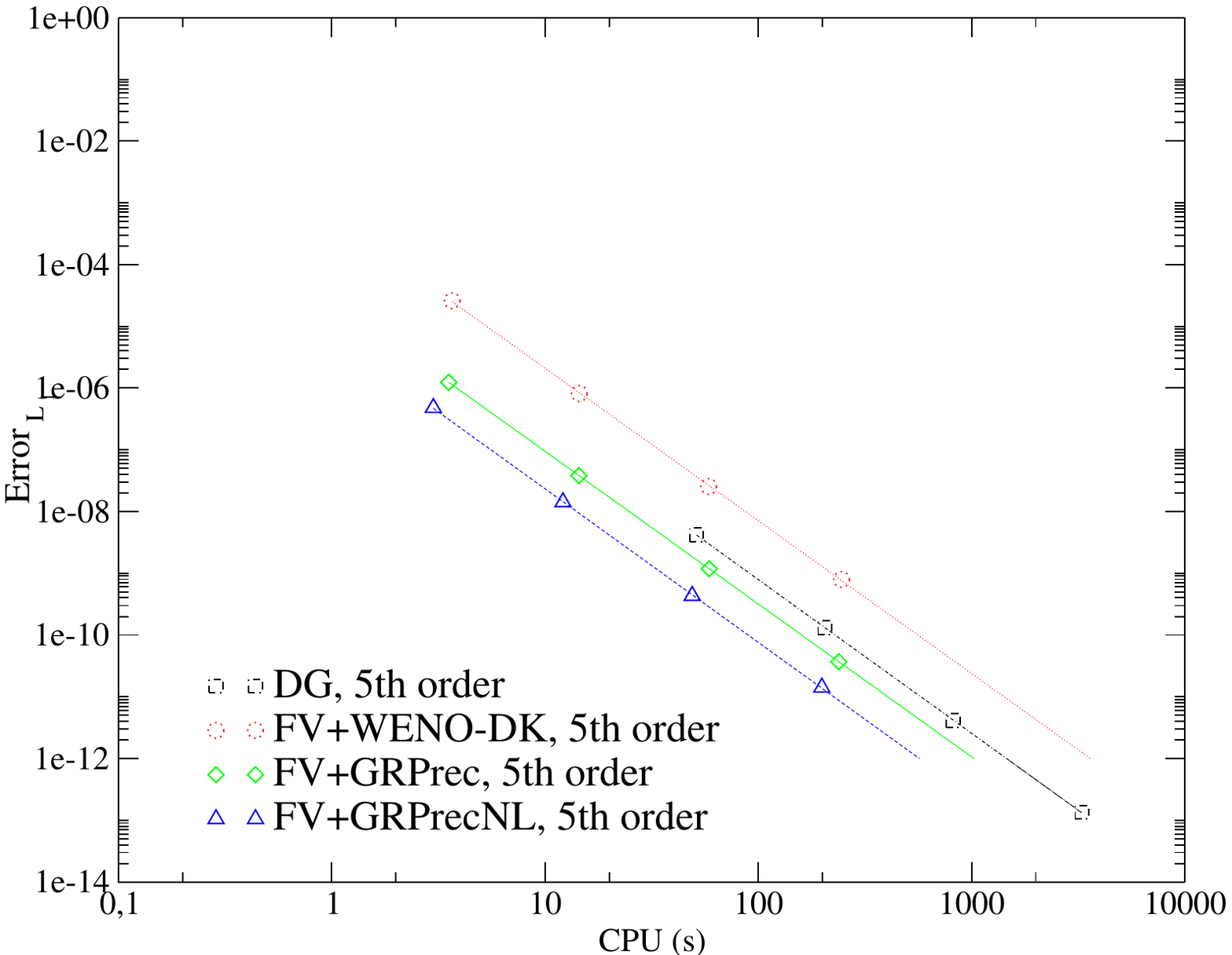}} 

\caption{Efficiency plots for the Euler equations: $L_1-$error vs CPU time for second, third, fourth and fifth order of accuracy. Parameters $t_{End} = 12$, $C_{CFL}=0.9$ ( $C_{CFL}=0.7$ for the fifth order only).}
\label{figure:EULER-Efficiency}
\end{figure}

\begin{figure}
\centering
\includegraphics[scale=0.5]{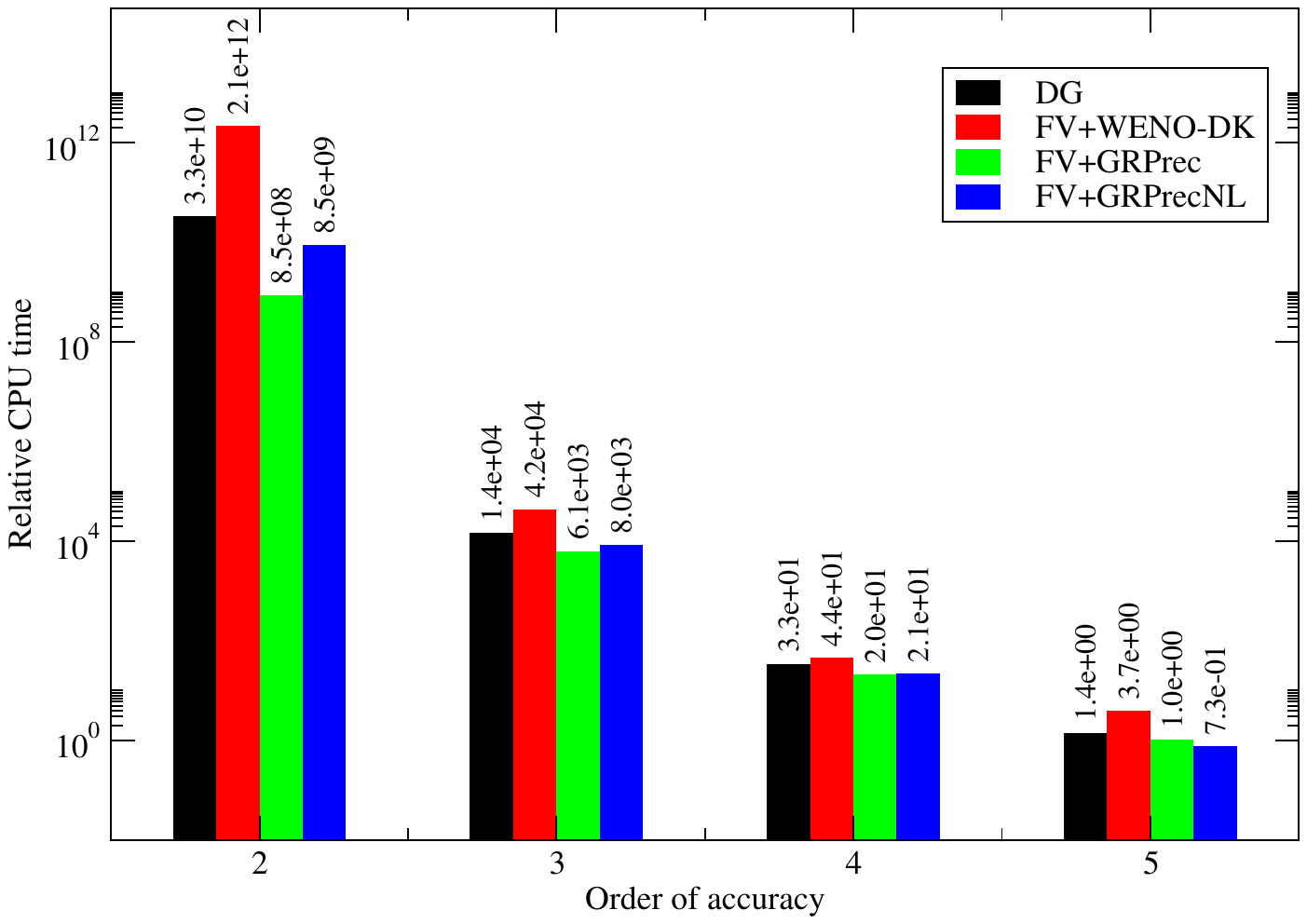}

\caption{Euler equations: Relative CPU time with respect to fifth-order FV+GRPrec for an $L_1-$error of $10^{-16}$ for second-, third-, fourth- and fifth-order schemes. Parameter: $C_{CFL}=0.9$  ( $C_{CFL}=0.7$ for the fifth order only).}
\label{figure:Euler-Efficiency-Barplots}
\end{figure}

Tabs. \ref{Table-Euler-MinimalRec-linear}  and \ref{Table-Euler-MinimalRec-NonLinear} show empirical convergence rates for the schemes FV+GRPrec and FV+GRPrecNL. Furthermore, we provide Tabs. \ref{Table-Euler-WENO-DK} and \ref{Table-Euler-ADER-DG} in Section \ref{sec:app_conv_euler} of the Appendix, which contain empirical convergence rate results for FV+WENO-DK and DG schemes. We can observe how the proposed reconstruction method results in schemes that provide the expected order of accuracy for all considered orders. Fig. \ref{figure:EULER-Efficiency} provides efficiency plots for this test. As for the linear advection case, also here we have that FV+GRPrec turns out to be the most efficient scheme, followed in this case by FV+GRPrecNL and then by DG. In this test we see that the least efficient scheme is FV+WENO-DK. This is expected since WENO-DK is the most computationally intensive reconstruction method. The efficiency of our novel FV+GRPrec scheme can be further appreciated in Fig. \ref{figure:Euler-Efficiency-Barplots}, which provides relative CPU times with respect to that needed by the fifth-order FV+GRPrec to obtain an $L_1-$error of $10^{-16}$.  Importantly,  as for the linear advection equation, CPU times for such a small error were not directly computed,  but extrapolated from computed data used to produce the efficiency plots shown in Fig. \ref{figure:EULER-Efficiency}.  Relative CPU times show that, especially for lower orders of accuracy, FV+GRPrec is orders-of-magnitude more efficient than other schemes considered in the exercise.



%
%


\begin{table}
\begin{center}
Theoretical order : 2, GRPrec \\
\begin{tabular}{cccccccc} 

\hline
\hline 
Mesh  & $L_\infty$ - ord  & $L_\infty$ - err  & $L_1$ - ord  & $L_1$ - err  & $L_2$ - ord  & $L_2$ - err  & CPU  \\  
\hline

    40  &  -  &$  5.64e-0 2$& - &$  7.19e-0 2$&  -  &$  5.65e-0 2$  &1.0588\\
    80  &  1.99  &$  1.42e-0 2$&  1.99  &$  1.81e-0 2$&  1.99  &$  1.42e-0 2$  &0.9665\\
   160  &  2.02  &$  3.51e-0 3$&  2.02  &$  4.48e-0 3$&  2.02  &$  3.51e-0 3$  &3.9799\\
   320  &  2.01  &$  8.75e-0 4$&  2.01  &$  1.11e-0 3$&  2.01  &$  8.75e-0 4$  &16.7539\\
   640  &  2.00  &$  2.18e-0 4$&  2.00  &$  2.78e-0 4$&  2.00  &$  2.19e-0 4$  &75.5874\\


 \hline
 \\
\end{tabular}
\\
Theoretical order : 3, GRPrec \\
\begin{tabular}{cccccccc} 

\hline
\hline 
Mesh  & $L_\infty$ - ord  & $L_\infty$ - err  & $L_1$ - ord  & $L_1$ - err  & $L_2$ - ord  & $L_2$ - err  & CPU  \\  
\hline

    40  &  -  &$  1.37e-0 3$&  -  &$  1.74e-0 3$&  -  &$  1.37e-0 3$  &1.0884\\
    80  &  3.23  &$  1.46e-0 4$&  3.23  &$  1.86e-0 4$&  3.23  &$  1.46e-0 4$  &1.5466\\
   160  &  3.07  &$  1.74e-0 5$&  3.07  &$  2.21e-0 5$&  3.07  &$  1.74e-0 5$  &6.4663\\
   320  &  3.02  &$  2.14e-0 6$&  3.02  &$  2.73e-0 6$&  3.02  &$  2.14e-0 6$  &25.4359\\
   640  &  3.00  &$  2.67e-0 7$&  3.01  &$  3.40e-0 7$&  3.01  &$  2.67e-0 7$  &102.5059\\


 \hline
 \\
\end{tabular}
\\
Theoretical order : 4, GRPrec\\
\begin{tabular}{cccccccc} 

\hline
\hline 
Mesh  & $L_\infty$ - ord  & $L_\infty$ - err  & $L_1$ - ord  & $L_1$ - err  & $L_2$ - ord  & $L_2$ - err  & CPU  \\  
\hline

    40  &  - &$  3.85e-0 4$&  -  &$  4.93e-0 4$&  -  &$  3.87e-0 4$  &1.3020\\
    80  &  4.30  &$  1.96e-0 5$&  4.30  &$  2.50e-0 5$&  4.30  &$  1.96e-0 5$  &2.4040\\
   160  &  4.11  &$  1.14e-0 6$&  4.11  &$  1.45e-0 6$&  4.11  &$  1.14e-0 6$  &9.6728\\
   320  &  4.03  &$  6.98e-0 8$&  4.03  &$  8.88e-0 8$&  4.03  &$  6.98e-0 8$  &38.0699\\
   640  &  4.01  &$  4.34e-0 9$&  4.01  &$  5.52e-0 9$&  4.01  &$  4.34e-0 9$  &145.9692\\


 \hline
 \\
\end{tabular}
\\
Theoretical order : 5, GRPrec \\
\begin{tabular}{cccccccc} 

\hline
\hline 
Mesh  & $L_\infty$ - ord  & $L_\infty$ - err  & $L_1$ - ord  & $L_1$ - err  & $L_2$ - ord  & $L_2$ - err  & CPU  \\  
\hline

    40  &  -  &$  6.24e-0 6$&  -  &$  7.96e-0 6$&  -  &$  6.25e-0 6$  &2.1063\\
    80  &  5.96  &$ 10.00e-0 8$&  5.97  &$  1.27e-0 7$&  5.97  &$  1.00e-0 7$  &5.8938\\
   160  &  5.92  &$  1.65e-0 9$&  5.92  &$  2.10e-0 9$&  5.92  &$  1.65e-0 9$  &24.4946\\
   320  &  5.75  &$  3.06e-011$&  5.76  &$  3.89e-011$&  5.76  &$  3.06e-011$  &97.4251\\
   640  &  5.30  &$  7.75e-013$&  5.49  &$  8.64e-013$&  5.49  &$  6.80e-013$  &368.7965\\


 \hline

\end{tabular}
\end{center}

\caption{Euler equations.  Solution obtained with FV+GRPrec for second, third, fourth and fifth order of accuracy. Output  time  $t_{out} = 4$, with  $C_{CFL}= 0.9$.}

  \label{Table-Euler-MinimalRec-linear}
\end{table}
%


\begin{table}
\begin{center}
Theoretical order : 2, GRPrecNL\\
\begin{tabular}{cccccccc} 

\hline
\hline 
Mesh  & $L_\infty$ - ord  & $L_\infty$ - err  & $L_1$ - ord  & $L_1$ - err  & $L_2$ - ord  & $L_2$ - err  & CPU  \\  
\hline

    40  &  -  &$  6.58e-0 2$&  -  &$  9.49e-0 2$& -&$  7.23e-0 2$  &0.9683\\
    80  &  1.86  &$  1.81e-0 2$&  2.08  &$  2.24e-0 2$&  2.05  &$  1.75e-0 2$  &1.0463\\
   160  &  1.69  &$  5.62e-0 3$&  1.96  &$  5.78e-0 3$&  1.92  &$  4.61e-0 3$  &4.4493\\
   320  &  1.60  &$  1.85e-0 3$&  1.89  &$  1.56e-0 3$&  1.93  &$  1.21e-0 3$  &19.0080\\
   640  &  1.58  &$  6.19e-0 4$&  1.98  &$  3.95e-0 4$&  1.93  &$  3.19e-0 4$  &80.5597\\


 \hline
 \\
\end{tabular}
\\
Theoretical order : 3, GRPrecNL \\
\begin{tabular}{cccccccc} 

\hline
\hline 
Mesh  & $L_\infty$ - ord  & $L_\infty$ - err  & $L_1$ - ord  & $L_1$ - err  & $L_2$ - ord  & $L_2$ - err  & CPU  \\  
\hline

    40  &  -  &$  1.28e-0 3$&  -  &$  1.63e-0 3$&  -  &$  1.29e-0 3$  &1.1142\\
    80  &  3.55  &$  1.09e-0 4$&  3.55  &$  1.40e-0 4$&  3.56  &$  1.10e-0 4$  &1.3888\\
   160  &  3.22  &$  1.17e-0 5$&  3.22  &$  1.50e-0 5$&  3.22  &$  1.17e-0 5$  &6.0211\\
   320  &  3.06  &$  1.40e-0 6$&  3.07  &$  1.79e-0 6$&  3.07  &$  1.40e-0 6$  &24.9329\\
   640  &  3.02  &$  1.73e-0 7$&  3.02  &$  2.21e-0 7$&  3.02  &$  1.73e-0 7$  &103.2699\\


 \hline
 \\
\end{tabular}
\\
Theoretical order : 4, GRPrecNL\\
\begin{tabular}{cccccccc} 

\hline
\hline 
Mesh  & $L_\infty$ - ord  & $L_\infty$ - err  & $L_1$ - ord  & $L_1$ - err  & $L_2$ - ord  & $L_2$ - err  & CPU  \\  
\hline

    40  &  - &$  3.86e-0 4$&  -  &$  4.99e-0 4$&  -  &$  3.90e-0 4$  &1.2993\\
    80  &  4.29  &$  1.98e-0 5$&  4.30  &$  2.53e-0 5$&  4.30  &$  1.98e-0 5$  &1.9955\\
   160  &  4.10  &$  1.15e-0 6$&  4.11  &$  1.46e-0 6$&  4.11  &$  1.15e-0 6$  &7.7314\\
   320  &  4.03  &$  7.04e-0 8$&  4.03  &$  8.96e-0 8$&  4.03  &$  7.04e-0 8$  &34.3729\\
   640  &  4.01  &$  4.38e-0 9$&  4.01  &$  5.57e-0 9$&  4.01  &$  4.38e-0 9$  & 142.3519\\


 \hline
 \\
\end{tabular}
\\
Theoretical order : 5, GRPrecNL \\
\begin{tabular}{cccccccc} 

\hline
\hline 
Mesh  & $L_\infty$ - ord  & $L_\infty$ - err  & $L_1$ - ord  & $L_1$ - err  & $L_2$ - ord  & $L_2$ - err  & CPU  \\  
\hline

    40  &  -  &$  2.79e-0 5$& -  &$  3.54e-0 5$&  -  &$  2.79e-0 5$  &2.1213\\
    80  &  4.95  &$  9.01e-0 7$&  4.95  &$  1.15e-0 6$&  4.95  &$  9.01e-0 7$  &5.2640\\
   160  &  4.99  &$  2.84e-0 8$&  4.99  &$  3.61e-0 8$&  4.99  &$  2.84e-0 8$  &23.859\\
   320  &  5.00  &$  8.89e-010$&  5.00  &$  1.13e-0 9$&  5.00  &$  8.89e-010$  &102.0512\\
   640  &  3.76  &$  6.58e-011$&  4.92  &$  3.73e-011$&  4.88  &$  3.02e-011$  &390.4961\\


 \hline

\end{tabular}
\end{center}
\caption{Euler equations.  Solution obtained with FV+GRPrecNL for second, third, fourth and fifth order of accuracy. Output  time  $t_{out} = 4$, with  $C_{CFL}= 0.9$.}
  \label{Table-Euler-MinimalRec-NonLinear}
\end{table}
\subsubsection{Riemann Problems for Euler equations}

Here we solve Riemann problems for the Euler equations (\ref{eq:Euler-gen})-(\ref{eq:Euler-variables}),  with initial conditions
\begin{eqnarray}
\begin{array}{c}
\mathbf{Q}(x,0) = \left\{
\begin{array}{cc}
(\rho_L, \rho_L u_L, E_L)\;,& x \leq x_c \;, \\
(\rho_R, \rho_R u_L, E_R)\;,& x > x_c \;.
\end{array}
\right.
\end{array}
\end{eqnarray}
Tab. \ref{table:initial-cond-EulerRP} shows the initial conditions for the two tests considered here. The first one is Sod's problem \cite{Sod:1978a} and the second one is the so-called {\it 123 problem} \cite{Einfeld:1991a}.

\begin{table}
\begin{center}

\begin{tabular}{c|ccc|ccc|c}
Test & $\rho_L$ & $u_L$ & $p_L$ & $\rho_R$ & $u_R$ & $p_R$ & $x_c$ \\
\hline
Sod problem & 1 & 0 & 1 & 0.125 & 0 & 0.1 & 0.3 \\
123 problem & 1 & -2 & 0.4 & 1 & 2 & 0.4 & 0.5 \\
\hline
\end{tabular}
\end{center}
\caption{Data for Riemann problems for the Euler equations. Initial conditions in terms of density, velocity and pressure ($L/R$); $x_c$ denotes the position of the  initial discontinuity.}
\label{table:initial-cond-EulerRP}
\end{table}

\begin{figure}
\centering
\subfloat[][2nd order]{ \includegraphics[scale=0.22]{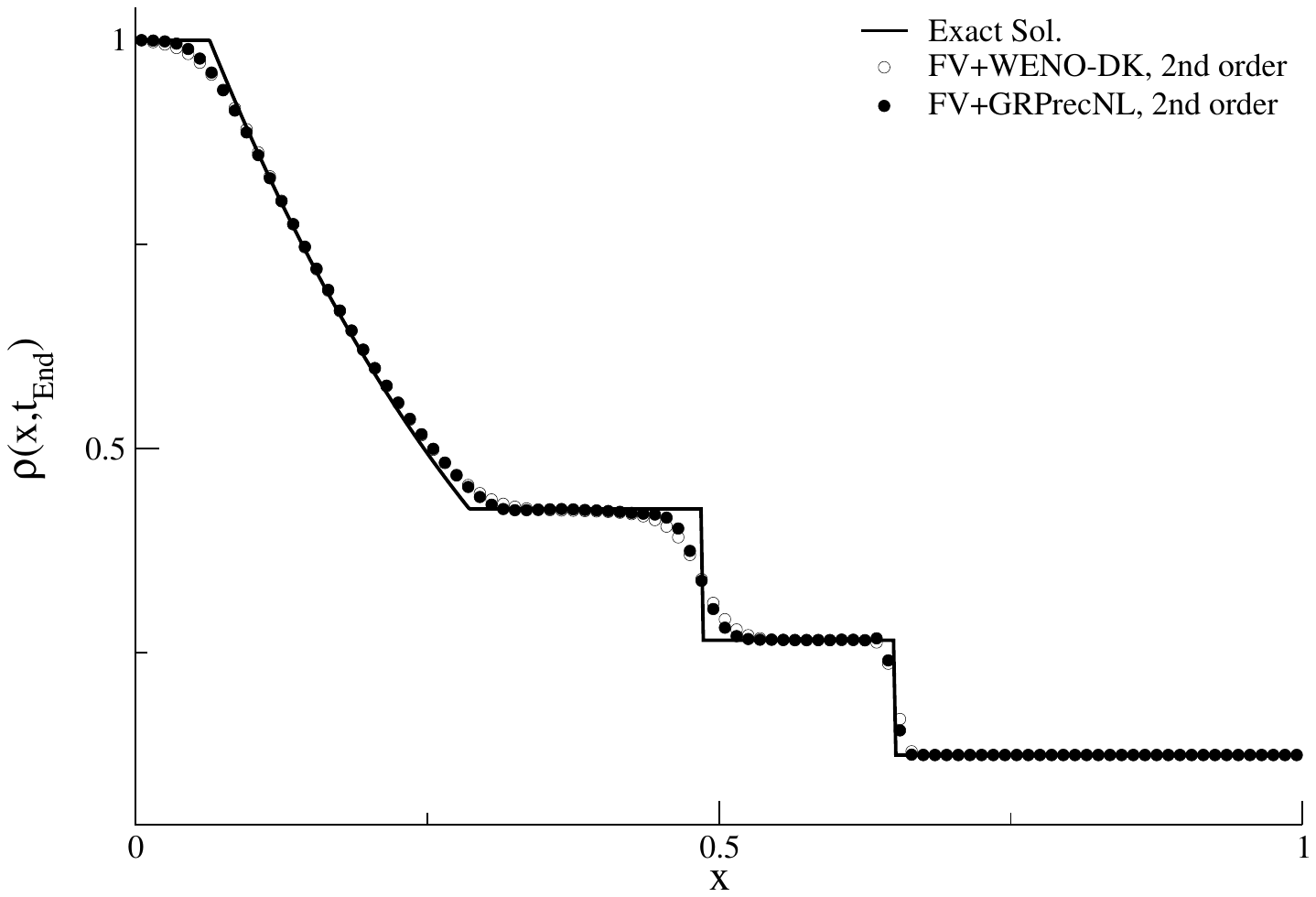}}
\subfloat[][3rd order]{ \includegraphics[scale=0.22]{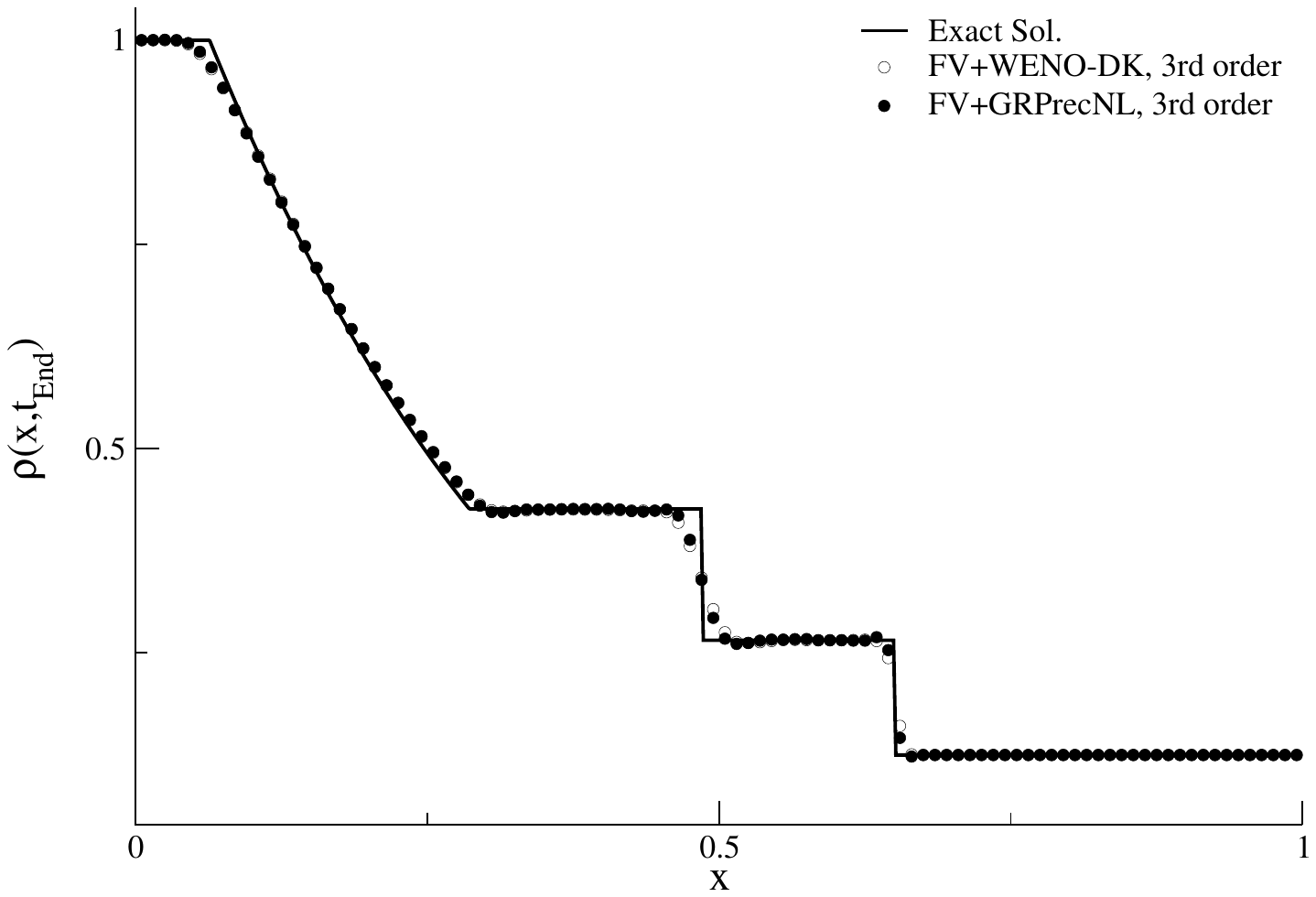}} 
\quad
\subfloat[][4th order]{\includegraphics[scale=0.22]{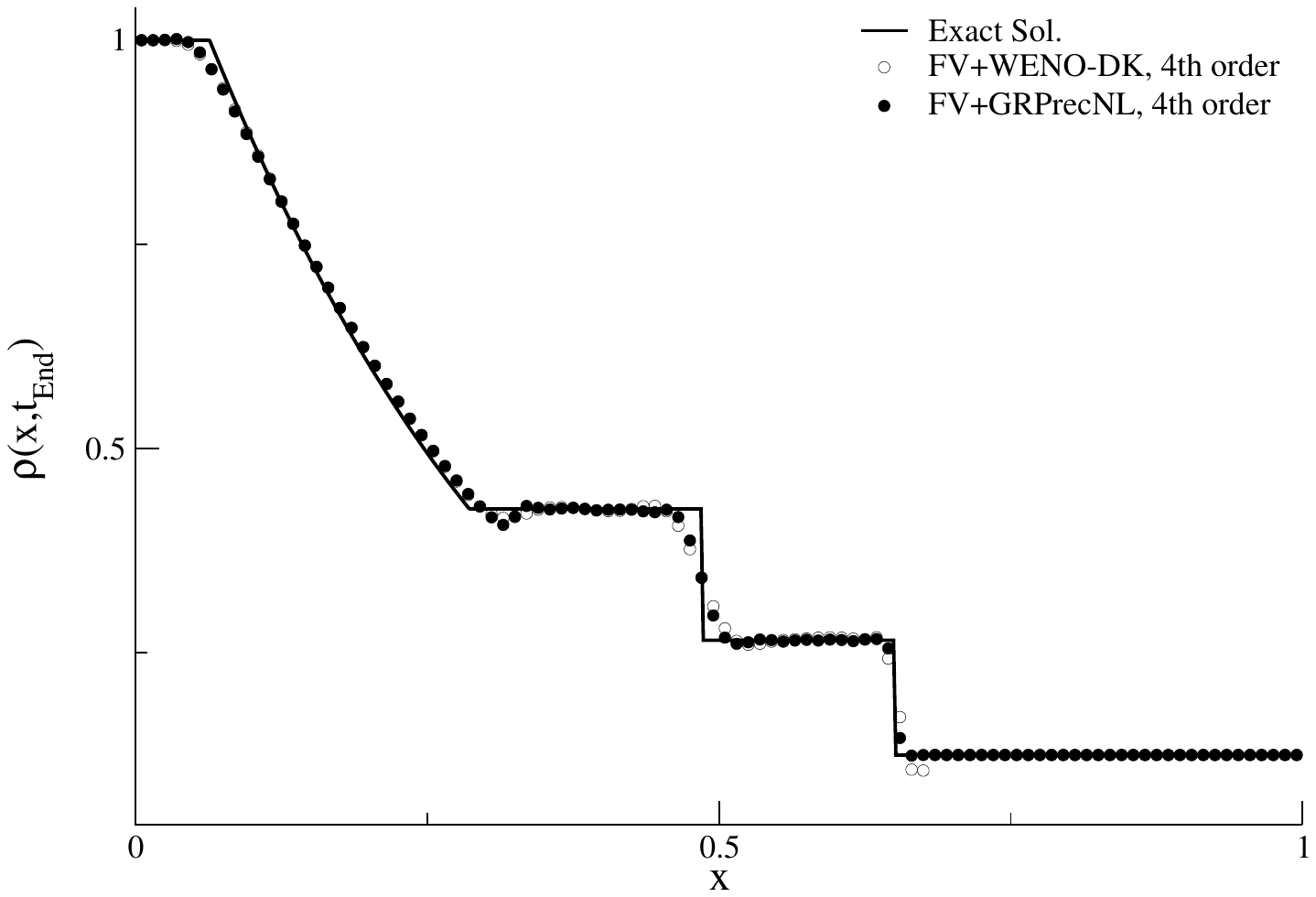}}
\subfloat[][5th order]{ \includegraphics[scale=0.22]{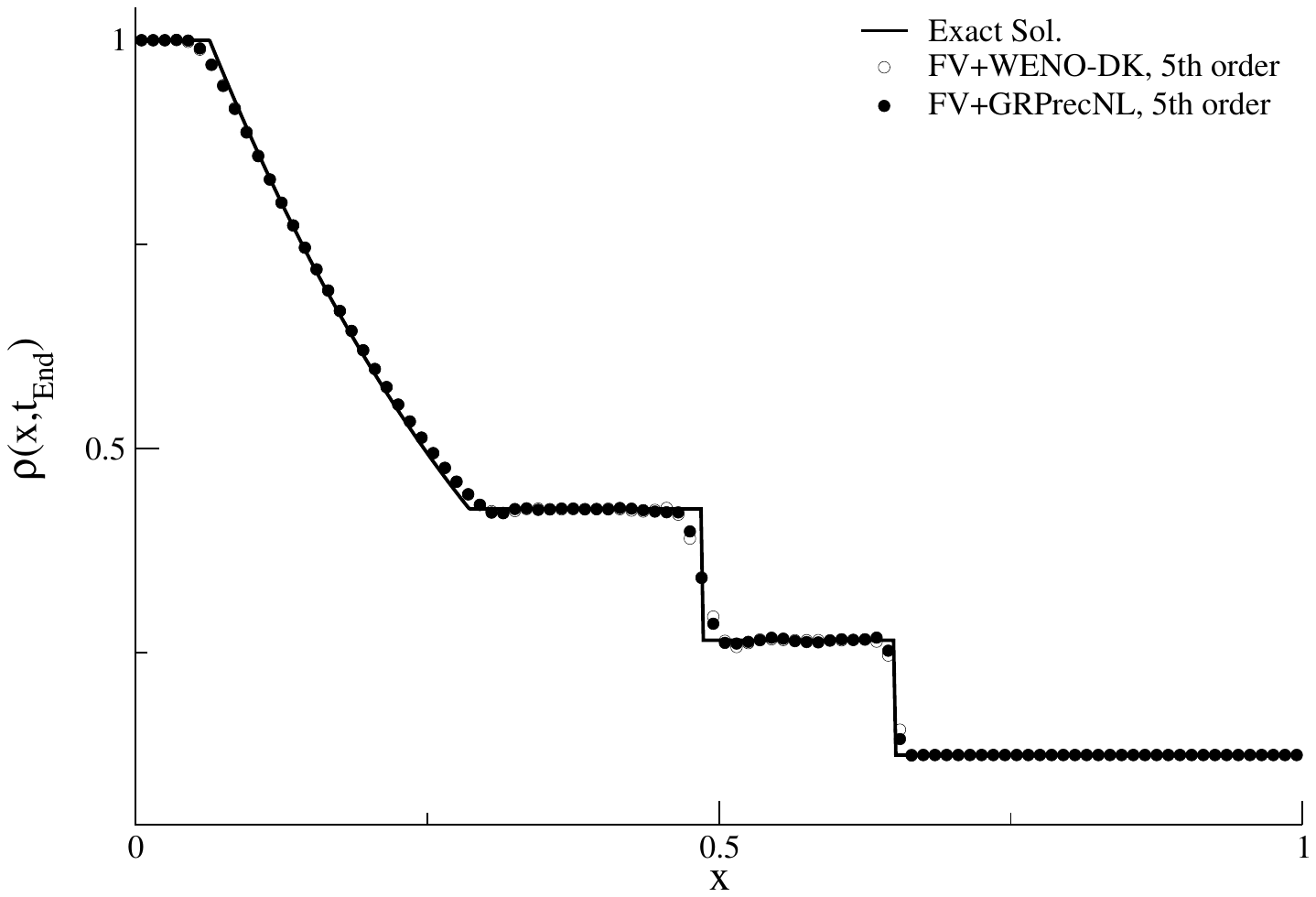}} 

\caption{Euler-Sod problem: $L_1-$error vs CPU time for second, third, fourth and fifth order of accuracy. Parameters $t_{End} = 12$, $C_{CFL}=0.9$.}
\label{figure:EULER-Sod-problem}
\end{figure}

Fig. \ref{figure:EULER-Sod-problem} shows the density for the Sod problem at $t_{out}=0.2s$ obtained with FV+WENO-DK and FV+GRPrecNL reconstructions. The solution of this problem contains three waves:  a left rarefaction, an intermediate contact discontinuity and a right shock.  We observe that FV schemes using WENO-DK and GRPrecNL generate very similar results, with $C_{CFL}=0.9$.

\begin{figure}
\centering
\subfloat[][2nd order]{ \includegraphics[scale=0.23]{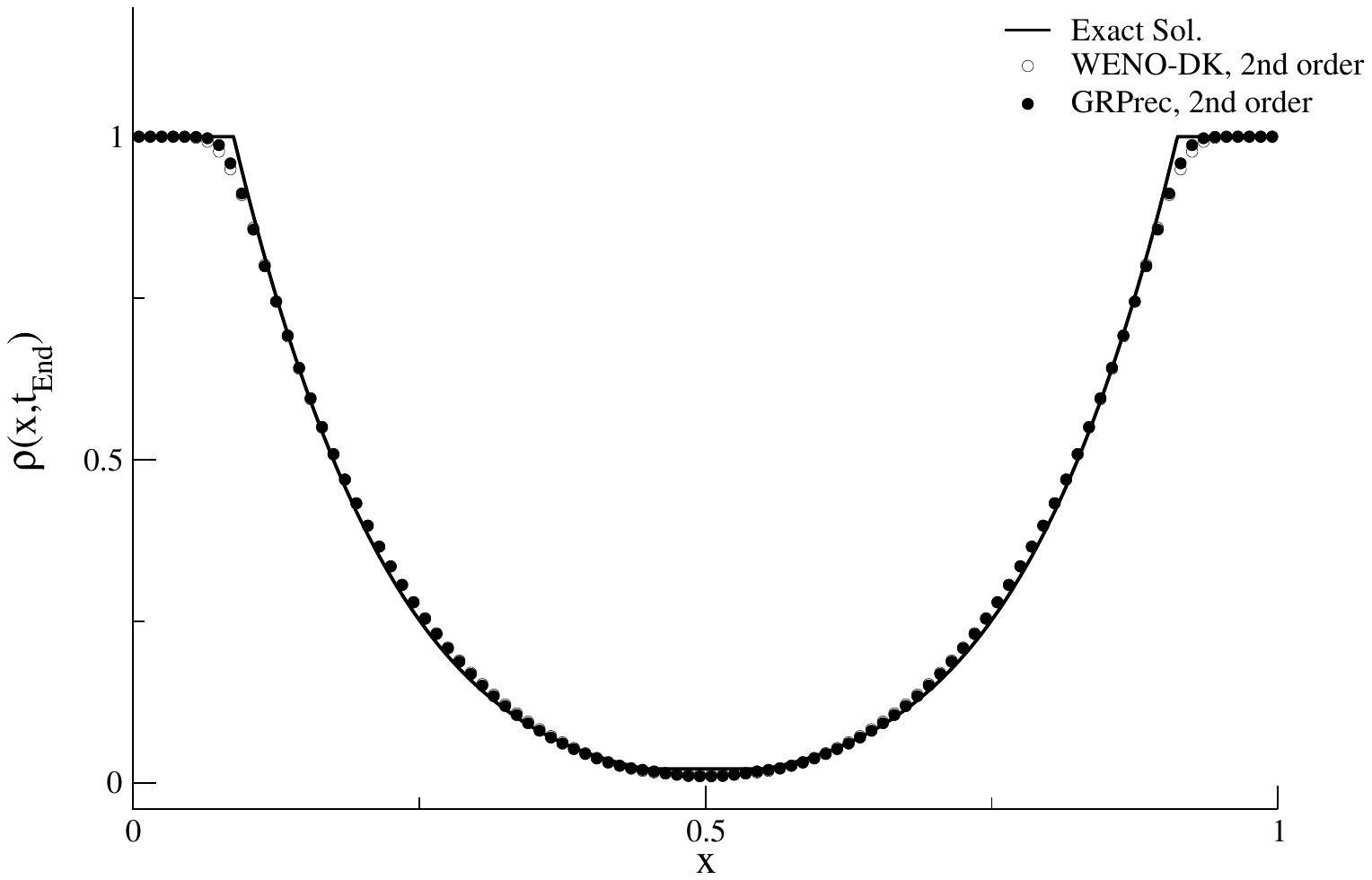}}
\subfloat[][3rd order]{ \includegraphics[scale=0.23]{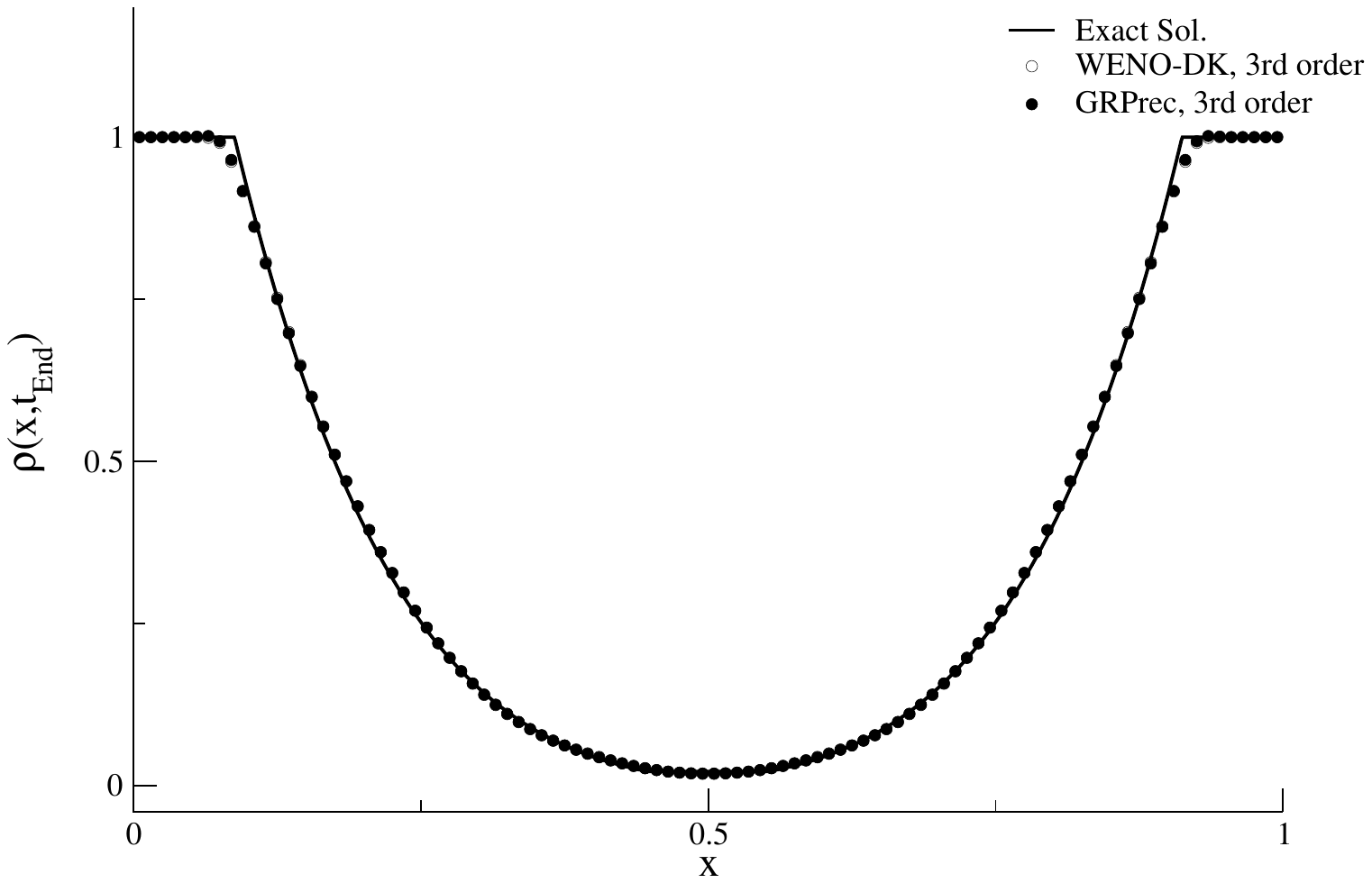}} 
\quad
\subfloat[][4th order]{\includegraphics[scale=0.23]{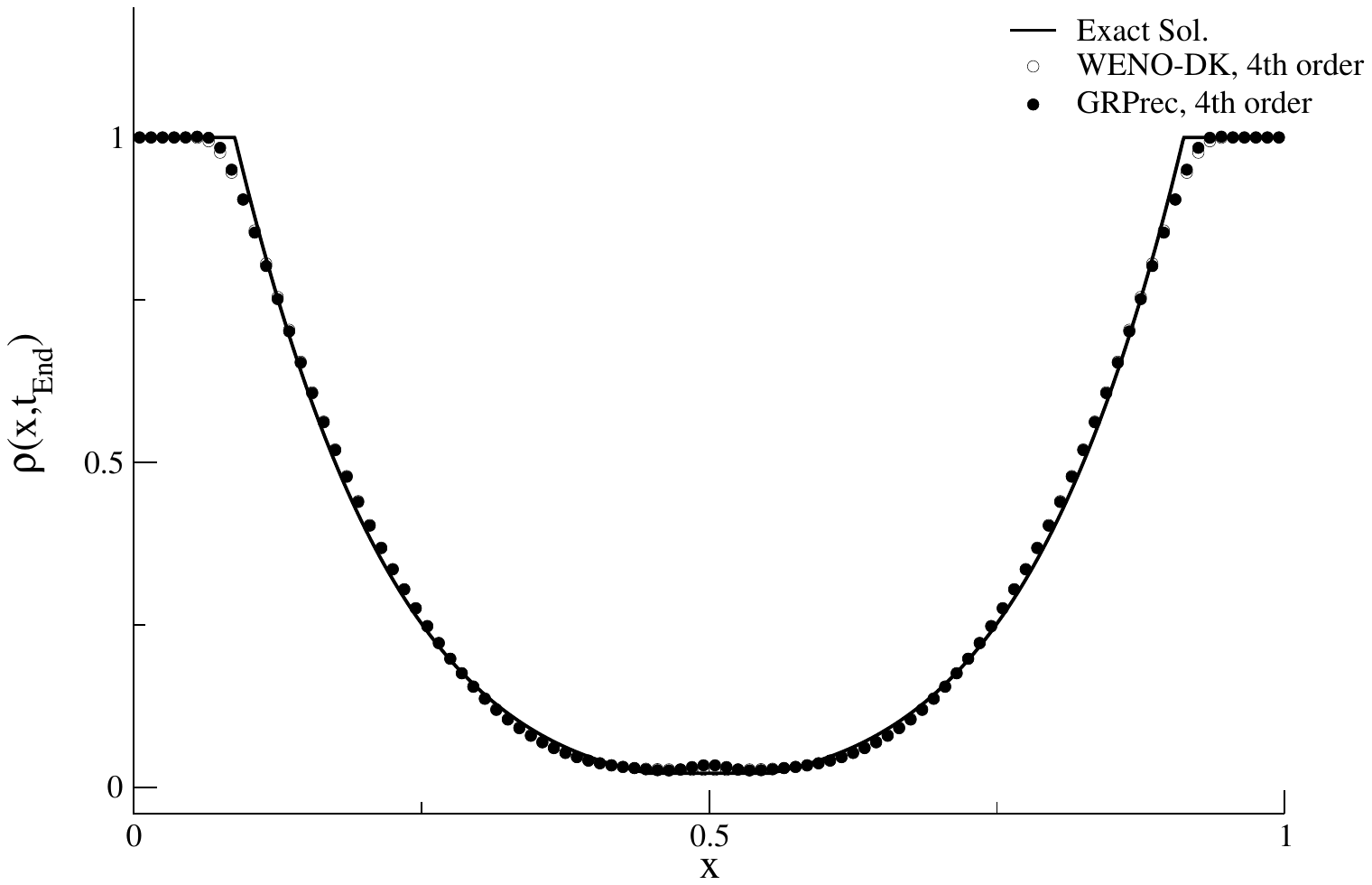}}
\subfloat[][5th order]{ \includegraphics[scale=0.23]{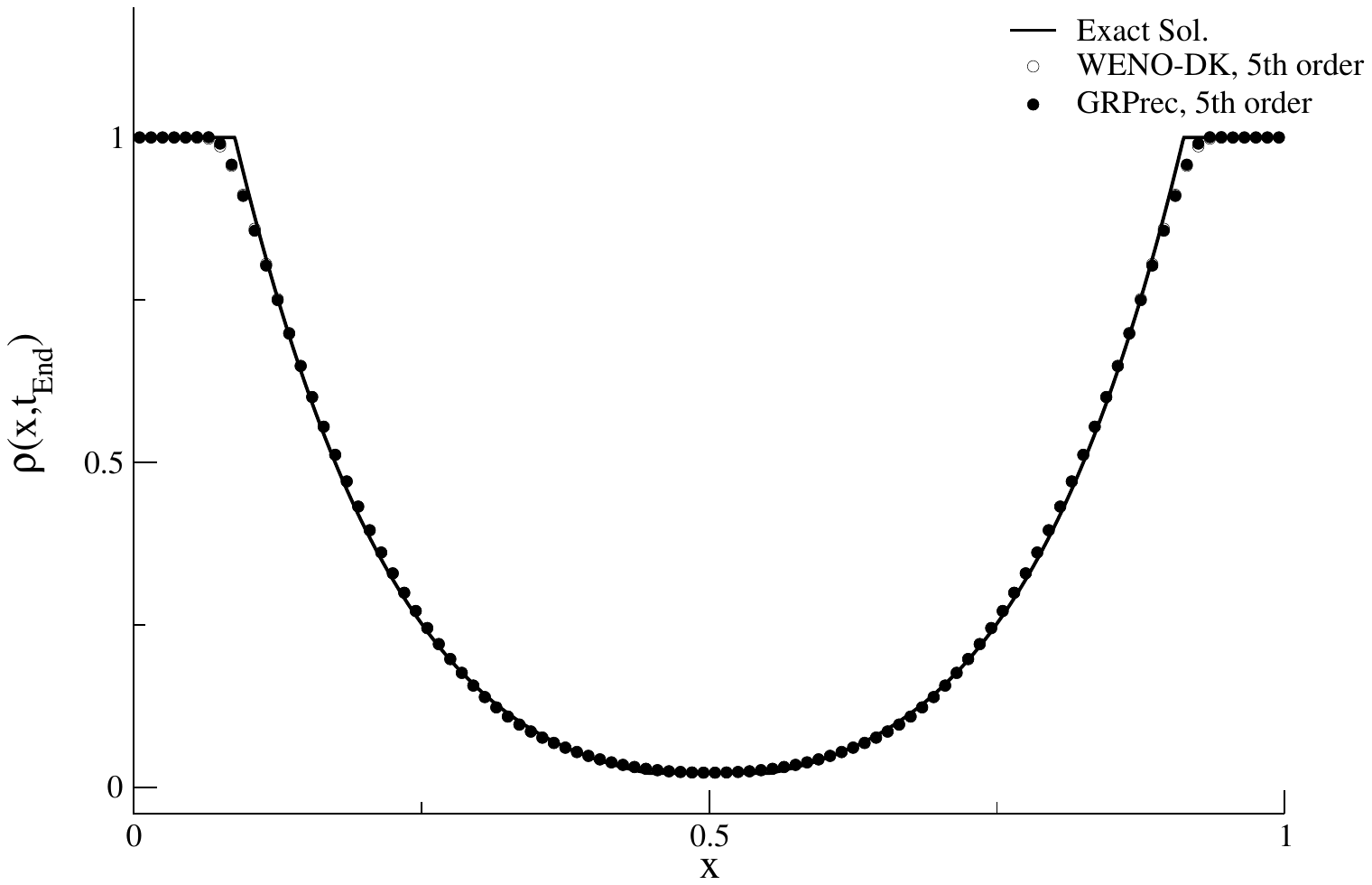}} 

\caption{Euler-123 problem: $L_1-$error vs CPU time for second, third, fourth and fifth orders of accuracy. Parameters $t_{End} = 12$, $C_{CFL}=0.9$ ( $C_{CFL}=0.7$ for the fourth and fifth orders).}
\label{figure:EULER-123-problem}
\end{figure}

Fig. \ref{figure:EULER-123-problem} shows the density for the {\it 123 problem} at $t_{out}=0.15s$, obtained with FV+WENO-DK and FV+GRPrecNL reconstructions. The solution of this test consists of two strong rarefactions and a trivial stationary contact discontinuity. Notably, pressure and density reach very small values, which represents a numerical challenge.  We observe that the solutions obtained with FV+GRPrecNL is very similar to the one computed with FV+WENO-DK. We also note that for fourth and fifth order of accuracy all  schemes used the reduced CFL coefficient,  namely $C_{CFL}=0.7$; for $C_{CFL}=0.9$ negative densities are generated by the schemes.  


\section{Conclusions}\label{sec:conclusions}

We have presented a new spatial polynomial reconstruction procedure,  called GRPrec,  suitable for high-order finite volume methods in both the semi-discrete and the fully discrete frameworks.
The novelty of the reconstruction scheme consists in adding new data at cell interfaces to the classical cell averages, resulting in very compact stencils for prescribed orders of accuracy.  For example,  for schemes of fifth order of accuracy  the stencil of the corresponding fourth degree polynomial involves the cell of interest and its immediate neighbours.  This property makes GRPrec  analogous to the classical DG scheme. However, finite volume schemes based on GRPrec  have the much more generous stability limit of unity for all orders, instead of 
$\frac{1}{(2m+ 1)}$,  which depends on the order of accuracy $m+1$,  or the degree $m$ of the underlying polynomial.


Our reconstruction scheme GRPrec has been combined with a fully discrete  high-order ADER finite volume scheme, whose building block is the solution of the generalized Riemann problem GRP$_m$.  In this paper we have used the solver DET to compute an approximate solution to GRP$_m$ at the interface in order to determine the numerical flux, noting that DET makes it possible to reconcile high order of accuracy with stiffness of source terms, if present.
Fully discrete one step schemes of  up to fifth order of accuracy in space and time have been implemented for both 
GRPrec  and its non-linear version GRPrecNL,  for solving the linear advection equation and the Euler equations of gas dynamics.   

For smooth solutions, convergence rate studies have shown that the theoretically expected orders of accuracy are attained by the new schemes,  while for solutions including shock waves the schemes are seen to be sufficiently robust and to cope satisfactorily with discontinuities.  Computed results have been compared with exact solutions and with numerical solutions obtained from two established high-order ADER methods, one in the finite volume framework and another in the DG framework.   Regarding efficiency,  that is error against CPU cost,  GRPrec has been shown to be
the most efficient of all schemes tested.

The design of  non-linear versions of the reconstruction scheme  remains a challenge.  Surprisingly,   results from our preliminary version GRPrecNL applied to the Euler equations are quite satisfactory,  but not so much for the linear advection equation.  This aspect of the methodology will be the subject of future work,  which may also consider the alternative MOOD approach  \cite{Diot:2013a} as a distinct possibility, as well as limiting techniques applied to the conventional semi-discrete DG schemes  \cite{Wei:2024a}.The inclusion of additional interface knots in the reconstruction procedure, a key contribution of this paper,  opens a wide range of possibilities, also for multidimensional meshes, 
structured or unstructured.
The current version of GRPrec uses approximate solutions of generalized Riemann problems
GRP$_m$  at the previous time level.   Simplifications of this aspect of the methodology are desirable,  as are
modifications of the approach to produce well-balanced versions of the schemes.

Here we have limited ourselves to implementations of GRPrec to high-order finite volume schemes in the fully-discrete,  one-step ADER framework; in future works we shall implement and asses GRPrec for high-order finite volume schemes in the framework of semi-discrete,  multistep  methods,  such as in  \cite{Chu:2025a} and \cite{Mecalizzi:2025a}.

\section*{Acknowledgments}

Lucas O. Müller is a member of the ”Gruppo Nazionale per il Calcolo Scientifico dell’ Istituto Nazionale di Alta Matematica” (INdAM-GNCS, Italy). He also acknowledges funding by the European Union under NextGenerationEU, Mission 4, Component 2 - PRIN 2022 (D.D. 104/22), project title: Immersed methods for multIscale and multiphysics problems, CUP: E53D23005920006.

\bibliographystyle{plain}
\bibliography{ref}

\appendix

\section{Empirical convergence rates for the linear advection equation} \label{sec:app_conv_lae}

In this section we report empirical convergence rates for the test presented in Sec. \ref{sec:conv_lae}, for FV+WENO-DK and DG schemes.


\begin{table}
\begin{center}
Theoretical order : 2, WENO-DK \\
\begin{tabular}{cccccccc} 

\hline
\hline 
Mesh  & $L_\infty$ - ord  & $L_\infty$ - err  & $L_1$ - ord  & $L_1$ - err  & $L_2$ - ord  & $L_2$ - err  & CPU  \\  
\hline

    16  &  -  &$  1.95e-0 1$&  - &$  1.88e-0 1$&  -  &$  1.57e-0 1$  &0.0004\\
    32  &  0.90  &$  1.04e-0 1$&  1.72  &$  5.70e-0 2$&  1.37  &$  6.05e-0 2$  &0.0014\\
    64  &  1.18  &$  4.61e-0 2$&  1.53  &$  1.98e-0 2$&  1.49  &$  2.15e-0 2$  &0.0051\\
   128  &  1.35  &$  1.81e-0 2$&  1.93  &$  5.19e-0 3$&  1.68  &$  6.72e-0 3$  &0.0198\\
   256  &  1.46  &$  6.55e-0 3$&  2.07  &$  1.24e-0 3$&  1.83  &$  1.88e-0 3$  &0.0801\\




 \hline
 \\
\end{tabular}
\\
Theoretical order : 3, WENO-DK \\
\begin{tabular}{cccccccc} 

\hline
\hline 
Mesh  & $L_\infty$ - ord  & $L_\infty$ - err  & $L_1$ - ord  & $L_1$ - err  & $L_2$ - ord  & $L_2$ - err  & CPU  \\  
\hline

    16  &  -  &$  1.27e-0 1$&  -  &$  1.25e-0 1$&  -  &$  1.04e-0 1$  &0.0005\\
    32  &  2.15  &$  2.86e-0 2$&  2.31  &$  2.52e-0 2$&  2.29  &$  2.12e-0 2$  &0.0017\\
    64  &  2.82  &$  4.07e-0 3$&  2.70  &$  3.87e-0 3$&  2.73  &$  3.19e-0 3$  &0.0064\\
   128  &  2.97  &$  5.19e-0 4$&  2.84  &$  5.39e-0 4$&  2.89  &$  4.30e-0 4$  &0.0264\\
   256  &  2.98  &$  6.57e-0 5$&  3.01  &$  6.70e-0 5$&  2.99  &$  5.40e-0 5$  &0.1022\\




 \hline
 \\
\end{tabular}
\\
Theoretical order : 4, WENO-DK\\
\begin{tabular}{cccccccc} 

\hline
\hline 
Mesh  & $L_\infty$ - ord  & $L_\infty$ - err  & $L_1$ - ord  & $L_1$ - err  & $L_2$ - ord  & $L_2$ - err  & CPU  \\  
\hline

    16  &  -  &$  9.00e-0 2$& -  &$  9.22e-0 2$&  -  &$  7.59e-0 2$  &0.0006\\
    32  &  3.50  &$  7.95e-0 3$&  4.14  &$  5.23e-0 3$&  4.00  &$  4.75e-0 3$  &0.0020\\
    64  &  3.86  &$  5.46e-0 4$&  4.17  &$  2.90e-0 4$&  4.18  &$  2.62e-0 4$  &0.0079\\
   128  &  2.73  &$  8.26e-0 5$&  3.85  &$  2.02e-0 5$&  3.32  &$  2.62e-0 5$  &0.0317\\
   256  &  3.03  &$  1.01e-0 5$&  3.86  &$  1.39e-0 6$&  3.46  &$  2.38e-0 6$  &0.1279\\




 \hline
 \\
\end{tabular}
\\
Theoretical order : 5, WENO-DK \\
\begin{tabular}{cccccccc} 

\hline
\hline 
Mesh  & $L_\infty$ - ord  & $L_\infty$ - err  & $L_1$ - ord  & $L_1$ - err  & $L_2$ - ord  & $L_2$ - err  & CPU  \\  
\hline

    16  &  0.00  &$  5.19e-0 2$&  0.00  &$  5.70e-0 2$&  0.00  &$  4.41e-0 2$  &0.0010\\
    32  &  2.49  &$  9.22e-0 3$&  3.38  &$  5.46e-0 3$&  3.15  &$  4.97e-0 3$  &0.0039\\
    64  &  6.42  &$  1.08e-0 4$&  5.40  &$  1.29e-0 4$&  5.61  &$  1.02e-0 4$  &0.0139\\
   128  &  4.97  &$  3.45e-0 6$&  4.97  &$  4.13e-0 6$&  4.97  &$  3.25e-0 6$  &0.0531\\
   256  &  5.00  &$  1.08e-0 7$&  5.00  &$  1.29e-0 7$&  5.00  &$  1.01e-0 7$  &0.2126\\




 \hline

\end{tabular}
\end{center}
\caption{Linear advection - quartic sinus wave.  Solution obtained with FV+WENO-DK for second, third, fourth and fifth order of accuracy. Output time  $t_{out} = 4$, with  $C_{CFL}= 0.9$.}

  \label{Table-LAE-Quartic-WENO-DK}
\end{table}
%


\begin{table}
\begin{center}
Theoretical order : 2, DG \\
\begin{tabular}{cccccccc} 

\hline
\hline 
Mesh  & $L_\infty$ - ord  & $L_\infty$ - err  & $L_1$ - ord  & $L_1$ - err  & $L_2$ - ord  & $L_2$ - err  & CPU  \\  
\hline

    16  &  -  &$  4.93e-0 1$&  -  &$  5.70e-0 1$&  -   &$  4.51e-0 1$  &0.0005\\
    32  &  0.76  &$  2.91e-0 1$&  1.17  &$  2.54e-0 1$&  1.00  &$  2.26e-0 1$  &0.0018\\
    64  &  1.30  &$  1.18e-0 1$&  1.18  &$  1.12e-0 1$&  1.27  &$  9.34e-0 2$  &0.0070\\
   128  &  1.87  &$  3.22e-0 2$&  1.82  &$  3.16e-0 2$&  1.84  &$  2.62e-0 2$  &0.0274\\
   256  &  1.98  &$  8.19e-0 3$&  1.97  &$  8.07e-0 3$&  1.97  &$  6.67e-0 3$  &0.1041\\



 \hline
 \\
\end{tabular}
\\
Theoretical order : 3, DG \\
\begin{tabular}{cccccccc} 

\hline
\hline 
Mesh  & $L_\infty$ - ord  & $L_\infty$ - err  & $L_1$ - ord  & $L_1$ - err  & $L_2$ - ord  & $L_2$ - err  & CPU  \\  
\hline

    16  &  -  &$  5.43e-0 2$&  -  &$  7.17e-0 2$&  -  &$  5.38e-0 2$  &0.0017\\
    32  &  2.36  &$  1.05e-0 2$&  2.64  &$  1.15e-0 2$&  2.56  &$  9.12e-0 3$  &0.0065\\
    64  &  2.88  &$  1.43e-0 3$&  2.95  &$  1.49e-0 3$&  2.93  &$  1.20e-0 3$  &0.0248\\
   128  &  2.98  &$  1.82e-0 4$&  2.98  &$  1.89e-0 4$&  2.99  &$  1.51e-0 4$  &0.0996\\
   256  &  3.00  &$  2.28e-0 5$&  3.00  &$  2.37e-0 5$&  3.00  &$  1.89e-0 5$  &0.4008\\



 \hline
 \\
\end{tabular}
\\
Theoretical order : 4, DG\\
\begin{tabular}{cccccccc} 

\hline
\hline 
Mesh  & $L_\infty$ - ord  & $L_\infty$ - err  & $L_1$ - ord  & $L_1$ - err  & $L_2$ - ord  & $L_2$ - err  & CPU  \\  
\hline

    16  & -  &$  5.50e-0 3$&  -  &$  5.71e-0 3$& -  &$  5.10e-0 3$  &0.0044\\
    32  &  3.98  &$  3.48e-0 4$&  3.83  &$  4.03e-0 4$&  3.98  &$  3.23e-0 4$  &0.0172\\
    64  &  4.01  &$  2.17e-0 5$&  3.99  &$  2.54e-0 5$&  4.01  &$  2.00e-0 5$  &0.0698\\
   128  &  4.00  &$  1.35e-0 6$&  4.00  &$  1.59e-0 6$&  4.00  &$  1.25e-0 6$  &0.2799\\
   256  &  4.00  &$  8.45e-0 8$&  4.00  &$  9.92e-0 8$&  4.00  &$  7.82e-0 8$  &1.1235\\

  

 \hline
 \\
\end{tabular}
\\
Theoretical order : 5, DG \\
\begin{tabular}{cccccccc} 

\hline
\hline 
Mesh  & $L_\infty$ - ord  & $L_\infty$ - err  & $L_1$ - ord  & $L_1$ - err  & $L_2$ - ord  & $L_2$ - err  & CPU  \\  
\hline

    16  &  -  &$  2.71e-0 4$&  -  &$  3.04e-0 4$&  -  &$  2.59e-0 4$  &0.0116\\
    32  &  4.84  &$  9.48e-0 6$&  4.77  &$  1.11e-0 5$&  4.85  &$  8.96e-0 6$  &0.0493\\
    64  &  4.95  &$  3.06e-0 7$&  4.93  &$  3.65e-0 7$&  4.96  &$  2.89e-0 7$  &0.1905\\
   128  &  4.99  &$  9.65e-0 9$&  4.98  &$  1.16e-0 8$&  4.99  &$  9.10e-0 9$  &0.7531\\
   256  &  5.00  &$  3.03e-010$&  5.00  &$  3.63e-010$&  5.00  &$  2.85e-010$  &2.9046\\



 \hline

\end{tabular}
\end{center}
\caption{Linear advection - quartic sinus wave.  Solution obtained with DG for second, third, fourth and fifth order of accuracy. Output time  $t_{out} = 4$, with  $C_{CFL}= 0.9$.}

  \label{Table-LAE-Quartic-ADER-DG}
\end{table}

\section{Empirical convergence rates for the Euler equations} \label{sec:app_conv_euler}

In this section we report empirical convergence rates for the test presented in Sec. \ref{sec:conv_euler}, for FV+WENO-Dk and DG schemes.


\begin{table}
\begin{center}
Theoretical order : 2, WENO-DK \\
\begin{tabular}{cccccccc} 

\hline
\hline 
Mesh  & $L_\infty$ - ord  & $L_\infty$ - err  & $L_1$ - ord  & $L_1$ - err  & $L_2$ - ord  & $L_2$ - err  & CPU  \\  
\hline

    40  &  -  &$  1.89e-0 1$&  - &$  2.43e-0 1$&  -  &$  1.91e-0 1$  &1.0936\\
    80  &  0.70  &$  1.17e-0 1$&  0.88  &$  1.32e-0 1$&  0.82  &$  1.08e-0 1$  &1.2252\\
   160  &  1.39  &$  4.46e-0 2$&  1.72  &$  4.02e-0 2$&  1.66  &$  3.42e-0 2$  &5.0603\\
   320  &  1.31  &$  1.80e-0 2$&  1.45  &$  1.47e-0 2$&  1.39  &$  1.30e-0 2$  &20.5157\\
   640  &  1.37  &$  6.93e-0 3$&  2.00  &$  3.68e-0 3$&  1.73  &$  3.94e-0 3$  &86.7846\\


 \hline
 \\
\end{tabular}
\\
Theoretical order : 3, WENO-DK \\
\begin{tabular}{cccccccc} 

\hline
\hline 
Mesh  & $L_\infty$ - ord  & $L_\infty$ - err  & $L_1$ - ord  & $L_1$ - err  & $L_2$ - ord  & $L_2$ - err  & CPU  \\  
\hline

    40  &  - &$  4.71e-0 2$&  -  &$  6.08e-0 2$&  -  &$  4.76e-0 2$  &1.1132\\
    80  &  2.79  &$  6.82e-0 3$&  2.80  &$  8.72e-0 3$&  2.80  &$  6.84e-0 3$  &1.8443\\
   160  &  2.97  &$  8.72e-0 4$&  2.97  &$  1.11e-0 3$&  2.97  &$  8.72e-0 4$  &7.3015\\
   320  &  3.00  &$  1.09e-0 4$&  3.00  &$  1.39e-0 4$&  3.00  &$  1.09e-0 4$  &27.0575\\
   640  &  3.00  &$  1.37e-0 5$&  3.00  &$  1.74e-0 5$&  3.00  &$  1.37e-0 5$  &104.8304\\


 \hline
 \\
\end{tabular}
\\
Theoretical order : 4, WENO-DK\\
\begin{tabular}{cccccccc} 

\hline
\hline 
Mesh  & $L_\infty$ - ord  & $L_\infty$ - err  & $L_1$ - ord  & $L_1$ - err  & $L_2$ - ord  & $L_2$ - err  & CPU  \\  
\hline

    40  &  -  &$  8.59e-0 3$&  -  &$  1.11e-0 2$&  -  &$  8.67e-0 3$  &1.4206\\
    80  &  4.64  &$  3.44e-0 4$&  4.66  &$  4.39e-0 4$&  4.65  &$  3.45e-0 4$  &2.4708\\
   160  &  4.41  &$  1.62e-0 5$&  4.41  &$  2.07e-0 5$&  4.41  &$  1.62e-0 5$  &9.4315\\
   320  &  4.16  &$  9.09e-0 7$&  4.16  &$  1.16e-0 6$&  4.16  &$  9.09e-0 7$  &37.1746\\
   640  &  4.05  &$  5.51e-0 8$&  4.05  &$  7.01e-0 8$&  4.05  &$  5.51e-0 8$  &155.4573\\


 \hline
 \\
\end{tabular}
\\
Theoretical order : 5, WENO-DK \\
\begin{tabular}{cccccccc} 

\hline
\hline 
Mesh  & $L_\infty$ - ord  & $L_\infty$ - err  & $L_1$ - ord  & $L_1$ - err  & $L_2$ - ord  & $L_2$ - err  & CPU  \\  
\hline

    40  &  -  &$  1.43e-0 3$&  -  &$  1.85e-0 3$&  - &$  1.44e-0 3$  &2.5085\\
    80  &  4.95  &$  4.61e-0 5$&  4.97  &$  5.90e-0 5$&  4.96  &$  4.63e-0 5$  &6.9937\\
   160  &  4.99  &$  1.45e-0 6$&  4.99  &$  1.85e-0 6$&  4.99  &$  1.45e-0 6$  &25.7010\\
   320  &  5.00  &$  4.55e-0 8$&  5.00  &$  5.79e-0 8$&  5.00  &$  4.55e-0 8$  &103.7334\\
   640  &  5.00  &$  1.42e-0 9$&  5.00  &$  1.81e-0 9$&  5.00  &$  1.42e-0 9$  &401.6974\\


 \hline

\end{tabular}
\end{center}
\caption{Euler equations.  Solution obtained with FV+WENO-DK for second, third, fourth and fifth order of accuracy. Output  time  $t_{out} = 4$, with  $C_{CFL}= 0.9$ ( $C_{CFL}= 0.7$ for fifth order only).}

  \label{Table-Euler-WENO-DK}
\end{table}
%


\begin{table}
\begin{center}
Theoretical order : 2, DG \\
\begin{tabular}{cccccccc} 

\hline
\hline 
Mesh  & $L_\infty$ - ord  & $L_\infty$ - err  & $L_1$ - ord  & $L_1$ - err  & $L_2$ - ord  & $L_2$ - err  & CPU  \\  
\hline

     40  &  -  &$  3.16e-0 3$&  -  &$  3.71e-0 3$&  -  &$  2.92e-0 3$  &2.3277\\
    80  &  2.26  &$  6.61e-0 4$&  2.21  &$  8.03e-0 4$&  2.20  &$  6.34e-0 4$  &6.0261\\
   160  &  2.08  &$  1.57e-0 4$&  2.05  &$  1.94e-0 4$&  2.05  &$  1.53e-0 4$  &24.4942\\
   320  &  2.02  &$  3.85e-0 5$&  2.01  &$  4.80e-0 5$&  2.01  &$  3.79e-0 5$  & 103.3990\\
   640  &  2.01  &$  9.59e-0 6$&  2.00  &$  1.20e-0 5$&  2.00  &$  9.45e-0 6$  &416.2141\\


 \hline
 \\
\end{tabular}
\\
Theoretical order : 3, DG \\
\begin{tabular}{cccccccc} 

\hline
\hline 
Mesh  & $L_\infty$ - ord  & $L_\infty$ - err  & $L_1$ - ord  & $L_1$ - err  & $L_2$ - ord  & $L_2$ - err  & CPU  \\  
\hline

    40  &  -  &$  1.53e-0 4$&  -  &$  1.92e-0 4$& -  &$  1.52e-0 4$  &3.9881\\
    80  &  2.99  &$  1.93e-0 5$&  2.99  &$  2.43e-0 5$&  2.99  &$  1.91e-0 5$  &12.47077\\
   160  &  3.00  &$  2.41e-0 6$&  2.99  &$  3.05e-0 6$&  3.00  &$  2.39e-0 6$  &55.7678\\
   320  &  3.00  &$  3.02e-0 7$&  3.00  &$  3.82e-0 7$&  3.00  &$  3.00e-0 7$  &219.3321\\
   640  &  3.00  &$  3.78e-0 8$&  3.00  &$  4.77e-0 8$&  3.00  &$  3.75e-0 8$  &885.8146\\


 \hline
 \\
\end{tabular}
\\
Theoretical order : 4, DG\\
\begin{tabular}{cccccccc} 

\hline
\hline 
Mesh  & $L_\infty$ - ord  & $L_\infty$ - err  & $L_1$ - ord  & $L_1$ - err  & $L_2$ - ord  & $L_2$ - err  & CPU  \\  
\hline

    40  &  -  &$  8.04e-0 7$&  -  &$  1.00e-0 6$&  -  &$  7.91e-0 7$  &7.3087\\
    80  &  4.01  &$  4.99e-0 8$&  4.01  &$  6.21e-0 8$&  4.01  &$  4.91e-0 8$  &27.1410\\
   160  &  4.00  &$  3.11e-0 9$&  4.00  &$  3.88e-0 9$&  4.00  &$  3.06e-0 9$  &109.0513\\
   320  &  4.00  &$  1.95e-010$&  4.00  &$  2.42e-010$&  4.00  &$  1.91e-010$  &461.3226\\
   640  &  4.00  &$  1.22e-011$&  4.00  &$  1.51e-011$&  4.00  &$  1.20e-011$  &2100.1705\\
  

 \hline
 \\
\end{tabular}
\\
Theoretical order : 5, DG \\
\begin{tabular}{cccccccc} 

\hline
\hline 
Mesh  & $L_\infty$ - ord  & $L_\infty$ - err  & $L_1$ - ord  & $L_1$ - err  & $L_2$ - ord  & $L_2$ - err  & CPU  \\  
\hline

    40  &  -  &$  2.51e-0 8$&  -  &$  3.12e-0 8$&  -  &$  2.46e-0 8$  &33.5383\\
    80  &  4.90  &$  8.39e-010$&  4.89  &$  1.05e-0 9$&  4.90  &$  8.25e-010$  &127.3769\\
   160  &  4.97  &$  2.67e-011$&  4.97  &$  3.35e-011$&  4.97  &$  2.63e-011$  &490.2949\\
   320  &  4.94  &$  8.72e-013$&  4.99  &$  1.06e-012$&  4.99  &$  8.29e-013$  &1790.9428\\
   640  &  2.61  &$  1.43e-013$&  3.85  &$  7.31e-014$&  3.73  &$  6.26e-014$  &6413.4435\\


 \hline

\end{tabular}
\end{center}
\caption{Euler equations.  Solution obtained with DG for second, third, fourth and fifth order of accuracy. Output  time  $t_{out} = 4$, with  $C_{CFL}= 0.9$ ( $C_{CFL}= 0.7$ for fifth order only).}

  \label{Table-Euler-ADER-DG}
\end{table}

\end{document}